\documentclass[times,authoryear]{elsarticle}

\usepackage{jasr}
\usepackage{framed,multirow}
\usepackage{amssymb}
\usepackage{latexsym}
\usepackage{url}
\usepackage{xcolor}
\definecolor{newcolor}{rgb}{.8,.349,.1}
\usepackage[citebordercolor=white]{hyperref}
\usepackage[nolist]{acronym}
\usepackage{amsmath}
\usepackage{subcaption}
\usepackage{textcomp}
\usepackage{graphicx}
\usepackage[version=4]{mhchem}
\usepackage{siunitx}
\usepackage{longtable,tabularx}
\usepackage{lineno}
\usepackage{multirow}
\usepackage{pifont}
\usepackage{comment}
\usepackage{algorithmic}  
\usepackage{algorithm} 
\usepackage{import}
\usepackage{ulem}

\usepackage{booktabs, array, tabularx, enumitem}
\newlist{compactitem}{itemize}{1}
\setlist[compactitem]{nosep, topsep=0pt, left=1em}

\journal{Advances in Space Research}

\begin{document}

\verso{Given-name Surname \textit{etal}}

\begin{frontmatter}

\title{A Stochastic Dynamic Network Model of the Space Environment}

\author[1]{Yirui \snm{Wang}}
\author[1]{Pietro \snm{De Marchi}}
\author[1]{Massimiliano \snm{Vasile}\corref{cor1}}
\cortext[cor1]{Corresponding author}
\ead{Massimiliano.vasile@strath.ac.uk}

\affiliation[1]{organization={Aerospace Centre of Excellence, University of Strathclyde},
                addressline={75 Montrose Street},
                city={Glasgow},
                postcode={G11XJ},
                country={UK}}



\begin{abstract}
This work proposes to model the space environment as a stochastic dynamic network where each node is a group of objects of a given class, or species, and their relationship is represented by stochastic links. 
A set of stochastic dynamic equations, governing the evolution of the network, are derived from the network structure and topology. 
It will be shown that the proposed system of stochastic dynamic equations well reproduces existing results on the evolution of the space environment. 
The analysis of the structure of the network and relationships among node can help to understand which species of objects and orbit regimes are more critical and affect the most the future evolution of the space environment. 
In analogy with ecological networks, we develop a theory of the carrying capacity of space based on the stability of equilibria of the network dynamics.

Some examples are presented starting from the current population of resident objects and different launch traffic forecast models.
It will be shown how the proposed network model can be used to study the effect of the adoption of different policies on the execution of collision avoidance and post-mission disposal manoeuvres. 
\end{abstract}


\end{frontmatter}


\section{Introduction}
The space industry is one of today’s most growing sectors, and as a consequence, the number of resident space objects is continuously rising. According to ESA’s Space Environment Report 2023, at the reference epoch November 1st 2016, the estimated number of debris objects in orbit in the different size ranges is: 34,000 objects greater than 10 cm, 900,000 objects from 1 cm to 10 cm, and 128 million objects from 1 mm to 1 cm. A collision among space objects could generate a cloud of space debris that can trigger a cascade of collisions that could make the use of space dangerous for future missions \citep{1978Collision, mcknight1989collision}. Large constellation is an important factor that affects the evolution of the space debris environment \citep{10121609}, and untracked debris will lead to potentially dangerous on-orbit collisions on a regular basis due to the large number of satellites within mega-constellation orbital shells \citep{2021Satellite}. The growth rate of collisional debris would exceed the natural decay rate in $\approx$ 50 years \citep{2009Space, 2009A}.
Moreover, when the number of objects increases, chain effects of 
collision may be triggered, which is known as Kessler Syndrome \citep{1978Collision}. The Kessler Syndrome refers to a scenario in
which Earth orbits inevitably become so polluted with satellite-related orbital debris that a self-reinforcing collisional cascade, which destroys satellites in orbit and
makes orbital space unusable, is inevitable \citep{2018An}. It is therefore essential to develop the necessary techniques to model, understand and predict the current space environment and its future evolution \citep{2019Space}. An effective Space Traffic Management (STM) is pivotal to expand the population of space objects safely and to reduce the burden
on space operators \citep{murakami2019space}.

There are many existing efforts to model the long-term evolution of the space environment. Space agencies have developed their own models, such as the NASA LEGEND environment model \citep{liou2004legend} and ESA’s DELTA model \citep{virgili2016delta}. Both of these examples have high-fidelity propagators that apply perturbations to all objects in the simulation. Their approaches to modelling collision rates calculate probabilities of collisions within control volumes, which is based on the widely used CUBE method \citep{liou2003new}. This method is also used in MEDEE, developed by CNES \citep{jc2013introducing}, and SOLEM, developed by CNSA \citep{wang2019introduction}. Another well-established tool using a similar approach is SDM \citep{anselmo_pardini_SDM3}. This approach to modelling environment evolution can be computationally expensive as it often requires propagating each object. Other environment models use some simplifying assumptions to reduce their computational requirements. It is common to divide the environment into discrete bins based on orbital parameters and define the evolution based on the statistics of each bin. One earlier example is IDES, which divides objects into bins based on some spatial and physical parameters \citep{walker1999enhancement}. INDEMN makes a further simplification of only considering densities of objects in discrete orbital shells in terms of altitude \citep{lucken2019collision}. These approaches allow the generation of much faster estimates of the space environment evolution. MIT Orbital Capacity Tool (MOCAT) is operated by propagating all the Anthropogenic Space Objects (ASOs) forward in time, where MOCAT-MC \citep{jang2023monte} stands for the model based on Monte-Carlo method, and MOCAT-SSEM \citep{dambrosio2023} stands for the model based on a source-sink model, and MOCAT-ML \citep{rodriguez2024towards} stands for the model based on Machine Learning. 
 
When performing long-term simulations through an environment model, the results are usually presented in terms of the number of objects and the cumulative number of catastrophic collisions. Besides this basic information, several metrics have been proposed to understand the health of the space environment. For example, Environmental Consequences
 of Orbital Breakups (ECOB)  focused on the evolution of the consequences (i.e. the effects) of a fragmentation \citep{LETIZIA20161255, letizia2017extending}. Criticality of Spacecraft Index (CSI) is proposed to measure the environmental impact of large bodies \citep{rossi2015criticality}. Undisposed Mass Per Year is a measure as the amount of mass left in orbit due to the failure of PMD \citep{henning2019impacts}. A maximum orbital capacity calculated by treating the model as a dynamical system and calculating the equilibrium point \citep{d2022capacity, CarryingCapacity2024}. Though these metrics provide valuable insights into the characteristics of the space environment from different scopes, there is not a single commonly-accepted metric to measure the health of the space environment.

This paper proposes to model the space environment as a stochastic dynamic network, where each node is a group of objects, belonging to a given class, or species, and their relationship is represented by stochastic links. In the remainder of the paper, this model is called NESSY, or NEtwork model for Space SustainabilitY. This modelling approach is different from previous works on the modelling of the space environment in that it enables a direct and explicit analysis of the relationships among different species of objects, it allows identifying communities of objects with similar characteristics and their impact on the rest of the environment. It also differs from previous works on the use of network theory applied to the space environment \citep{Romano_2024} in that it derives the evolutionary equations from the network structure. 
The results in this paper build upon previous work by the authors \cite{acciarini2020multi,acciarini2021network,wang2023multi} and introduces five new main contributions: i) an improvement of the completeness of the network model introduced in \cite{acciarini2020multi,acciarini2021network,wang2023multi}, with the ability to model more species of space objects and interactions; ii) a comparison of NESSY against another well-established and validated model; iii) a study of the effect of launch traffic, collision avoidance manoeuvres and post mission disposal policies on the evolution of the space environment; iv) a preliminary study of centrality and eigenvalues of the network dynamics and v) a theory of the carrying capacity of the space environment, starting from the stability of the equilibria of the network dynamics. 

The paper is structured as follows. In Section 2, we introduce the definition of our proposed network model and some validation tests against known evolutionary models. In Section 3, we present some examples where the network model is used to predict the effect of launch traffic, collision avoidance manoeuvres and post-mission disposal policy. In Section 4, we introduce some structural properties of the network like centrality and in Section 5 we introduce a carrying capacity theory. Finally, in Section 6 we conclude with some remarks and recommendations for future work. 

\section{A Network Model of the Space Environment}
\begin{figure*}[h]
    \centering
    \begin{minipage}[t]{0.49\textwidth}
        \centering
        \includegraphics[width=\linewidth]{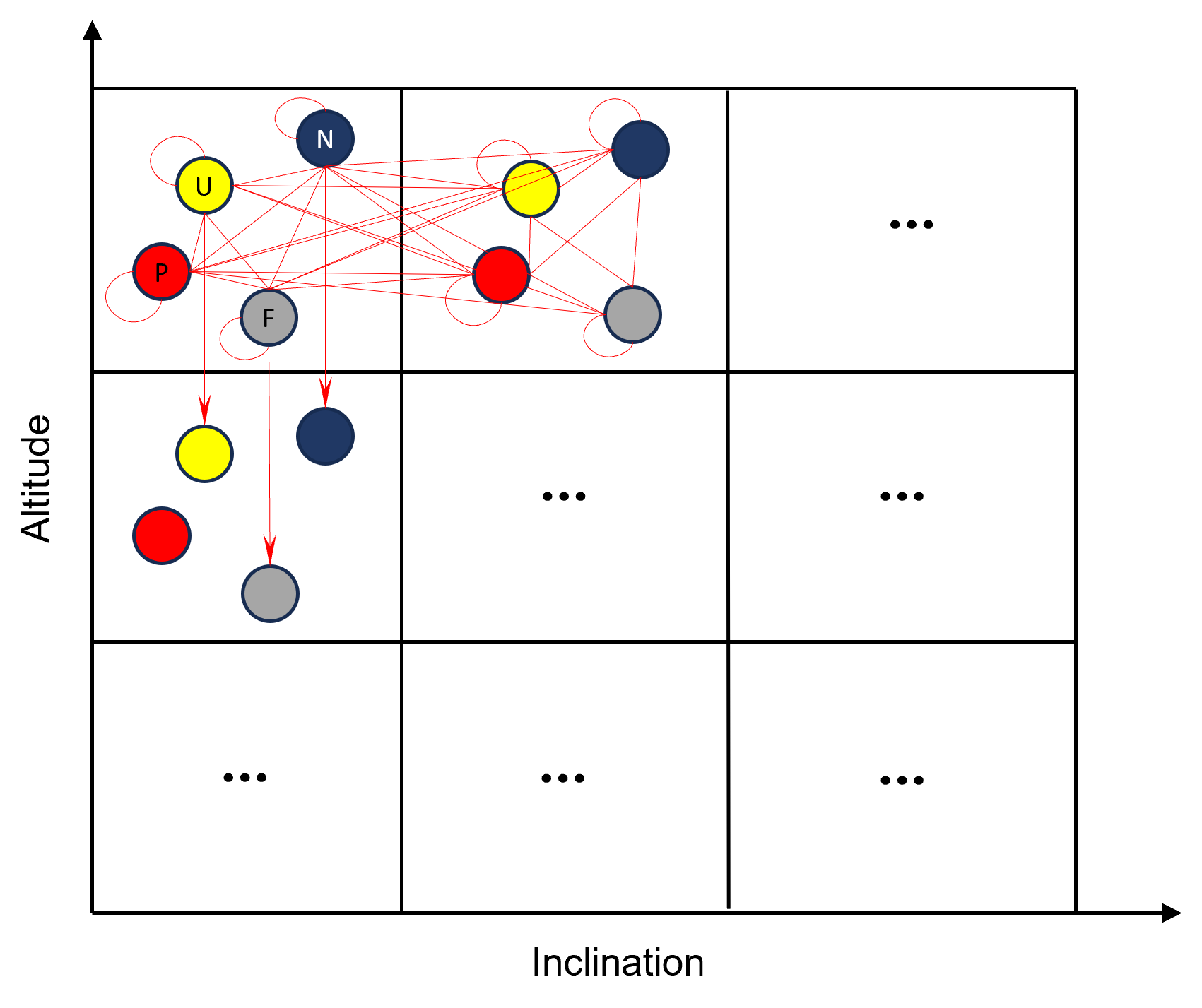}
        \caption{Illustration of the network model of the space environment proposed in this paper.}
        \label{fig:Diagram of network model in the physical layer}
    \end{minipage}
    \hspace{0.04\textwidth}
    \begin{minipage}[t]{0.41\textwidth} 
        \centering
        \includegraphics[width=\linewidth]{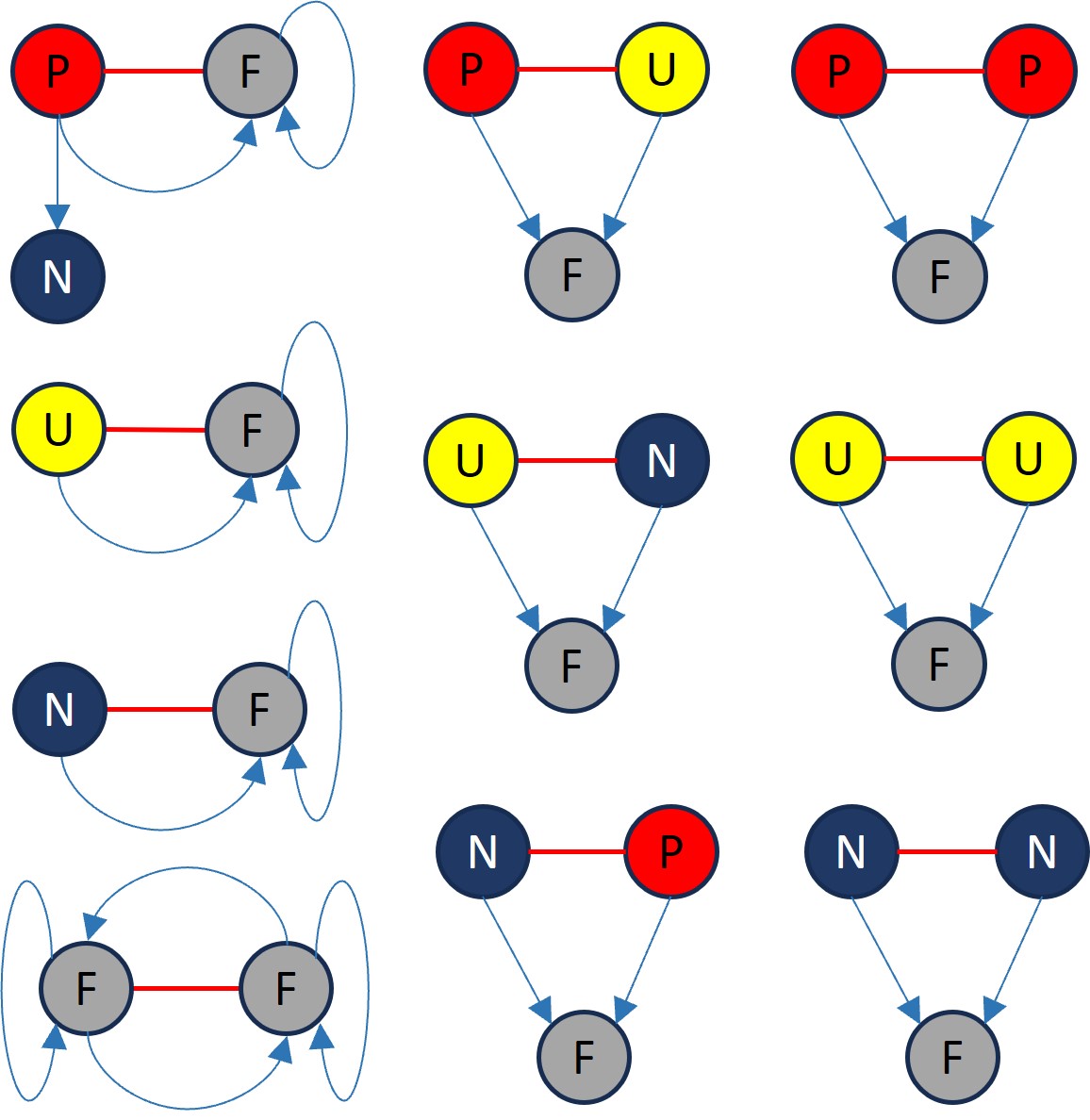}
        \caption{List of modelled interactions among the nodes due to collisions.}
        \label{fig:Patterns of interactions}
    \end{minipage}
\end{figure*}

In this work, the space environment is modelled with a network consisting of $N_n$ connected nodes. Each node $S_i$ represents $x_{S_i}(t_k)$ objects of a given class, or species, $S_i$ at time $t_k$. 
Each node belongs to an orbit site. In the general case, an orbit site is defined by a set of orbital parameters but, in this work, an orbit site is defined only by a combination of two parameters: altitude and inclination. There are $n$ orbit sites with $i=1,...,n$.
Figure \ref{fig:Diagram of network model in the physical layer} shows an illustration of the network model studied in this paper, where four classes of space objects are considered: Payloads ($P$), Upper stages ($U$), Fragments ($F$) and Non-manoeuvrable satellites ($N$). Therefore, in this representation, $S$ takes any value from the set $\mathcal{S}=\{P,U,N,F\}$. 

The number of objects $x_{S_i}(t_k)$ in a given node $S_i$ can change because of the combination of different events such as collisions, explosions, natural dynamics, manoeuvres, post-mission disposal strategies, operational lifetime duration, and new launches. Some of these events are continuous processes, like the orbital decay due to atmospheric drag, others are better represented with discrete jumps in the number of objects, like collisions and explosions.

In the network model, each event represents an interaction between node $S_i$ and other nodes or with itself (self-dynamics). We model these interactions with links connecting the nodes. Figure \ref{fig:Patterns of interactions} showcases the nine possible interactions related to collisions among the nodes we modelled in this paper. The red undirected links represent collisions and the blue directed links represent the flows of objects resulting from different types of events: collisions or explosions, that produce fragments, failed post mission disposal operations and impacts with small fragments, that produce non-manoeuvrable satellites, and natural orbit decay due to drag. Other events can be modelled but in this work we limited the model only to the ones listed above.   
For instance, with reference to Figure \ref{fig:Patterns of interactions}, a typical interaction comes from the aftermath of a collision between a $U$ and an $F$ node: the fragments generated during the collision will flow from both nodes to a set of $F$ nodes in different sites. The flow of fragments to nodes $F$ is due to the velocity distribution post collision, as it will be explained in the remainder of the paper. Another interaction occurs after a lethal non-trackable fragment collides with an active payload. We model this case with what we will call \textit{small collisions} between a node $P$ and a node $F$, resulting in objects flowing from $P$ to $N$ within the same site. 
It is important to underline that the network model allows for multiple representations of the space environment with further nodes of different nature, in addition to the ones that were considered in this study, and orbit sites that, in the general case, would be defined by combinations of intervals of orbital elements. For example, one could insert nodes representing entire constellations or single satellites, or nodes representing fragments or satellites in elliptical orbits. In this paper, we focus our attention only on four classes, or species, of objects and sites defined only by altitude and inclination.

\subsection{Network Dynamics}\label{sec: network dynamics}

As stated in the previous section the number of Payloads, $x_{P_i}(t_k)$, Upper stages, $x_{U_i}(t_k)$, Non-manoeuvrable satellites, $x_{N_i}(t_k)$, and Fragments, $x_{F_i}(t_k)$, at time $t_k$ in the $i$-th site, changes due to a combination of continuous and discrete events, represented by links. Discrete events like collisions and explosions are driven by random jump processes.
Therefore, once the structure and topology of the network are defined, one can write, for each orbit site, four stochastic evolutionary equations governing the variation of the number of objects of each species.

The number of Payloads, $x_{P_i}(t_k)$, in a given site $i$, decreases due to collisions with other objects, post mission disposal manoeuvres (PMD) and natural orbit decay, and increases due to new launches and payloads decaying from higher altitudes. Thus the stochastic equation governing the variation in the number of payloads can be written as: 
\begin{equation}\label{eq: payload dynamics}
    x_{P_i}(t_{k+1})=x_{P_i}(t_{k})-\sum_{j}^n\left(\Delta^{CAM}_{P_iU_j} + \Delta^{CAM}_{P_iN_j}+\Delta^{CAM}_{P_iF_j}+(1-s_{CAM})\Delta^{CAM}_{P_iP_j}\right) 
         -\sum_{j}^n\Delta^{\kappa}_{P_iF_j}- \Delta^{EOL}_{P_i}  + \varepsilon^+_{P_i} - \varepsilon^-_{P_i} +\Lambda_{P_i}(t_k),\\
\end{equation}
where $\sum_{j}^n\left(\Delta^{CAM}_{P_iU_j} + \Delta^{CAM}_{P_iN_j}+\Delta^{CAM}_{P_iF_j}+(1-s_{CAM})\Delta^{CAM}_{P_iP_j}\right)$ is the decrease due to possible collisions with upper stage, non-manoeuvrable satellites, other payloads and fragments in all orbit sites affecting node $P_i$, $\sum_{j}^n\Delta^{\kappa}_{P_iF_j}$ is the decrease due to the collisions with small fragments in all orbit sites affecting node $P_i$, $\Delta^{EOL}_{P_i}$ is the decrease due to post a mission disposal manoeuvre, $\varepsilon^+_{P_i}$ is the increase due to payloads decaying from higher altitudes within a time step $\Delta t_k$, $\varepsilon^-_{P_i}$ is the decrease due to orbital decay and $\Lambda_{P_i}(t_k)$ is the increase due to new launches within the same time step $\Delta t_k$.
The factor $s_{CAM}$ defines the fraction of payloads that perform a successful collision avoidance manoeuvre. The post mission disposal term $\Delta^{EOL}_{S_i}$ removes a payloads from the population at the end of a predefined period of time from its launch.
The general Poisson jump process $\Delta^{CAM}_{S_iS_j}$ defining the number of collisions between $S_i$ and $S_j$ subject to the successful execution of collision avoidance manoeuvres is defined as:
\begin{equation}\label{eq:coll_cam}
\Delta^{CAM}_{S_iS_j}=\int_{t_k}^{t_{k+1}} d\mathcal{N}((1-s_{CAM})\tau_{S_iS_j}\Delta t_k),
\end{equation}
and the number of collisions with small fragments is defined as:
\begin{equation}\label{eq:coll_small_frag}  
\Delta^{\kappa}_{S_iS_j}=\int_{t_k}^{t_{k+1}} d\mathcal{N}(\kappa \tau_{S_iS_j}\Delta t_k),
\end{equation}
where the small collision factor $\kappa$, quantifies the impacts of small fragments between 1 cm and 10 cm on payloads, leading to a transition of objects from $P$ nodes to $N$ nodes. Following \cite{2021Design}, we assumed that collisions with small fragments are $\kappa$ times more frequent than normal collisions with larger objects.
The use of Poisson jump processes was studied by \cite{pardini2014review} who demonstrated that it well captures the historical distribution of collisions. 
The number of Upper stages, $x_{U_i}(t_k)$, in site $i$, decreases due to catastrophic collisions, active debris removal (ADR), and natural orbit decay and increases due to new launches and to upper stages decaying from higher altitudes. Thus the stochastic equation governing the variation in the number of upper stages can be written as: 
\begin{equation}\label{eq:upper dynamics}
x_{U_i}(t_{k+1})=x_{U_i}(t_{k}) -\sum_{j}^n\Delta^{CAM}_{U_iP_j} - \sum^n_{j}(\Delta_{U_iN_j}+\Delta_{U_iF_j}+\Delta_{U_iU_j})- \Delta^{ADR}_{U_i} + \varepsilon^+_{U_i} - \varepsilon^-_{U_i}  + \Lambda_{U_i}(t_k),
\end{equation}
where the terms $\sum^n_{j}\Delta^{CAM}_{U_iP_j}$, $\sum^n_{j}(\Delta_{U_iN_j}+\Delta_{U_iF_j}+\Delta_{U_iU_j})$ and $\Delta^{ADR}_{U_i}$ are the variations due collisions with other species in all sites affecting node $U_i$ and the decrease due to ADR respectively. The terms $\varepsilon^+_{U_i}$ and $\varepsilon^-_{U_i}$ are the increase due to decay of upper stages from higher altitudes and the decrease due to orbital decay from site $i$, while $\Lambda_{U_i}(t_k)$ is the increase due to new launches. The general jump process due to collisions between nodes $S_i$ and $S_j$ is defined as:
\begin{equation}\label{eq:collision_jump}
\Delta_{S_iS_j}=\int_{t_k}^{t_{k+1}} d\mathcal{N}(\tau_{S_iS_j}\Delta t_k)
\end{equation}
The number of Non-manoeuvrable satellites, $x_{N_i}(t_k)$, in site $i$, decreases due to catastrophic collisions, active debris removal (ADR), and natural orbit decay and increases due to the fact that payloads can become non-manoeuvrable satellites due to collisions with small fragments or failed post mission disposal. Thus the stochastic equation governing the variation in the number of non-manoeuvrable satellites can be written as:
\begin{equation}
x_{N_i}(t_{k+1})=x_{N_i}(t_{k}) -\sum^n_{j}\Delta^{CAM}_{N_iP_j} - \sum_{j}(\Delta_{N_iN_j}+\Delta_{N_iF_j}+\Delta_{N_iU_j}) + \sum^n_{j}\Delta^{\kappa}_{P_iF_j} +\Delta^{\gamma}_{P_i}- \Delta^{ADR}_{N_i} + \varepsilon^+_{N_i} - \varepsilon^-_{N_i} + \Lambda_{N_i}(t_k)
\end{equation}
where the terms $\sum^n_{j}\Delta^{CAM}_{N_iP_j}$ and $\sum^n_{j}(\Delta_{N_iN_j}+\Delta_{N_iF_j}+\Delta_{N_iU_j})$ are the variations due to collisions with payloads and other objects, the terms $\sum^n_{j}\Delta^{\kappa}_{P_iF_j}$ and $\Delta^{\gamma}_{P_i}$ are the variations due to payloads impacting with small debris and failing a post mission disposal manoeuvre, respectively, $\Delta^{ADR}_{N_i}$, $ \varepsilon^+_{N_i}$ and $\varepsilon^-_{N_i}$ are the variations due to ADR, decay of non-manoeuvrable satellites from higher altitudes and orbital decay from site $i$, respectively. $\Lambda_{N_i}(t_k)$ is the increase due to new launches. The jump processes due to failed post-mission disposal is defined as:
\begin{equation}\label{eq:failed_PMD_jump}
\Delta^{\gamma}_{S_i}=\int_{t_k}^{t_{k+1}} d\mathcal{N}(\gamma \Delta^{EOL}_{P_i}),\\
\end{equation}
where the factor $\gamma$ indicates the failure rate in performing PMD, resulting in objects in $P$ nodes transitioning to $N$ nodes. 
Finally, the number of Fragments, $x_{F_i}(t_k)$ in site $i$ decreases due to collisions with other objects, ADR, and natural decays and increases due to collisions, explosions and natural decay from higher altitudes. Note that in this work, we have not modelled explosions yet. Thus the stochastic equation governing the variation in the number of fragments can be written as:
\begin{equation}\label{eq:fragments_dynamics}
x_{F_i}(t_{k+1})=x_{F_i}(t_{k}) -\sum^n_{j}\Delta^{CAM}_{F_iP_j} - \sum^n_{j}(\Delta_{F_iN_j}+\Delta_{F_iF_j}+\Delta_{F_iU_j})-   \Delta^{ADR}_{F_i} + \varepsilon^+_{F_i} - \varepsilon^-_{F_i} + \sum^n_{j}\Omega^{F_i}_j,
\end{equation}
where the terms $\sum^n_{j}\Delta^{CAM}_{F_iP_j}$ and $\sum^n_{j}(\Delta_{F_iN_j}+\Delta_{F_iF_j}+\Delta_{F_iU_j})$ are the decrease due to collisions with objects in all sites affecting $F_i$, $\Delta^{ADR}_{F_i}$ is the decrease due to ADR, $\varepsilon^+_{F_i}$ and $\varepsilon^-_{F_i}$ are the increase due to decay from higher altitudes and the decrease due to orbit decay from site $i$, respectively. The term $\sum^n_{j}\Omega^{F_i}_j$ is the collection of all the flows of fragments generated by all collisions in all sites affecting node $F_i$ and is defined as follows:
\begin{equation}\label{eq:OmegaF}
\begin{array}{l}
  \Omega^{F_i}_j = \sum_l \left(\sum_s^{\Delta^{CAM}_{P_jU_l}}\zeta^{F_i}_{s,P_jU_l} + \sum_s^{\Delta^{CAM}_{P_jN_l}}\zeta^{F_i}_{s,P_jN_l} + \sum_s^{\Delta^{CAM}_{P_jF_l}}\zeta^{F_i}_{s,P_jF_l} + \sum_s^{\Delta^{CAM}_{P_jP_l}}\zeta^{F_i}_{s,P_jP_l}\right) +\\ 
     \sum_l \left(\sum_s^{\Delta_{U_jU_l}}\zeta^{F_i}_{s,U_jU_l} + \sum_s^{\Delta_{U_jN_l}}\zeta^{F_i}_{s,U_jN_l} + \sum_s^{\Delta_{U_jF_l}}\zeta^{F_i}_{s,U_jF_l} + \sum_s^{\Delta_{N_jN_l}}\zeta^{F_i}_{s,N_jN_l} + \sum_s^{\Delta_{N_jF_l}}\zeta^{F_i}_{s,N_jF_l} + \sum_s^{\Delta_{F_jF_l}}\zeta^{F_i}_{s,F_jF_l}\right),
     \end{array}
\end{equation}
where $\zeta^{F_i}_{S_jS_l}$ indicates the fragments generated by the catastrophic collisions between node $S_j$ and node $S_l$ that flow to node $F_i$. Therefore, $\sum_s^{\Delta_{S_jS_l}}\zeta^{F_i}_{s,S_jS_l}$ indicates the new fragments flowing to node $F_i$ after the collision between node $S_j$ and node $S_l$ in the given time interval $\Delta t_k$.
In summary for each site $i=1,...,n$, we write the system of stochastic equations:
\begin{equation}\label{eq: nodes dynamics}
\begin{array}{l}
    x_{P_i}(t_{k+1})=x_{P_i}(t_{k})-\sum_{j}\left(\Delta^{CAM}_{P_iU_j} + \Delta^{CAM}_{P_iN_j}+\Delta^{CAM}_{P_iF_j}+(1-s_{CAM})\Delta^{CAM}_{P_iP_j}\right) 
         -\sum_{j}\Delta^{\kappa}_{P_iF_j}- \Delta^{EOL}_{P_i}  + \varepsilon^+_{P_i} - \varepsilon^-_{P_i} +\Lambda_{P_i}(t_k),\\
    x_{U_i}(t_{k+1})=x_{U_i}(t_{k}) -\sum_{j}\Delta^{CAM}_{U_iP_j} - \sum_{j}(\Delta_{U_iN_j}+\Delta_{U_iF_j}+\Delta_{U_iU_j})- \Delta^{ADR}_{U_i} + \varepsilon^+_{U_i} - \varepsilon^-_{U_i}  + \Lambda_{U_i}(t_k),\\
    x_{N_i}(t_{k+1})=x_{N_i}(t_{k}) -\sum_{j}\Delta^{CAM}_{N_iP_j} - \sum_{j}(\Delta_{N_iN_j}+\Delta_{N_iF_j}+\Delta_{N_iU_j}) + \sum_{j}\Delta^{\kappa}_{P_iF_j} +\Delta^{\gamma}_{P_i}- \Delta^{ADR}_{N_i} + \varepsilon^+_{N_i} - \varepsilon^-_{N_i} + \Lambda_{N_i}(t_k),\\
     x_{F_i}(t_{k+1})=x_{F_i}(t_{k}) -\sum_{j}\Delta^{CAM}_{F_iP_j} - \sum_{j}(\Delta_{F_iN_j}+\Delta_{F_iF_j}+\Delta_{F_iU_j})-   \Delta^{ADR}_{F_i} + \varepsilon^+_{F_i} - \varepsilon^-_{F_i} + \sum_{j}\Omega^{F_i}_j.\\
\end{array}
\end{equation}
Let $\mathbf{X}_{k}=[x_{P_1},x_{U_1},x_{N_1},x_{F_1}...,x_{P_i},x_{U_i},x_{N_i},x_{F_i},...,x_{P_n},x_{U_n},x_{N_n},x_{F_n}]^T$ represent the vector that contains all the nodes $x_{S_i}$, with $S\in\{P,U,N,F\}$ and $i=1,...,n$, then a more compact form of the evolutionary equations in (\ref{eq: nodes dynamics}) is:
\begin{equation}
\label{eq: change in the nodes}
    \mathbf{X}_{k+1} = \mathbf{X}_{k}+ \mathbf{Y}(\mathbf{X}_k)+ \mathbf{g}(t_k),
\end{equation}
where $\mathbf{Y}(\mathbf{X}_k)$ contains the interactions among nodes (collisions and orbit change) and satellite failures, while $\mathbf{g}(t_k)$ contains all time dependent terms, such as the launch traffic model $\Lambda_i(t_k)$ and end of life $\Delta^{EOL}_{S_i}$. The launch traffic model, the collision model, and the decay model will be introduced in the following sections. 

As it will be explained later the nodes affected by ADR can be chosen so that one can remove the most critical objects.  The mean collision rate $\tau_{S_iS_j}$ between node $S_i$ and node $S_j$ and the fragments generated in a collision will be defined in section \ref{sec:collision_rate}, the function $\Lambda_{P_i}(t)$ will be defined in section \ref{sec:Launch Traffic Model} and the quantities $\varepsilon_{S_i}^-$ and $\varepsilon_{S_i}^+$ will be defined in section \ref{sec:Decay Model}. 

 The evolutionary model in (\ref{eq: nodes dynamics}) can be understood as a stochastic version of what is commonly known as source-sink model. The key difference with respect to previous works, like \cite{dambrosio2023}, is that we derived the evolutionary equations from the structure of the network. This also differentiates our approach from previous network representations of the space environment \citep{lewis2010new} that use network theory as a post-processing analysis tool.
In our work we first generate an abstract network representation of the physical space environment. We then derive the evolutionary equations of the abstract network representation. This modelling approach allows us to flexibly balance computational complexity and accuracy, because nodes can represent groups of objects of different size and can belong to sites with different volumes. We then study the links between pairs of nodes, which are orders of magnitude lower, in number, than the pairs of objects. More importantly, we can directly follow the evolution of the relationships among nodes and sites, and extract, alongside typical quantities like the total number of objects of a given species, information that quantifies the interdependency of objects and orbital sites. In the remainder of the paper we will show that even with a course network representation we can obtain results that are comparable with a more traditional cube based model and we can quantify the importance of the relationship between given orbit sites and the rest of the environment. By changing these relationships, one can control the evolution of the network and, therefore, of the environment.

\subsubsection{Launch Traffic Model}\label{sec:Launch Traffic Model}
\begin{figure*}
    \centering    \includegraphics[width=0.6\textwidth]{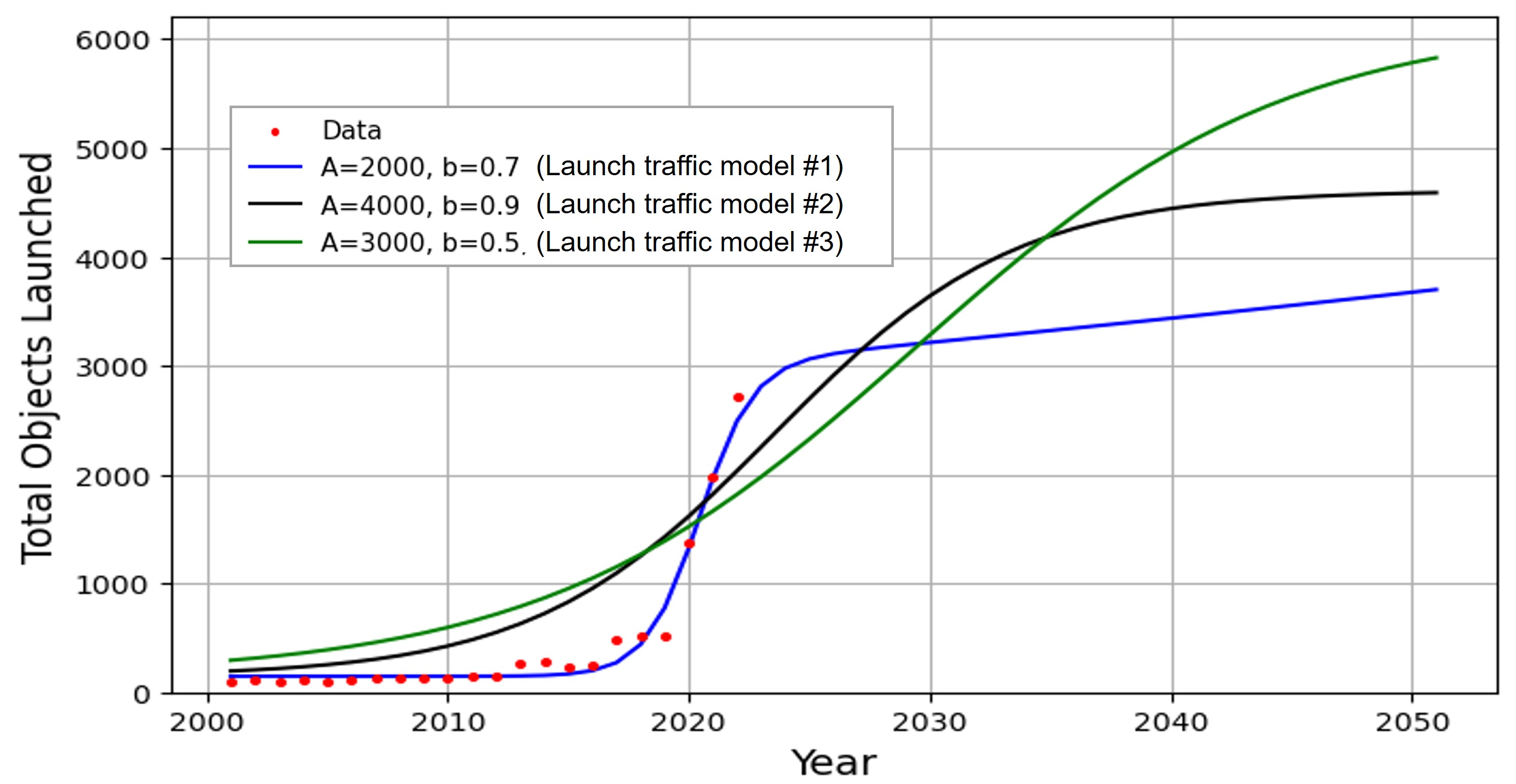}
    \caption{Three launch traffic forecasts of  objects launched into LEO.}
    \label{fig:Model of number of objects launched to LEO}
\end{figure*}

The launch traffic model is based on historical launch data and estimates the future launch rates. Historical data to fit the model are taken from DISCOSWeb. The total number of objects launched in a certain year is given by an exponential logistic curve of the form \citep{wilson2024modelling}
\begin{equation}\label{eq:traffic_model}
    \Lambda(t) = n_0 + \sum_{i=1}^M \frac{A_i e^{d_i(t-t_i)}}{b_i+e^{-c_i(t-t_i)}},
\end{equation}
where $t$ is expressed in years and the parameters $n_0$, $A_i$, $b_i$, $c_i$, $d_i$, and $t_i$ control the shape of the curve. These parameters can be selected based on historical trends and predicted future increases in launch traffic. Launched objects are either payloads, upper stages, or mission-related objects. The proportion of each is based on the most recent years of launches. For each object class, orbital and physical parameters of historical launches are fit to a Gaussian Mixture Model (GMM). 
{A separate GMM was fit to the distribution of semimajor axis and inclination, and to physical parameters like mass, cross-sectional area, and length. The increase in objects with a given set of physical properties sampled from each GMM, follows the model in (\ref{eq:traffic_model}). Hence, one can compute the expected number of new objects launched in a given site from the percentage of payloads and upper stages of a given size derived from the GMM fit, times the increase given by model (\ref{eq:traffic_model}).}  Figure \ref{fig:Model of number of objects launched to LEO} shows three example of launch traffic forecasts for $M=1$ and for objects launched into Low Earth Orbit. Also shown are the parameters used to model the different forecasts. The choice of parameters reflects different qualitative assumptions on how the number of objects launched will evolve in the future. The blue line suggests a steadier increase in new launches in future following the sharp upward trend since 2023 that representing a timely restriction on launch activity, whereas the black and red lines suggest that the number of objects reaching a plateau occurs later, at different levels, with the black line plateauing around 2032 and the red line around 2050. In Section \ref{Sec:Assessment of the space environment stability} we will use these three illustrative examples of launch traffic forecasts to propagate the space environment with NESSY. Note that for $M>1$ one can simulate the sudden increase in number of launches or launched objects, which can correspond to new mega-constellations or simply a rapid growth of the sector.
In the following, the GMM is fit to the altitude and inclination of historical launches. The launch traffic model samples from this distribution and adds a number of objects to each site, dictated by (\ref{eq:traffic_model}). 

\subsubsection{Collision Rate}\label{sec:collision_rate}
The mean collision rate between  node $S_i$ and node $S_j$ is defined as 
\begin{equation}
\label{eq: collision rate}
   \tau_{{S_i}{S_j}} =  \frac{n_{S_i} n_{S_j}\Delta v_c\sigma_{{S_i}{S_j}}}{\max(V_{S_i},V_{S_j})}, 
\end{equation}
see also \cite{pardini2014review}, where $n_{S_i}$ and $n_{S_j}$ are the number of objects respectively in nodes $S_i$ and $S_j$, $\sigma_{S_iS_j}$ is the average cross-section area of the colliding objects, $\Delta v_c$ is the average collision velocity and the collision is assumed to happen in the smaller between volumes $V_{S_i}$ and $V_{S_j}$. Note that when $S_i=S_j$ , $n_{S_j}=(n_{S_i}-1)/2$. This formulation is analogous to the one proposed in \cite{KESSLER1981}, although, as it will be explained in the remainder of this section, we use a different approach to compute the average collision velocity and cross-section area. We also assume that only the objects that reside in the smaller between volume $V_{S_i}$ and $V_{S_j}$ at a given time step, can collide. These objects include both the ones in a given orbit shell and those, on elliptical orbits, that are transiting through the same orbit shell. See section \ref{sec:Treatment of Elliptical Orbits} for more details on the treatment of objects on elliptical orbits.

The volume of a generic node $S$ defined by an orbit site with altitude thickness $dr$ and inclination bin $i\in[i_{min},i_{max}]$ is:
\begin{equation}
\label{eq: node volumn}
V_{S}=4\pi \sin(i_{max})\left[(r+dr)^3 - r^3\right]/3. 
\end{equation}
Note that if $\pi/2\in[i_{min},i_{max}]$ then $i_{max}=\pi/2$ and for $i_{max}\wedge i_{min}>\pi/2$ we use $i_{max}=\pi-i_{min}$. 
Given object $q$ and $p$ in the same node $S_i$, the average cross-section area $\sigma_{{S_i}{S_j}}$ is
\begin{align}
\label{eq: self-collision sigma}
    \sigma_{{S_i}{S_j}}=\frac{\pi}{4 (n^2_{S_i} - n_{S_i})}\sum^{n_{S_i}}_{q\neq p}\sum^{n_{S_i}}_{p \neq q} \left(d^q_{S_i}+d^p_{S_i}\right)^2 &=\frac{\pi}{4 (n^2_{S_i} - n_{S_i})} \sum^{n_{S_i}}_{q\neq p}\sum^{n_{S_i}}_{p \neq q} \left[(d^q_{S_i})^2 + (d^p_{S_i})^2 + 2d^q_{S_i}d^p_{S_i} \right] \notag\\
    &=\frac{\pi}{4 (n^2_{S_i} - n_{S_i})} \left[(n_{S_i}-1)\sum^{n_{S_i}}_{q}(d^q_{S_i})^2 + (n_{S_i}-1)\sum^{n_{S_i}}_{p}(d^p_{S_i})^2 + 2\sum^{n_{S_i}}_q\sum^{n_{S_i}}_p d^q_{S_i}d^p_{S_i}-2 \sum^{n_{S_i}}_q (d^q_{S_i})^2\right] \notag\\
    &=\frac{\pi}{4 (n^2_{S_i} - n_{S_i})}\left[ 2(n_{S_i}-2)\sum^{n_{S_i}}_q (d^q_{S_i})^2 + 2 (\sum^{n_{S_i}}_q d^q_{S_i})^2\right],
\end{align}
while given object $q$ and $p$ in two different nodes $S_i$ and $S_j$ the average cross-section area is
\begin{align}
\label{eq: cross-collision sigma}
    \sigma_{{S_i}{S_j}}=\frac{\pi}{4n_{S_i}n_{S_j}}\sum^{n_{S_i}}_{q=1}\sum^{n_{S_j}}_{p=1} \left(d^q_{S_i}+d^p_{S_j}\right)^2 
    &= \frac{\pi}{4n_{S_i}n_{S_j}} \sum^{n_{S_i}}_{q}\sum^{n_{S_j}}_{p} \left[(d^q_{S_i})^2 + (d^p_{S_j})^2 + 2d^q_{S_i}d^p_{S_j} \right] \notag\\
    &= \frac{\pi}{4n_{S_i}n_{S_j}} \left[ n_{S_j}\sum^{n_{S_i}}_q (d^q_{S_i})^2 + n_{S_i} \sum^{n_{S_j}}_p (d^p_{S_j})^2 + 2\sum^{n_{S_i}}_q d^q_{S_i}\sum^{n_{S_j}}_pd^p_{S_j}\right].
\end{align}

\begin{figure*}
    \centering    \includegraphics[width=0.5\textwidth]{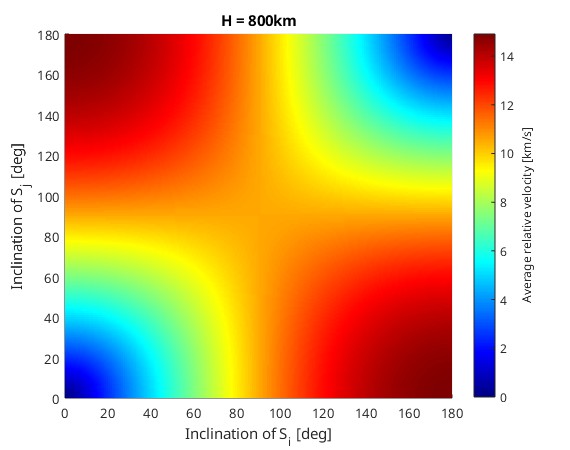}
    \caption{Average relative velocity for different combinations of inclinations at an altitude of 800 km}
    \label{fig:Average relative velocity of different combination of inclination}
\end{figure*}

The average relative collisional velocity $\Delta v_{c}$ is computed starting from
\begin{equation}
\label{eq: relative velocity}
    \Delta v_{c} \approx \sqrt{v_{S_i}^2 + v_{S_j}^2 - 2v_{S_i} v_{S_j} \mathbb{E}[\cos(\theta_{S_iS_j})]},
\end{equation}
where $v_{S_i}$ and $v_{S_j}$ are the representative velocities of nodes $S_i$ and $S_j$. The representative velocity of node $S_i$ is taken to be 
\begin{equation}
    v_{S_i} = \sqrt{\mu\left(\frac{2}{r_l} -\frac{1}{a_{S_i}}\right)},
\end{equation}
where $a_{S_i}$ is the average semi-major axis of nodes $S_i$ and $r_l$ is the mean radius of the corresponding orbit shell. Quantity $\mathbb{E}[\cos(\theta_{S_iS_j})]$ is the expected value of the cosine of the angle between the velocity vectors at collision. Given the inclination $i$, and Right Ascension of the Ascending Node $\Omega$, consider object $q$ in node $S_i$ and object $p$ in node $S_j$ then
\begin{equation}
\label{eq: cos thetaij}
    \cos(\theta^{qp}_{S_iS_j}) = \begin{bmatrix}
\sin i^q_{S_i} \cos \Omega^q_{S_i} \\
\sin i^q_{S_i} \sin \Omega^q_{S_i} \\
\cos i^q_{S_i}
\end{bmatrix}^T
\begin{bmatrix}
\sin i^p_{S_j} \cos \Omega^p_{S_j} \\
\sin i^p_{S_j} \sin \Omega^p_{S_j} \\
\cos i^p_{S_j}
\end{bmatrix} 
= \sin i^q_{S_i} \sin i^p_{S_j} \cos{(\Delta \Omega^{qp}_{S_iS_j})} + \cos i^q_{S_i} \cos i^p_{S_j},
\end{equation}
where $\Delta \Omega^{qp}_{S_iS_j} = \Omega^q_{S_i} - \Omega^p_{S_j}$. If the distributions of inclination and right ascension of the objects in node $S_i$ are independent and uncorrelated from the ones of the objects in node $S_j$ the expected value $\mathbb{E}[\cos(\theta_{S_iS_j})]$ reads
\begin{align}
    \mathbb{E}[\cos(\theta_{S_iS_j})] &= \mathbb{E}[\sin i^q_{S_i}] \mathbb{E}[\sin i^p_{S_j}] \mathbb{E}[\cos{(\Delta \Omega^{qp}_{S_iS_j})}]+ \mathbb{E}(\cos i^q_{S_i}) \mathbb{E}(\cos i^p_{S_j})  \notag \\
    & = \frac{1}{n_{S_i}}\sum_q\sin i^q_{S_i}\frac{1}{n_{S_j}}\sum_p \sin i^p_{S_j} \frac{1}{n_{S_i}n_{S_j}}\sum_q \sum_p \cos{(\Delta \Omega^{qp}_{S_iS_j})}+ \frac{1}{n_{S_i}n_{S_j}}\sum_q\sum_p\cos i^q_{S_i}\cos i^p_{S_j},
    \label{eq: expectation of cos_theta}
\end{align}
We now assume that $\Omega^q_{S_i}$ and $\Omega^p_{S_j}$ are independent and uniformly distributed over the interval $[0, 2\pi]$, then also $\Delta \Omega^{qp}_{S_iS_j}$ is uniformly distributed over the interval $[0, 2\pi]$. Therefore,
\begin{equation}
    \frac{1}{n_{S_i}n_{S_j}}\sum_q \sum_p \cos{\left(\Delta \Omega^{qp}_{S_iS_j}\right)} \approx \frac{1}{2\pi}\int_0^{2\pi} \cos\left(\Delta \Omega^{qp}_{S_iS_j}\right) \, d \Delta \Omega^{qp}_{S_iS_j} = 0,
\end{equation}
and if the inclination bins are sufficiently small 
\begin{equation}
    \frac{1}{n_{S_i}n_{S_j}}\sum_q\sum_p\cos i^q_{S_i}\cos i^p_{S_j} \approx \cos \left(\frac{1}{n_{S_i}} \sum_q i^q_{S_i}\right)\cos\left(\frac{1}{n_{S_j}}\sum_pi ^p_{S_j}\right),
\end{equation}
from which Eq.\eqref{eq: expectation of cos_theta} reduces to
\begin{equation}\label{eq:exp_costheta}
     \mathbb{E}[\cos(\theta_{S_iS_j})] \approx \cos i_{S_i} \cos i_{S_j},
\end{equation}
where $i_{S_i}$ and $i_{S_j}$ are the average inclinations of the objects in each node. Figure \ref{fig:Average relative velocity of different combination of inclination} shows the average relative collision velocity for different combinations of orbital inclinations, assuming both nodes are located at an altitude of 800 km. 
Table~\ref{table:validation of relative velocity}, shows the comparison between the relative collision velocities, computed with approximation (\ref{eq: relative velocity}), of three pairs of objects, located between 800 km and 810 km and with an eccentricity ranging from 0 to 0.1, and the same relative collision velocities computed by averaging all the relative velocities between 200 randomly generated objects in node $S_i$ colliding with other 200 randomly generated objects in $S_j$, for a total of 40,000 collision pairs. The simulation results demonstrate that the proposed approximation method achieves a relative error of less than 10\%, compared to an exact calculation of the average collision velocity, for inclination bins of 10 degrees. For larger bins the relative error is expected to be larger, however, note that, since the inclination goes from 0 to $\pi$, the relative approximation error in Eq. (\ref{eq: expectation of cos_theta}) remains around 10\% even for a bin size of 90 degrees.

\begin{table}[h]
\centering
\caption{Comparison of average relative collision velocity $\Delta v_c$ for three sample pairs of nodes in the same orbit shell.}
\label{table:validation of relative velocity}
\begin{tabularx}{\textwidth}{@{}XXccc@{}}
\toprule
\textbf{Inclination bin of $S_i$ [deg]} & 
\textbf{Inclination bin of $S_j$ [deg]} & 
\textbf{ exact $\Delta v_c$ [km/s]} & 
\textbf{approximated $\Delta v_c$ [km/s]} & 
\textbf{Relative error [\%]} \\
\midrule
$[0,10]$      & $[0,10]$        & 0.8915  & 0.8390   & 5.89 \\
$[0,10]$      & $[170,180]$     & 14.8880 & 14.8850  & 0.02    \\
$[60,70]$     & $[110,120]$     & 11.9507 & 11.4403  & 4.27    \\
\bottomrule
\end{tabularx}
\end{table}

\begin{figure*}[h!]
\centering
    \centering
    \begin{subfigure}{0.49\textwidth}
    \includegraphics[width=1\textwidth]{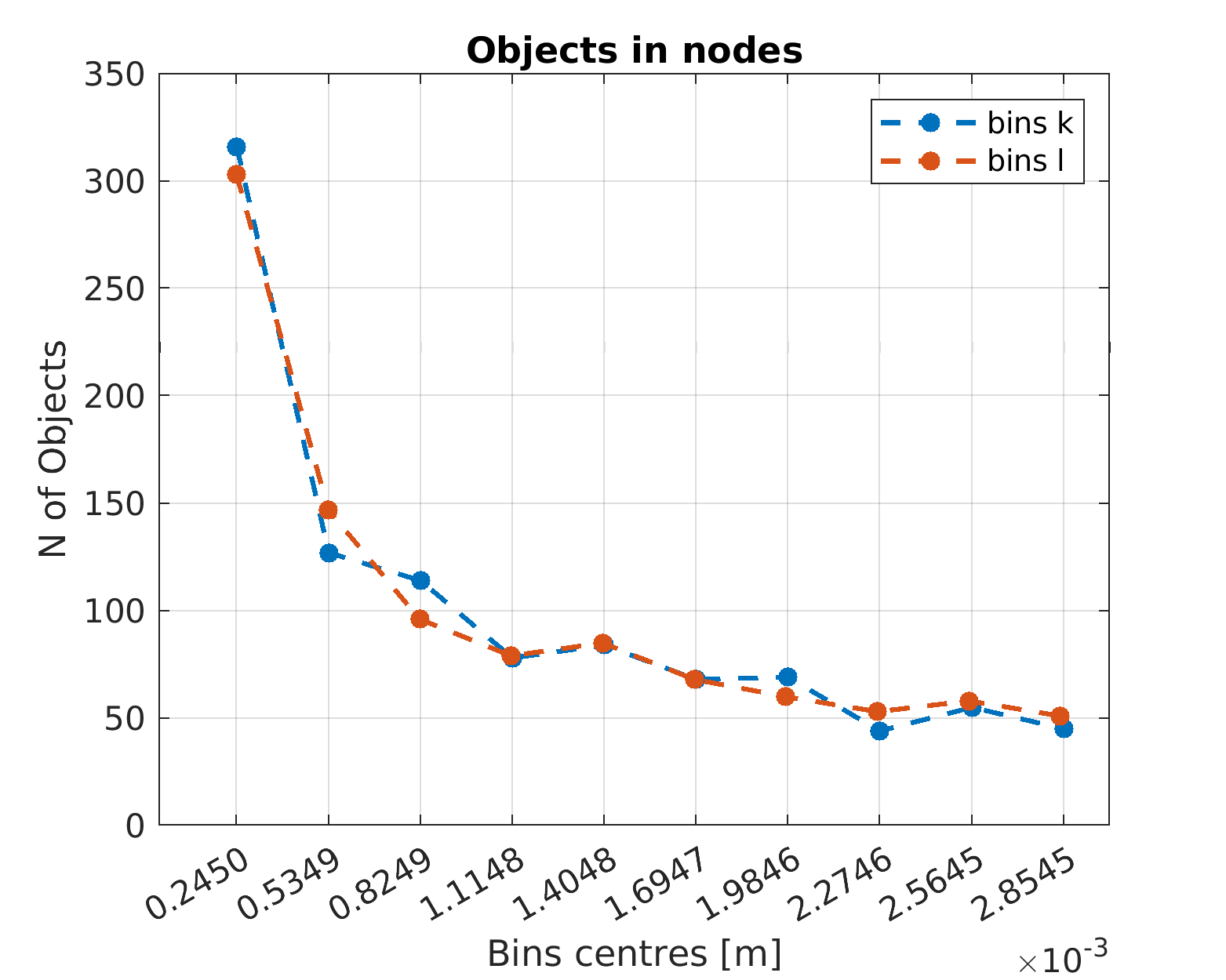}  \caption{\label{fig:pop_nodes}Distributions of the objects' radii for two examples of nodes $S_i$ and $S_j$.}
    \label{fig:objs_bins}
    \end{subfigure}
    \begin{subfigure}{0.495\textwidth}
    \centering
    \includegraphics[width=1\textwidth]{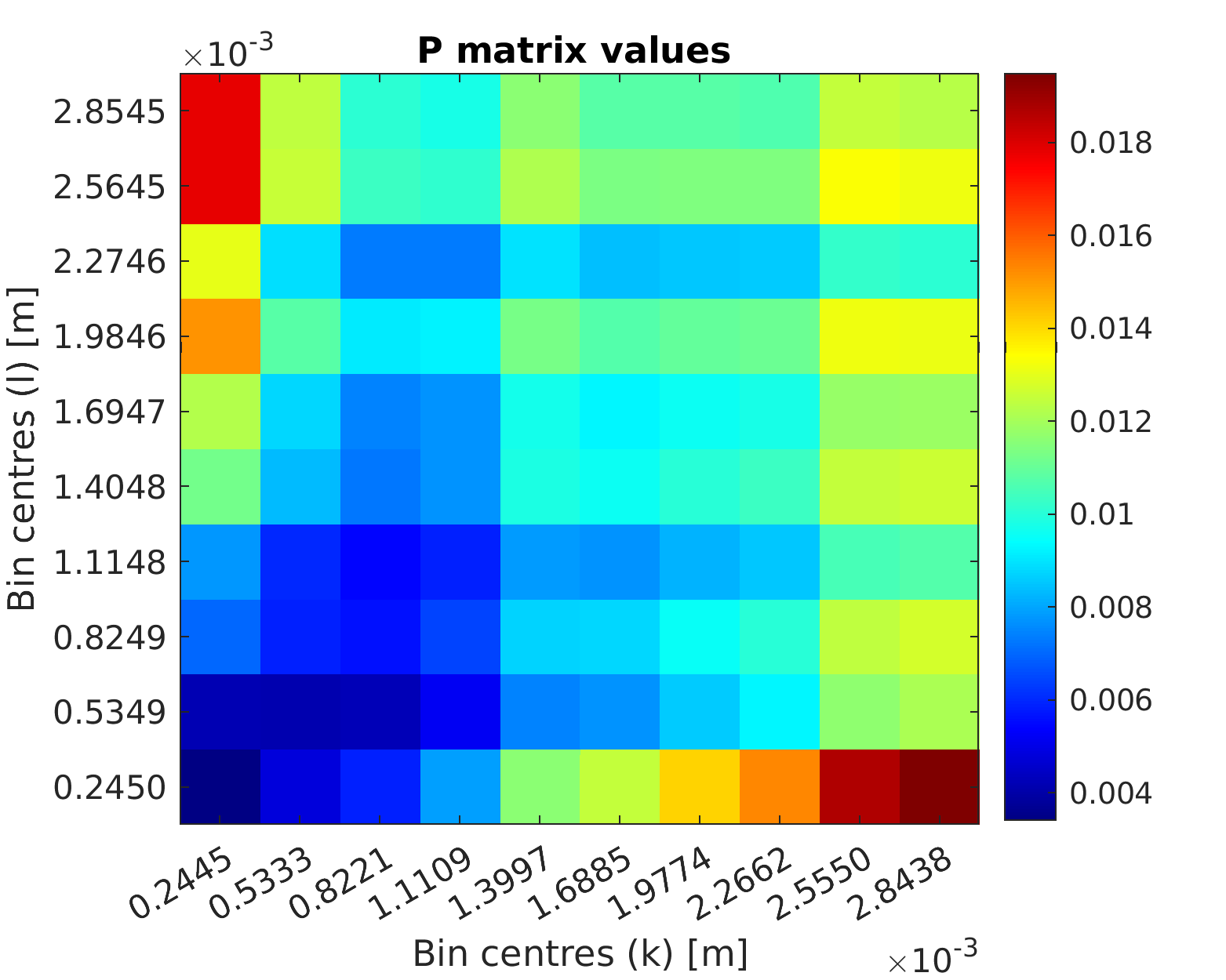}
    \caption{\label{fig:p_matrix} Probability matrix P.}
    \label{fig:sel_strg}
    \end{subfigure}
    \caption{Cross-collision test case between two test nodes $S_i$ and $S_j$: a) distribution of the radii in each node and b) probability matrix used to select the colliding objects.}
\end{figure*}

The collision rate (\ref{eq: collision rate}) is used in Eq.(\ref{eq:collision_jump}) to compute the number of collisions $N_C$ within the same node $S_i$, or self-collisions, or in between two nodes, or cross-collisions. Once the number of collisions $N_C$ is established the objects involved in a collision are selected as follows. The objects in each node are partitioned in $N_D$ bins according to their radius $R$, spanning the range $R_{S_i}^{min}$ and $R_{S_i}^{max}$ between the smallest and the largest objects in each node. We then draw $N_C$ samples from the distribution over the indexes of the bins:
\begin{equation}
    P(k,l) = 1 - \exp\left( -\frac{n_{k,S_i} n_{l,S_j}\Delta v_c \sigma_{kl}}{\max(V_{S_i},V_{S_j})}\Delta t\right),
\end{equation}
with $k=1,...,N_D$ and $l=1,...,N_D$, where $\sigma_{kl}$ is the average cross-sectional area for each couple of bins with index $k$ and $l$, computed in Eq.\eqref{eq: self-collision sigma} or Eq.\eqref{eq: cross-collision sigma}. The quantities $n_{k,S_i}$ and $n_{l,S_i}$ are the number of objects in each bin with index $k$ and $l$ respectively. Once the indices of two size bins are selected, the specific objects involved in the collisions are randomly drawn from a uniform distribution over the objects in those bins. Figures \ref{fig:objs_bins} and \ref{fig:p_matrix} show a simple application of the selection strategy for a cross collision case, using 10 bins per node and placing 1000 objects in each node. Objects' dimensions are chosen to favour smaller ones (between 0.1 and 3 meters). Figure \ref{fig:objs_bins} shows the distribution of the objects in the bins of the two nodes, while Figure \ref{fig:p_matrix} shows the values of the matrix $P(k,l)$. The time step $\Delta t$ is 30 days. The altitude range of the nodes considered in this example is 500-550 km, and the inclination range is 0-60 degrees for one node and 60-120 degrees for the other node. When each bin contains a single object, $P$ converges to the distribution over the size of all the pairs of objects in nodes $S_i$ and $S_j$. In the remainder of the paper, we used 50 bins per node, which can be shown to be a good compromise between accuracy and computational complexity.

The number of fragments generated by a given collision event is approximate with, see \cite{johnson2001nasa}:
\begin{equation}
\label{eq: number of fragments after the collision}
    \zeta_{{S_i}{S_j}}(L_c) = 0.1 M_{{S_i}{S_j}}^{0.75} L_c^{-1.71},
\end{equation}
where $L_c$ indicates the minimum size of the fragments. Given that a catastrophic collision occurs when the impact kinetic energy-to-target mass ratio exceeds \SI{40}{J/g}, $M_{{S_i}{S_j}}$ is defined as the mass of both objects in a catastrophic collision, while in the case of a non-catastrophic collision, the value of $M_{{S_i}{S_j}}$ is defined as the mass of the smaller object multiplied by the square of collision velocity \citep{krisko2011proper}. The physical parameter distribution and the velocity distribution of the resulting fragments can be computed by NASA's breakup model \citep{johnson2001nasa}.  Considering orbital elements of the fragments, only a portion of the fragments from each collision flows to $F_k$, therefore, the notation $\zeta^{F_k}_{S_iS_j}$ is introduced to represent the number of fragments, generated from the collision between node $S_i$ and node $S_j$, that flows to node $F_k$. The conservation of mass is considered to ensure that the total mass of the fragments remains consistent before and after the collision. Taking a catastrophic collision, as an example, if the total mass of the resulting fragments is less than the combined mass of the parent bodies, the remaining mass is redistributed across 2 to 8 randomly sampled fragments, based on the properties of the target \citep{jang2025new}.



\subsubsection{Decay Model}\label{sec:Decay Model}
Building on the averaging perturbation technique proposed by \cite{martinusi2015analytic}, a time-explicit analytic approximate solution, see \cite{jang2023monte}, for the motion of low-Earth-orbiting satellites is used. Given $[a_k, e_k]$ at the current time $t_k$, the $a_{k+1}$ after $\Delta t_k$ is
\begin{equation}
\begin{split}
    a_{k+1} &= \frac{a_k}{\beta_k^2} \tan^2 \left[ \arctan(\beta_k) - \beta_k n_k a_k C_k \Delta t_k \right],\\
    e_{k+1} &= \frac{2}{\sqrt{3}} \tan \left[ \arctan(\beta_k) - \beta_k n_k a_k C_k \Delta t_k \right],
\end{split}
\end{equation}
where
\begin{equation}
    C_k = \frac{1}{2} C_d \frac{A}{m} \rho, \qquad \bar{n}_k = \sqrt{\frac{\mu}{a_k^3}}, \qquad \beta_k = 
    \begin{cases}
        \frac{\sqrt{3}}{2} e_k, &e_k \neq 0\\
        1, &e_k = 0
    \end{cases},
\end{equation}
where $\rho$ indicates the atmospheric density. This paper employs the Jacchia-Bowman 2008 (JB2008) density model \citep{bowman2008new}, an empirical model derived from historical atmospheric density data, to represent the average statistical behaviour of the atmosphere under varying solar and geomagnetic influences \citep{jang2023monte}. 
Assuming the altitude of object $q$ is $h^q(t_{k+1})$ at $t_{k+1}$, the quantity $\omega^q_{{S_i}{S_{j}}}(t_{k+1})$ is an indicator that says if object $q$ in node $S_i$ can flow from node $S_i$ to the node $S_{j}$ within $\Delta t_k$: 
\begin{equation}
\label{Eq: indicator for decay}
    \omega^q_{{S_i}{S_{j}}}(t_{k+1}) = 
    \begin{cases}
        1, \quad if \;  \underline{h}_{Sj} < h^q(t_{k+1}) < \overline{h}_{Sj} \\
        0, \quad else
    \end{cases},
\end{equation}
where $\overline{h}_{S_j}$ indicates the upper bound of the orbital shell of node $S_j$, respectively. The flow from node $S_i$ to node $S_j$ is
\begin{equation}
    \varepsilon_{S_iS_j}(t_{k+1}) = \sum^{n_{S_i}}_{q=1} \omega^q_{{S_i}{S_{j}}}(t_{k+1}),
\end{equation}
Therefore, the outflow $\varepsilon_{S_i}^-$ and inflow $\varepsilon_{S_i}^+$ are
\begin{equation}
\begin{cases}
    \varepsilon^{-}_{S_i}=\sum^{n}_{j\neq i} \varepsilon_{S_iS_j}(t_{k+1})\\
    \varepsilon^{+}_{S_i}=\sum^{n}_{j\neq i} \varepsilon_{S_jS_i}(t_{k+1})\\
\end{cases}.
\end{equation}

\subsubsection{Treatment of Elliptical Orbits}\label{sec:Treatment of Elliptical Orbits}
Considering an object $q$ on an elliptical orbit, and assuming that the probability of its angular position is distribution along the entire orbit, we want to compute the fraction of that probability within each orbital shell at a given time, see \cite{lucken2019collision}. We first define the distribution function: 
\begin{equation}
    g^q(\theta, e) =  2 \arctan \left[ \left( \frac{1 - e}{1 + e} \right)^{\frac{1}{2}} \tan \frac{\theta}{2} \right]- \frac{e(1 - e^2)^{\frac{1}{2}} \sin \theta}{1 + e \cos \theta},
\end{equation}
where $e$ is the eccentricity. Then, the proportion of residence time of a given object in orbital shell $\nu$ is 
\begin{equation}
w^q_{\nu} = \frac{g^q(\theta_{\nu+1}, e_{k+1}) - g^q(\theta_{\nu}, e_{k+1})}{\pi},
\end{equation}
where $\theta_{\nu}$ is the phase angle corresponding to where the orbit crosses the lower edge of the $\nu$-th shell. Given the weights $w^q_{\nu}$, a sampling strategy without replacement is employed to select the orbital shell $\nu$ in which the object is most likely to be found at time $t_{k+1}$. Say that the mid-point altitude of the selected orbit shell is  
$h^q(t_{k+1})$, then by substituting $h^q(t_{k+1})$ into Eq.\eqref{Eq: indicator for decay},  one can compute an updated flow $\varepsilon^{-}_{S_i}$ and $\varepsilon^{+}_{S_i}$ with the inclusion of the
flow between orbital shells induced by orbital eccentricity of object $q$.

\section{Comparison of NESSY Against MOCAT-MC.}

\label{sec:validation_scenario}


In this section, we compare NESSY against an already well-tested evolution model: MOCAT-MC, \citep{jang2023monte, jang2025new}. The initial condition for both NESSY and MOCAT-MC is the population of resident objects in LEO at the beginning of 2023, available from SatCat TLEs \citep{stein_l_leo_2022} and ESA DISCOS \citep{esa_space_debris_office_discosweb_2024}. 
The portion of the space environment covered in this comparison spans a range of altitudes from 200 to 2200 km. The composition of the initial population is reported in Table \ref{table:pop_nessy}.  In the network model the population is divided in four nodes, corresponding to four species, for each inclination bin and orbit shell.

\begin{table}[h!]
\centering
\caption{Composition of the initial population in March 2023 (SatCat and ESA DISCOS data).}
\begin{tabular}{c c c c c} 
 \hline
 \textbf{Payloads} & \textbf{Upper-stages} & \textbf{Fragments} & \textbf{Non-maneuverable satellite} & \textbf{Total} \\ [0.5ex] 
 \hline
 5471 & 1111 & 9804 & 2440 & 18 826 \\ [1ex] 
 \hline
\end{tabular}
\label{table:pop_nessy}
\end{table}

\subsection{Long Term Evolution Results}

Both MOCAT-MC and NESSY use the same atmospheric model and for this comparison we did not include any future launch. For both NESSY and MOCAT-MC we assumed that each payload has a mission lifetime of 5 years. At the end of its mission lifetime, each payload is assumed to performs a Post-Mission Disposal (PMD) manoeuvre. A successful PMD results in the immediate removal of the object from the population, while a failed PMD transforms a payload into a non-manoeuvrable satellite, without changing orbit site. The small collision factor $\kappa$, which quantifies the impacts of small fragments between 1 cm and 10 cm on payloads, was set to 5.3, as in \cite{2021Design}. The success rate of collision avoidance manoeuvrers, $s_{CAM}$, was set to 99.9 \%. The new fragments have a diameter larger than 10 cm were considered. For an initial comparison with MOCAT-MC, a refined discretisation was tested with 10 km altitude shells (from 200 to 2200 km) and 10 deg inclination bins (from 0 deg to 180 deg). The settings used in NESSY and MOCAT-MC for this comparative test are summarised in Table \ref{table:params_valid}.

\begin{table}[h!]
\begin{center}
\caption{Settings of NESSY and MOCAT-MC for the comparison.}
\begin{tabular}{c c c} 
 \hline
 \textbf{Parameter} & \textbf{NESSY} & \textbf{MOCAT-MC} \\ [0.5ex] 
 \hline
 Time step (days) & 30 & 5 \\ 
 
 Propagation time [years] & 100 & 100 \\
 
 Inclination bin size [deg] & 10 & - \\
 
 Altitude shell size [km] & 10 & - \\
 
 $s_{CAM}$ [\%] & 99.99 & 99.99 \\ [1ex] 
 
 $\gamma$ [\%] & 5  & 5 \\ [1ex] 
 
 Payloads mission lifetime [years] & 5 & 5 \\ [1ex] 

  Small collisions factor & 5.3 & - \\ [1ex] 

  Explosion probability per day & 0 & 0 \\ [1ex] 
 \hline
\end{tabular}
\label{table:params_valid}
\end{center}
\end{table}

\begin{figure}[H]
    \centering
    \includegraphics[width=0.9\textwidth]{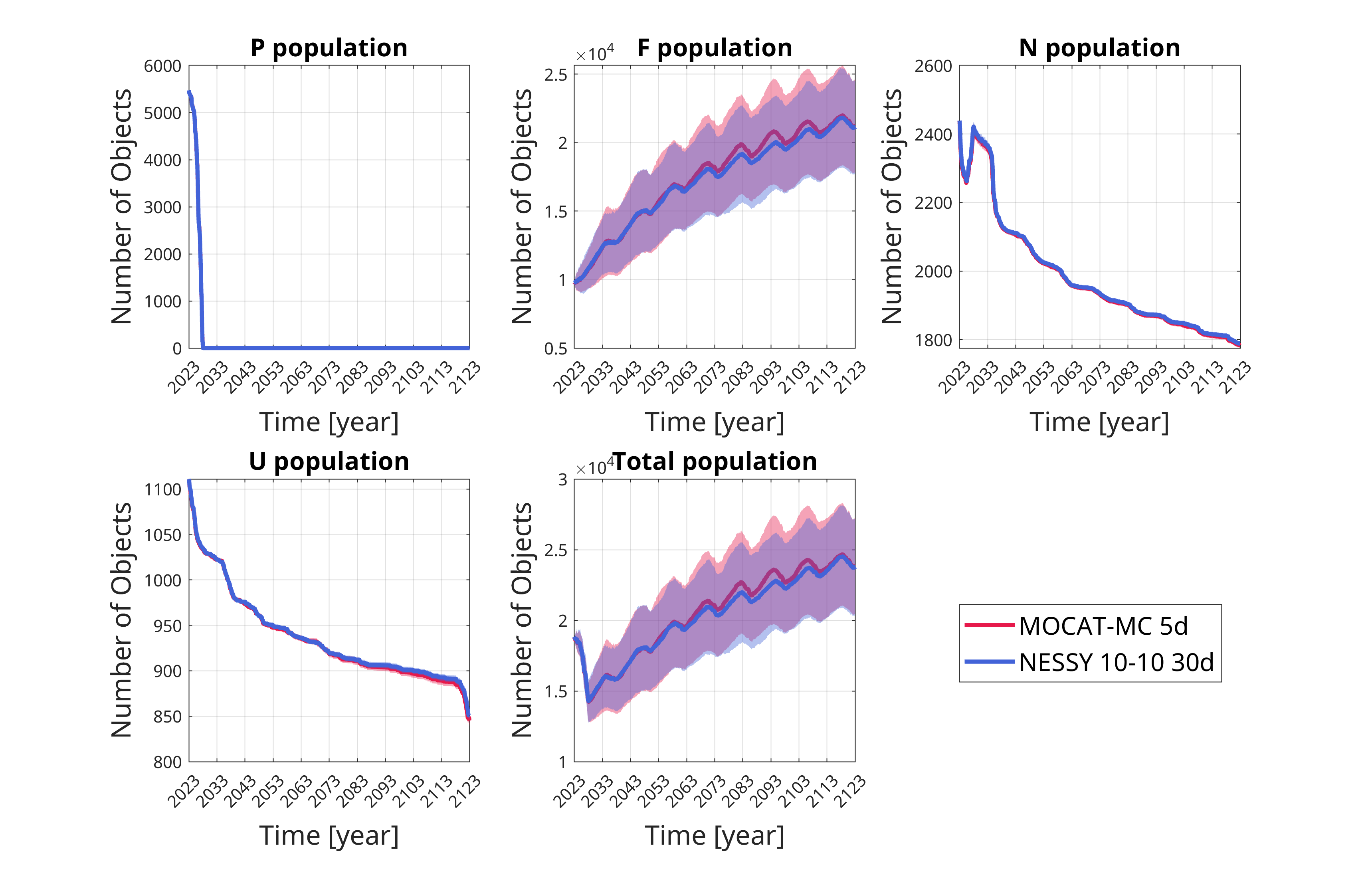}
    \caption{Evolution of the average population propagated with NESSY with 10 deg inclination bins, 10 km altitude shells and a 30 day time step compared with MOCAT-MC with 5 day timestep. The average trend (over 200 Monte Carlo runs) is plotted together with the 1-sigma standard deviation for both propagations.}
    \label{fig:nessy_ev}
\end{figure}

Figure \ref{fig:nessy_ev} compares the evolution of the four species of objects plus the whole population over 100 years. Each simulation was repeated 200 times, and the mean values are shown alongside the 1-sigma standard deviations. The overall trend of the different object populations returned by NESSY closely matches the one of MOCAT-MC. The total number of objects in the final population after 100 years is also similar, differing by a few hundred objects, primarily fragments. Oscillations in the evolution of all the species are due to the Jacchia-Bowman atmospheric model, which is characterised by a value of the atmospheric density that follows the solar cycle. The entire payload population disappears within the first 5 years of the simulation, because of small collisions and PMD. This sharp decrease in $P$ corresponds to an initial increase in non-manoeuvrable satellites in nodes $N$, primarily due to PMD failures. The population of fragments $F$ consistently increases because of collision events. 

\begin{figure*}[h!]
\centering
    \centering
    \begin{subfigure}{0.49\textwidth}
    \includegraphics[width=1\textwidth]{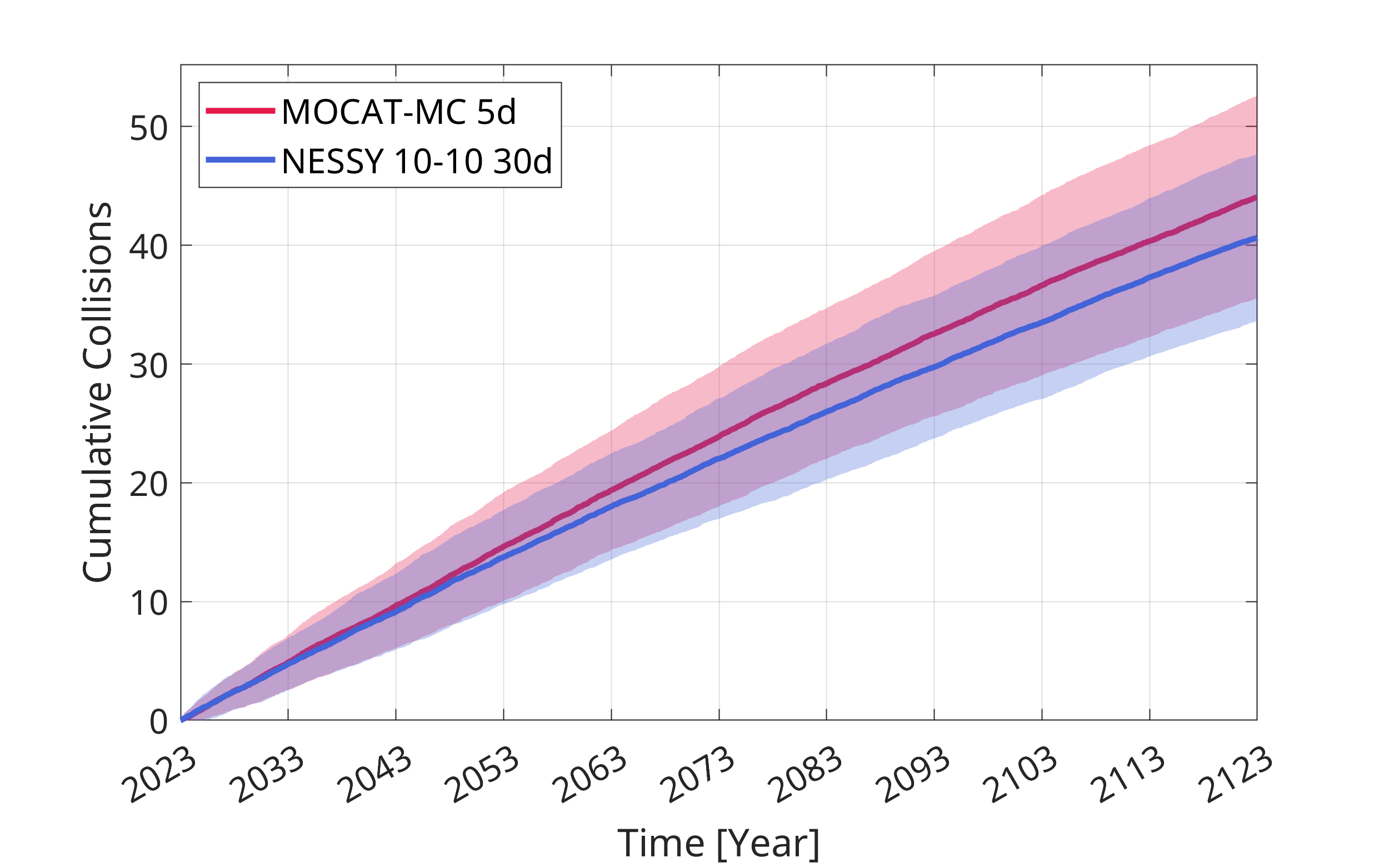}  \caption{\label{fig:pop_nodes}Cumulative number of catastrophic collisions.}
    \label{fig:nessy_cat_coll}
    \end{subfigure}
    \begin{subfigure}{0.49\textwidth}
    \centering
    \includegraphics[width=1\textwidth]{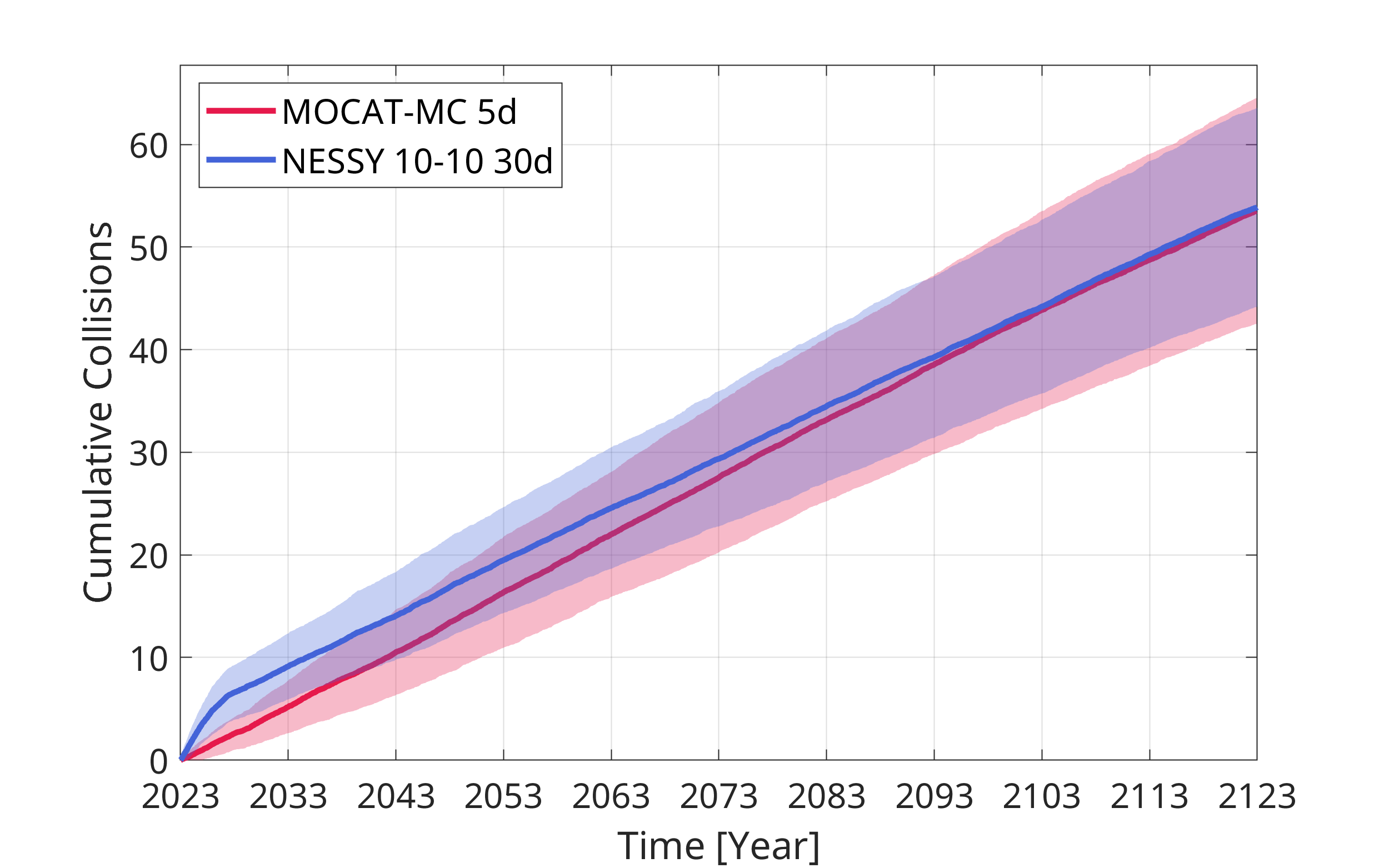}
    \caption{\label{fig:nn_nodes}Cumulative number of total collisions.}
    \label{fig:nessy_nocat_coll}
    \end{subfigure}
    \caption{Cumulative collisions comparison between NESSY with 10 deg inclination bins, 10 km altitude shells and 30 day time step, and MOCAT-MC with 5 days timestep. The average trend  (over 200 Monte Carlo runs) is plotted together with the 1-sigma standard deviation for both propagations.}
\end{figure*}

Figure \ref{fig:nessy_cat_coll} shows the cumulative number of catastrophic collisions over time. The simulations with MOCAT-MC result, on average, in approximately 3 more catastrophic collisions than NESSY over 100 years. Figure~\ref{fig:nessy_nocat_coll} shows the cumulative number of all collisions, including catastrophic, non-catastrophic, and small-scale collisions. The initial difference between the NESSY and MOCAT-MC arises from the approximated treatment of collisions between payloads and small fragments in NESSY. This discrepancy is evident during the early phase of the simulation, while payloads are still in orbit. After their removal, the number of cumulative collisions, in both NESSY and MOCAT-MC, shows a similar trend. It is also worth noting that NESSY is not tracking small fragments and their impact on defunct satellites. A preliminary implementation of a density model for the treatment of small fragments was tested in \cite{PietroSPS} to investigate the environmental impact of solar power satellites and will be introduced in NESSY in future work.

\begin{figure*}[htbp]
    \centering
    \begin{subfigure}{0.9\textwidth}
        \includegraphics[width=1.03\textwidth]{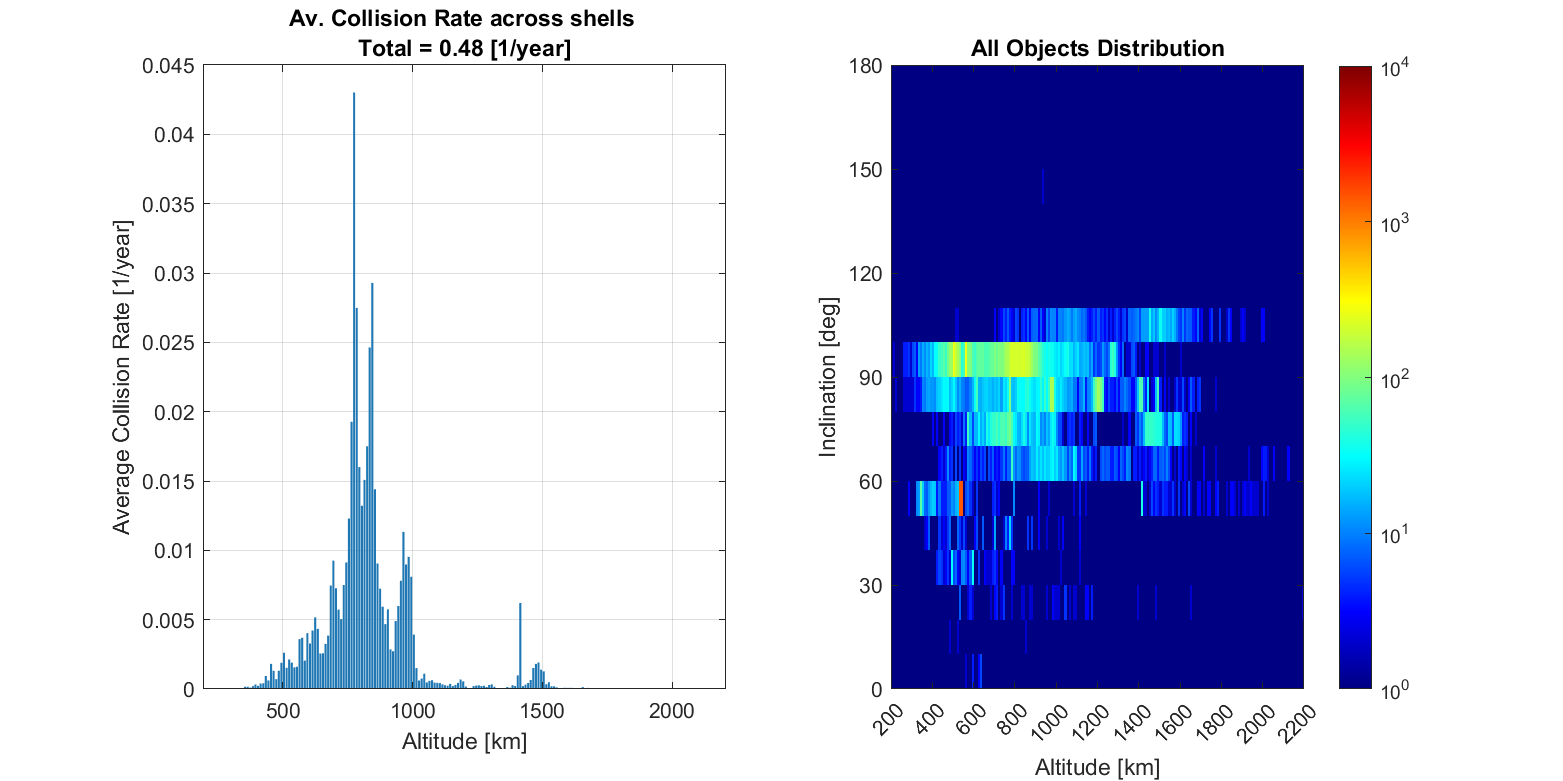}  
        \caption{\label{fig:cr_23}Collision rate distribution by shell and entire population distribution through nodes (averaged over 200 Monte Carlo runs).}
    \end{subfigure}
    \vspace{0.5em} 
    \begin{subfigure}{0.9\textwidth}
        \includegraphics[width=1.03\textwidth]{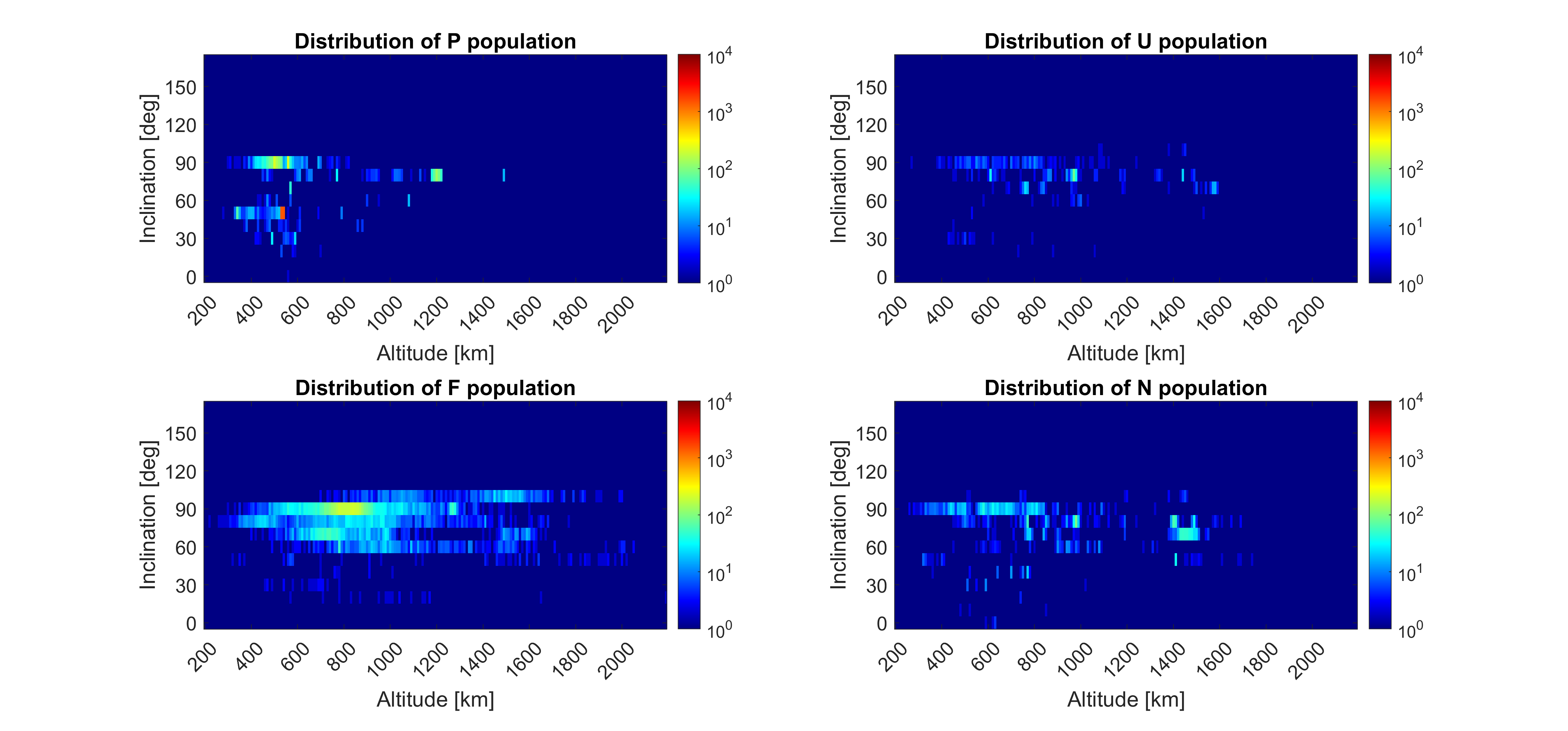}
        \caption{\label{fig:dens_23}Objects distribution by class through nodes (averaged over 200 Monte Carlo runs).}
    \end{subfigure}
    
    \caption{Distribution of the collision rate and number of objects by class in year 2023 (beginning of the propagation).}
\end{figure*}

\begin{figure*}[htbp]
    \centering
    \begin{subfigure}{0.9\textwidth}
        \includegraphics[width=1.03\textwidth]{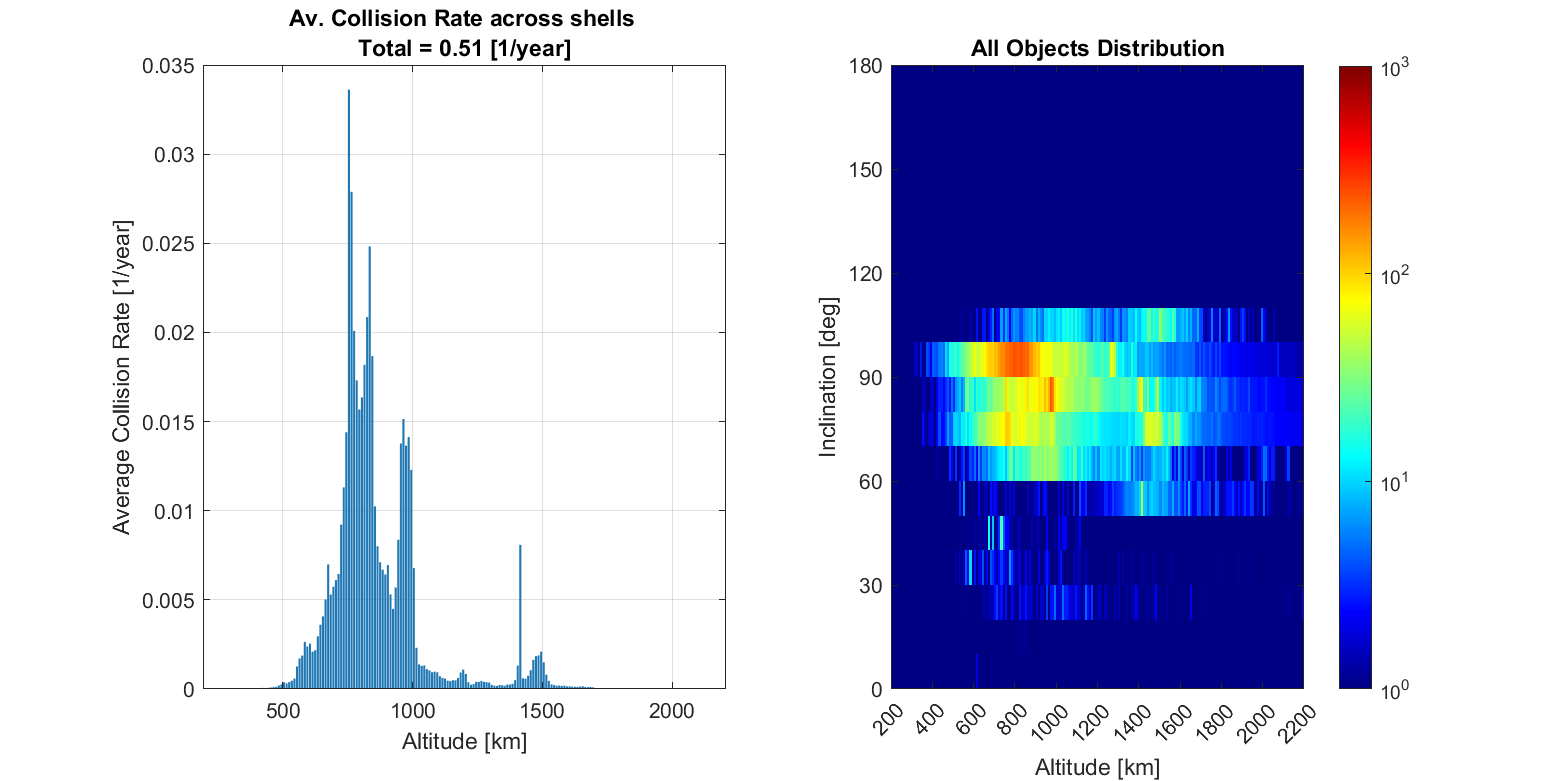}  
        \caption{\label{fig:cr_73}Collision rate distribution by shell and entire population distribution through nodes (averaged over 200 Monte Carlo runs).}
    \end{subfigure}
    \vspace{0.5em} 
    \begin{subfigure}{0.9\textwidth}
        \includegraphics[width=1.03\textwidth]{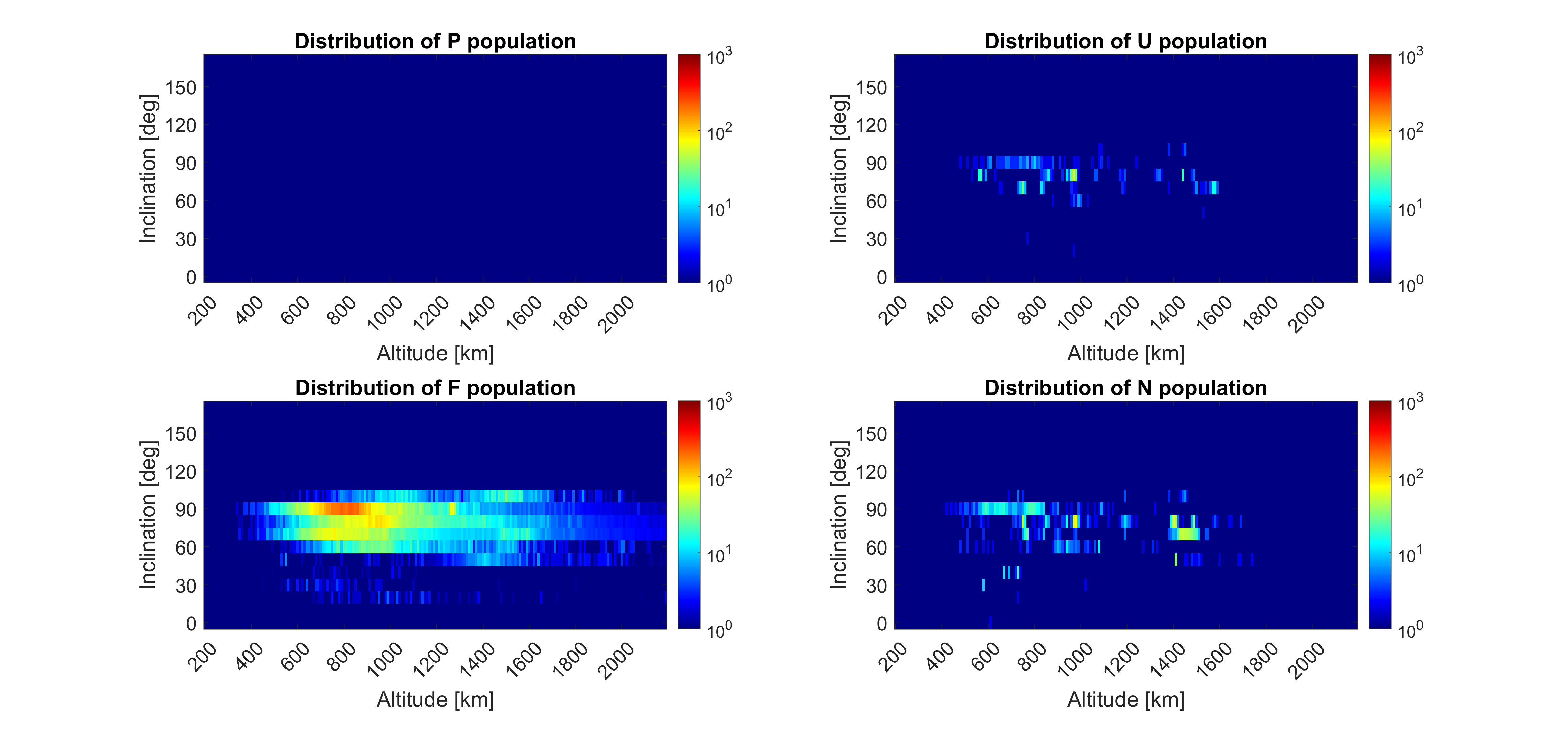}
        \caption{\label{fig:dens_73}Objects distribution by class through nodes (averaged over 200 Monte Carlo runs).}
    \end{subfigure}
    
    \caption{Distribution of the collision rate and number of objects by class in year 2073 (halftime of the propagation).}
\end{figure*}

\begin{figure*}[htbp]
    \centering
    \begin{subfigure}{0.9\textwidth}
        \includegraphics[width=1.03\textwidth]{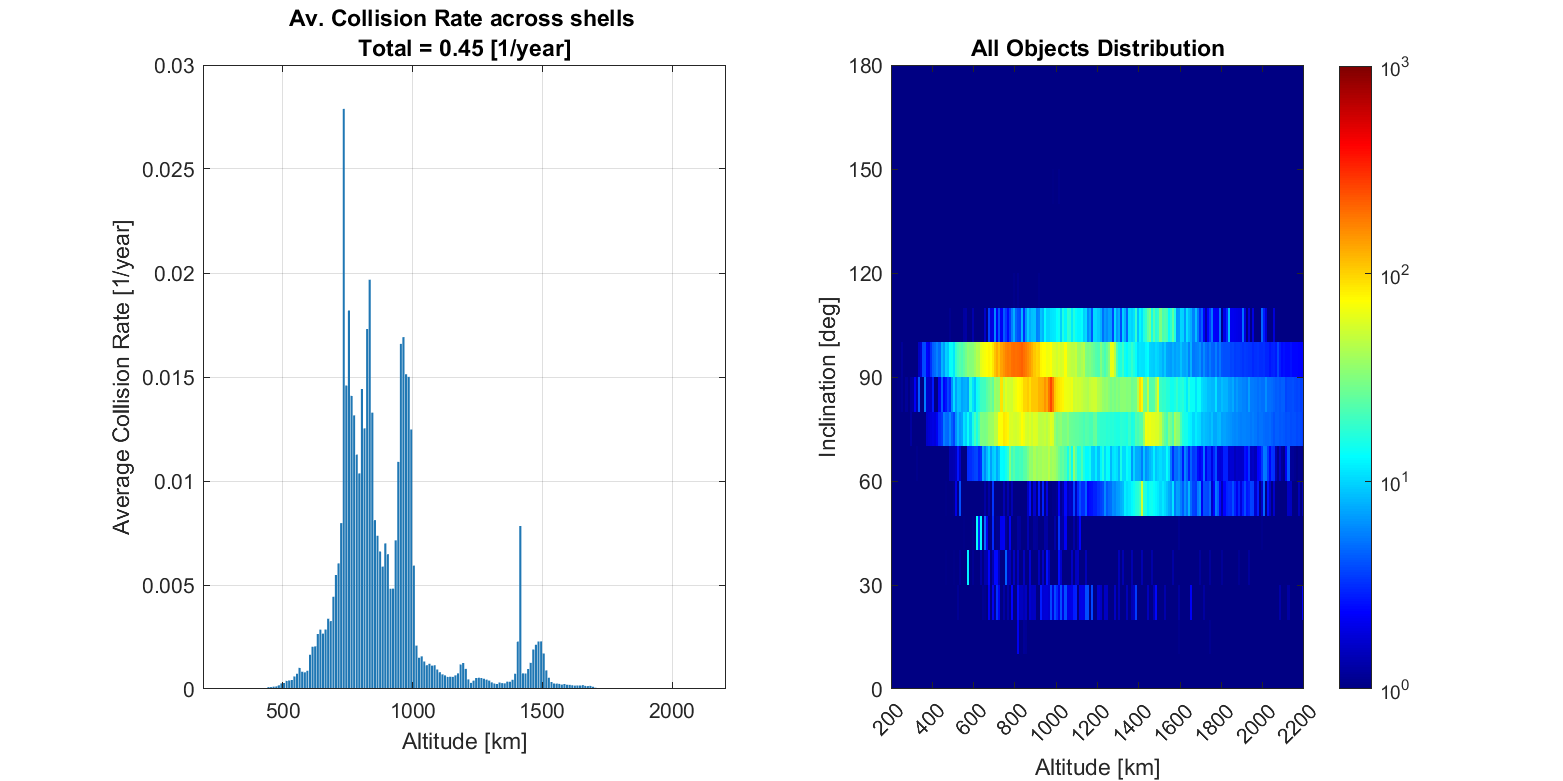}  
        \caption{\label{fig:cr_123}Collision rate distribution by shell and entire population distribution through nodes (averaged over 200 Monte Carlo runs).}
    \end{subfigure}
    \vspace{0.5em} 
    \begin{subfigure}{0.9\textwidth}
        \includegraphics[width=1.03\textwidth]{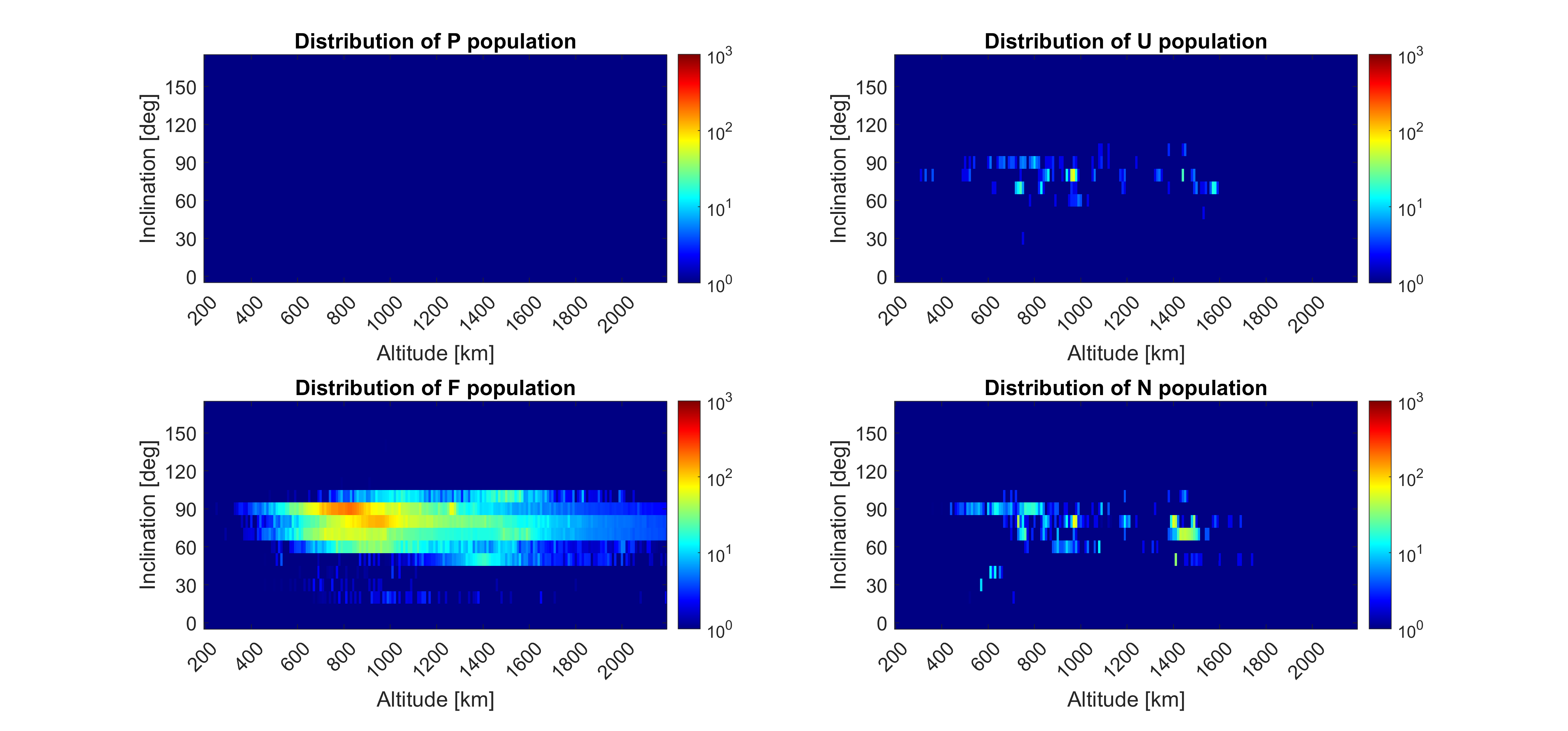}
        \caption{\label{fig:dens_123}Objects distribution by class through nodes (averaged over 200 Monte Carlo runs).}
    \end{subfigure}
    
    \caption{Distribution of the collision rate and number of objects by class in year 2123 (end of the propagation).}
\end{figure*}

Figures \ref{fig:cr_23}, \ref{fig:dens_23}, \ref{fig:cr_73}, \ref{fig:dens_73}, \ref{fig:cr_123}, \ref{fig:dens_123} show the estimated annual collision rate (catastrophic collisions, non-catastrophic collisions and small collisions) per 10 km and spatial distribution of objects, respectively, at 3 instants of time along the 100 year propagation (averaged over 200 Monte Carlo runs): year 2023, 2073 and 2123. These figures illustrate how the collision rate evolves in response to changes in the spatial distribution of objects. As it is clear from the plots, the main contribution to the collision rate comes from the fragment population.

\subsection{Sensitivity Analysis}
\begin{figure}[h!]
    \centering
    \includegraphics[width=0.9\textwidth]{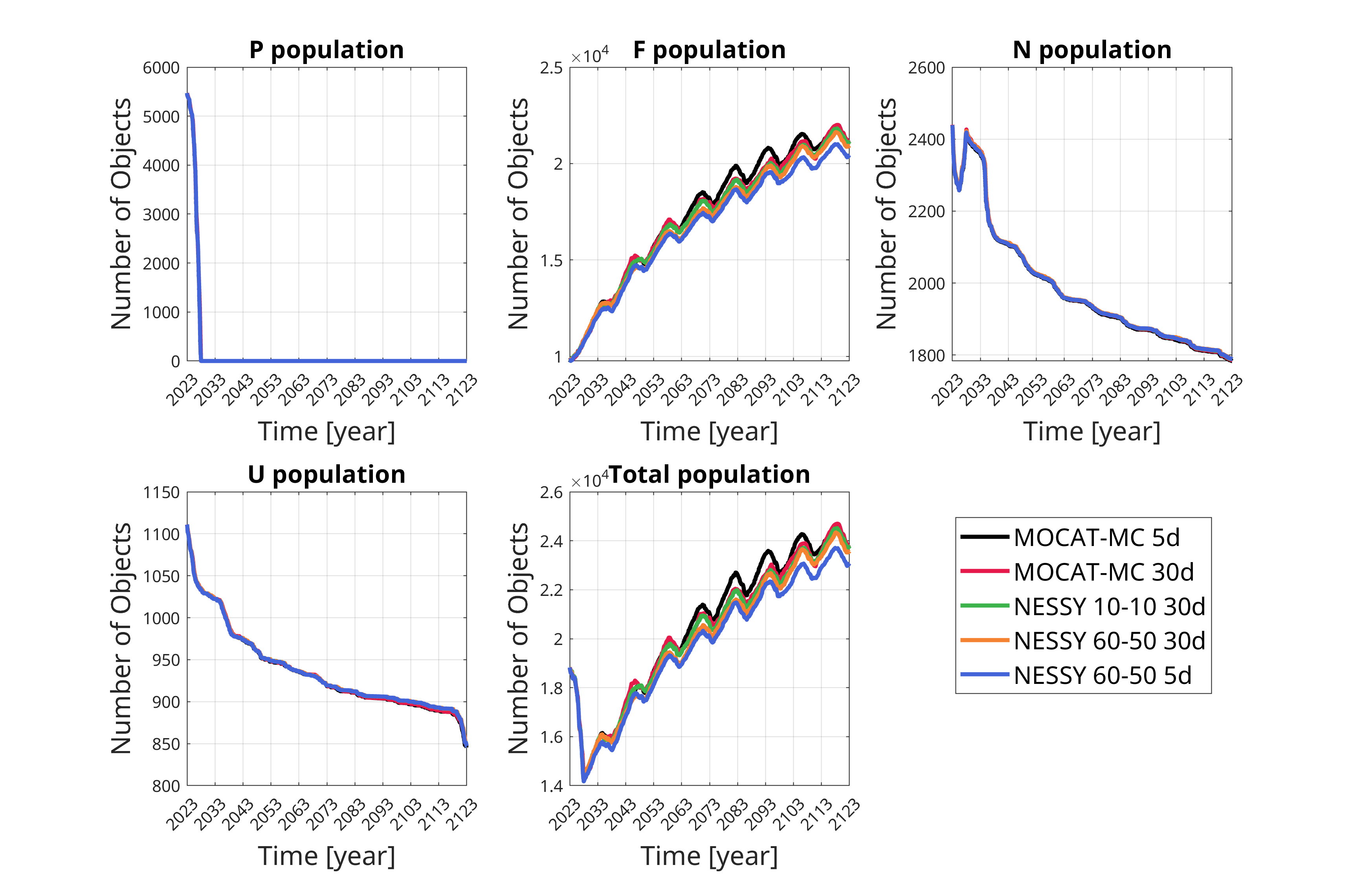}
    \caption{Evolution of the environment propagated with NESSY and MOCAT-MC for different network sizes and time steps (averaged over 200 Monte Carlo runs).}
    \label{fig:comparison_nessy_all}
\end{figure}

\begin{figure*}[h!]
\centering
    \centering
    \begin{subfigure}{0.49\textwidth}
    \includegraphics[width=1\textwidth]{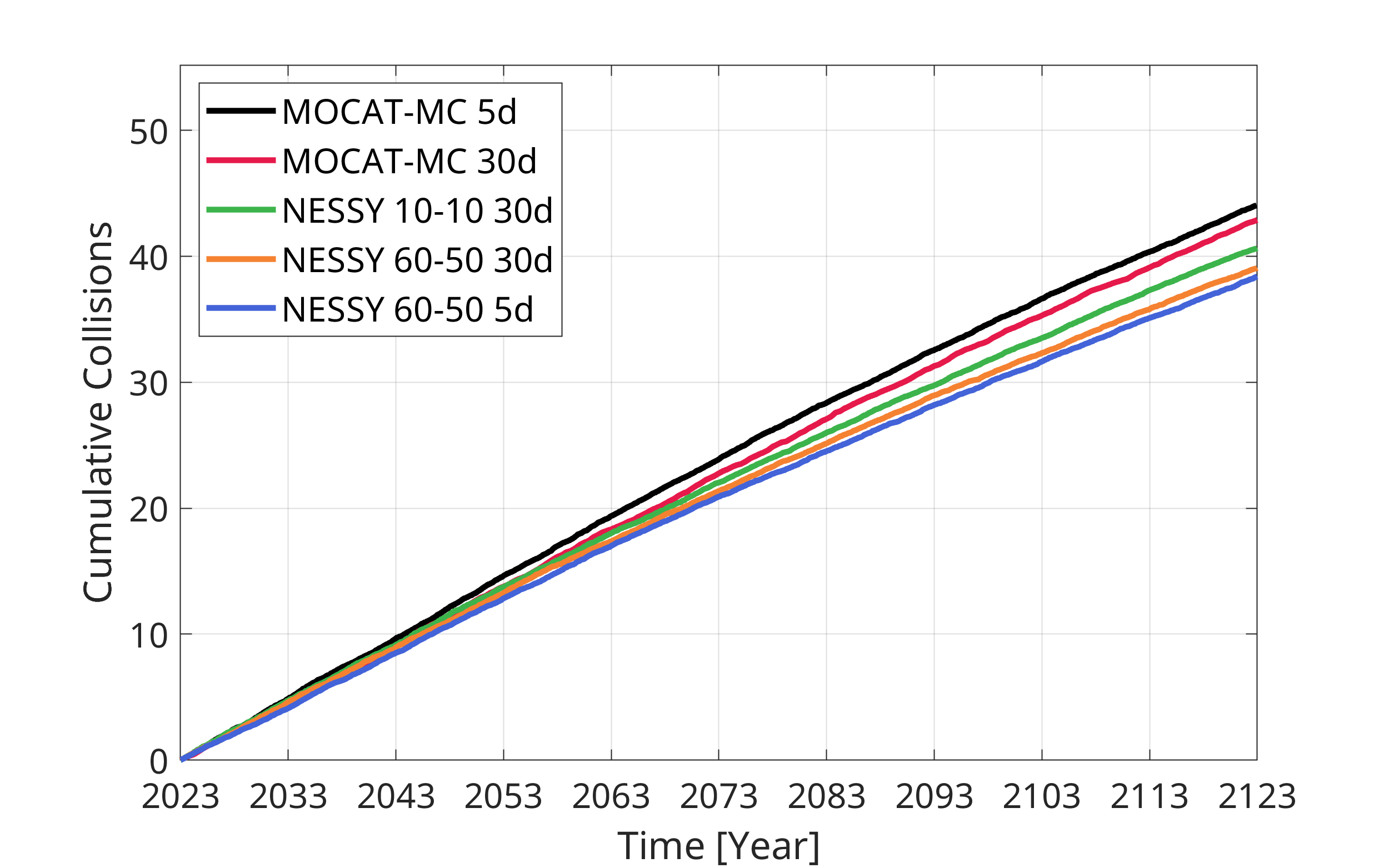}  \caption{\label{fig:pop_nodes}Cumulative number of catastrophic collisions.}
    \end{subfigure}
    \begin{subfigure}{0.49\textwidth}
    \centering
    \includegraphics[width=1\textwidth]{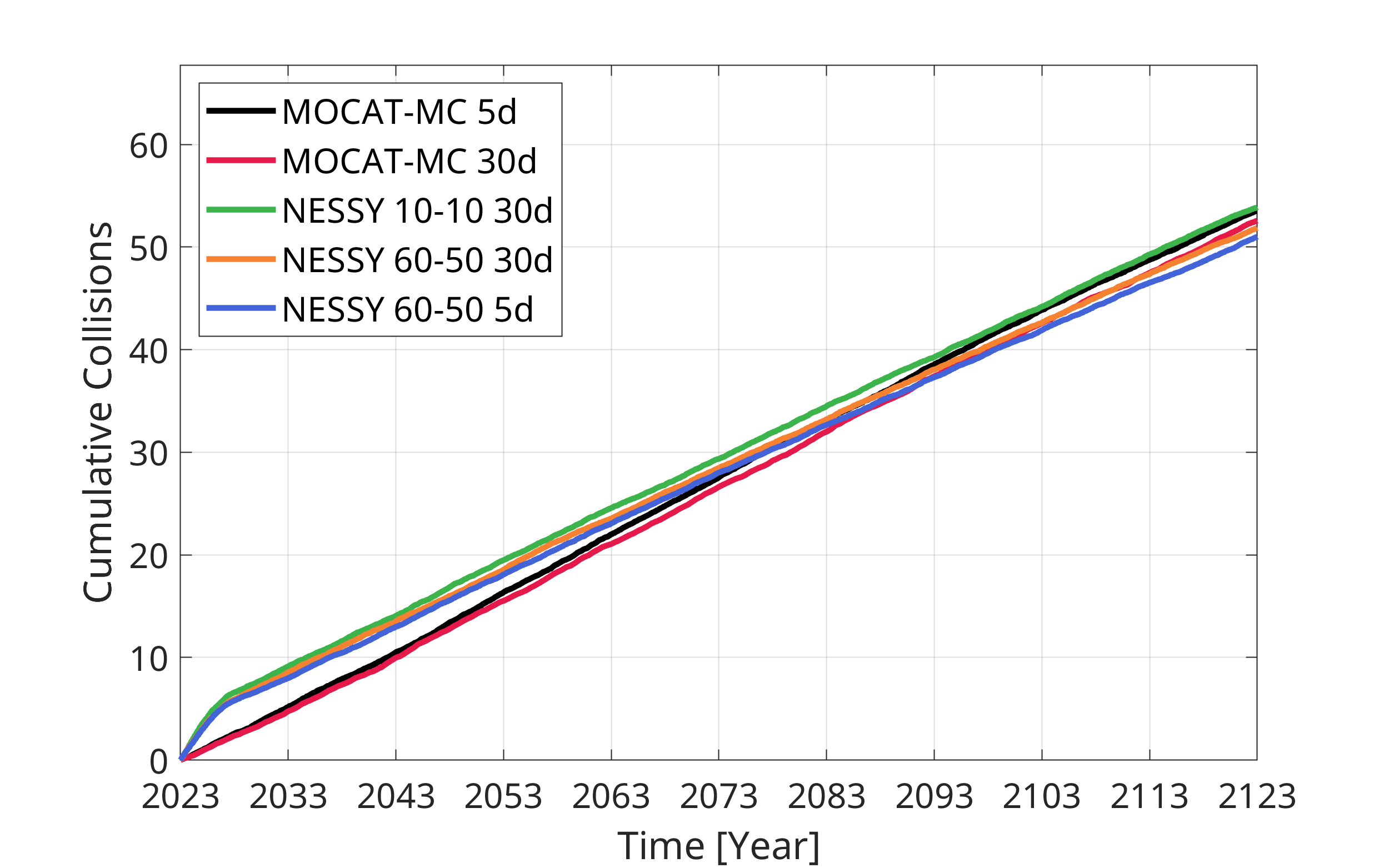}
    \caption{\label{fig:nn_nodes}Cumulative number of total collisions.}
    \end{subfigure}
    \caption{Comparison of the cumulative number of collisions between NESSY and MOCAT-MC for different network sizes and time steps (averaged over 200 Monte Carlo runs).\label{fig:cum_col_comp}}
\end{figure*}

In this section, we present an analysis of the sensitivity of the simulation results produced by NESSY to inclination bin size, time step $\Delta t_k$, and orbit shell size. The comparison of the evolution of each object species is reported in Figure \ref{fig:comparison_nessy_all} while the cumulative number of collisions is represented Figure \ref{fig:cum_col_comp}. 
We tested two time steps, 30-day and 5-day, for both NESSY and MOCAT-MC to assess the sensitivity to the discretisation in time. We tested multiple combinations of shell and inclination bin sizes, but in these figures we report only the following two combinations: 10 deg and 10 km and 60 deg and 50 km. 
The results in the figures suggest that, as the size of inclination bins and orbit shells decreases, the number of objects in the final population and a cumulative number of collisions of NESSY progressively converge toward those of MOCAT-MC. Moreover, the results obtained with a 30-day time step are highly consistent with those of the 5-day case, both in terms of population evolution and cumulative number of collisions.

Figure \ref{fig:pop_nodes}
shows the total number of objects in the final population after 100 years for different combinations of altitude shell and inclination bin size. In the figure, the first number next to each blue dot is the size of the inclination bin and the second number is the size of the altitude shell. The red line is a least square fit to the average final population, over 200 runs, for each combination of altitude shell and inclination bin size. The figure shows that as the number of nodes increases, the final population initially grows and then stabilises. Specifically, the final population size increases with the number of nodes up to approximately 1,000 nodes. Beyond this point, the total number of objects in the final population remains relatively stable, with minor fluctuations around 23,800 objects. All simulations were performed in MATLAB on an Ubuntu 24.04.2 system equipped with an Intel Xeon w9-3495X ×112 processor. The average computational time to propagate Eqs. (\ref{eq: nodes dynamics}) for a single node for 100 years is approximately 0.2 seconds. It is worth noting that the current implementation is not yet computationally optimal, and computational efficiency will be improved in future work.

\begin{figure}[h!]
    \centering
    \includegraphics[width=0.6\textwidth]{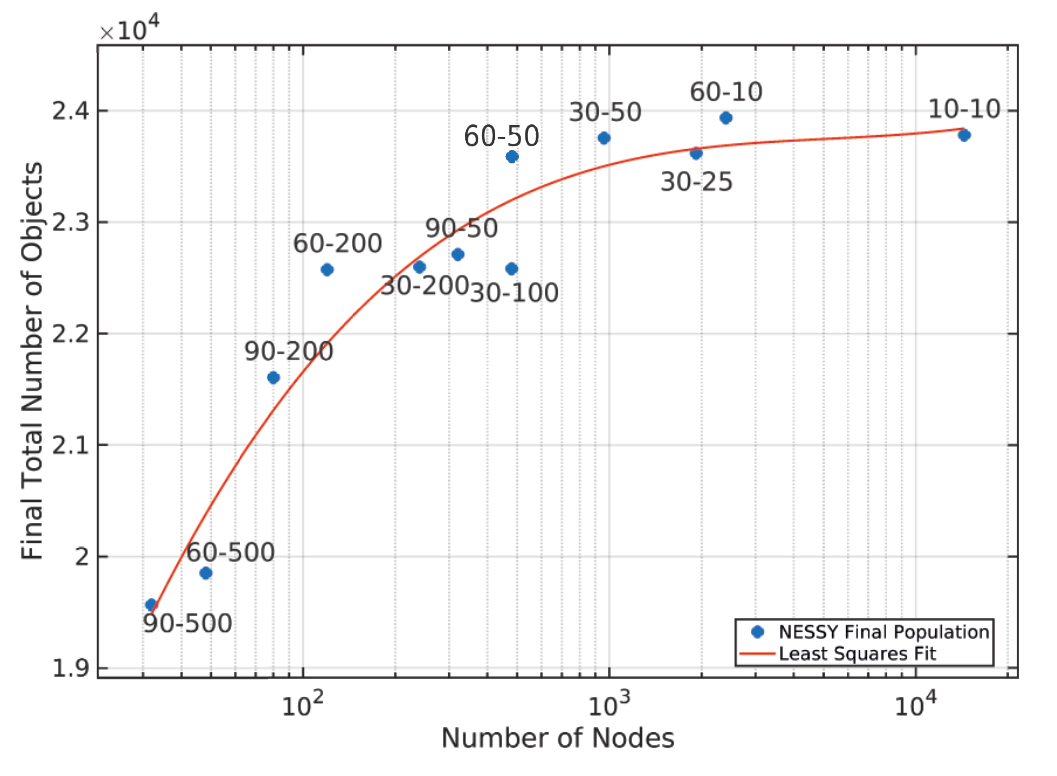}
    \caption{Average number of total objects in the population after 100 years 
 for different network discretisations (averaged over 200 Monte Carlo runs).}
    \label{fig:pop_nodes}
\end{figure}

\newpage

\section{Case Studies}\label{Sec:Assessment of the space environment stability}

This section demonstrates how NESSY can be applied to evaluate the impact of launch traffic and debris mitigation policies on the long-term evolution of the space environment. The initial population is identical to that used in previous sections (see Table~\ref{table:params_valid}), and the launch traffic model corresponds to the one introduced in Section~\ref{sec:Launch Traffic Model}. All simulations in this section were conducted using a discretisation where each node represents a 60 deg inclination bin and a 50 km altitude shell.

\subsection{Impacts of Future Launch Traffic}
\begin{figure}[H]
    \centering
    \includegraphics[width=0.7\textwidth]{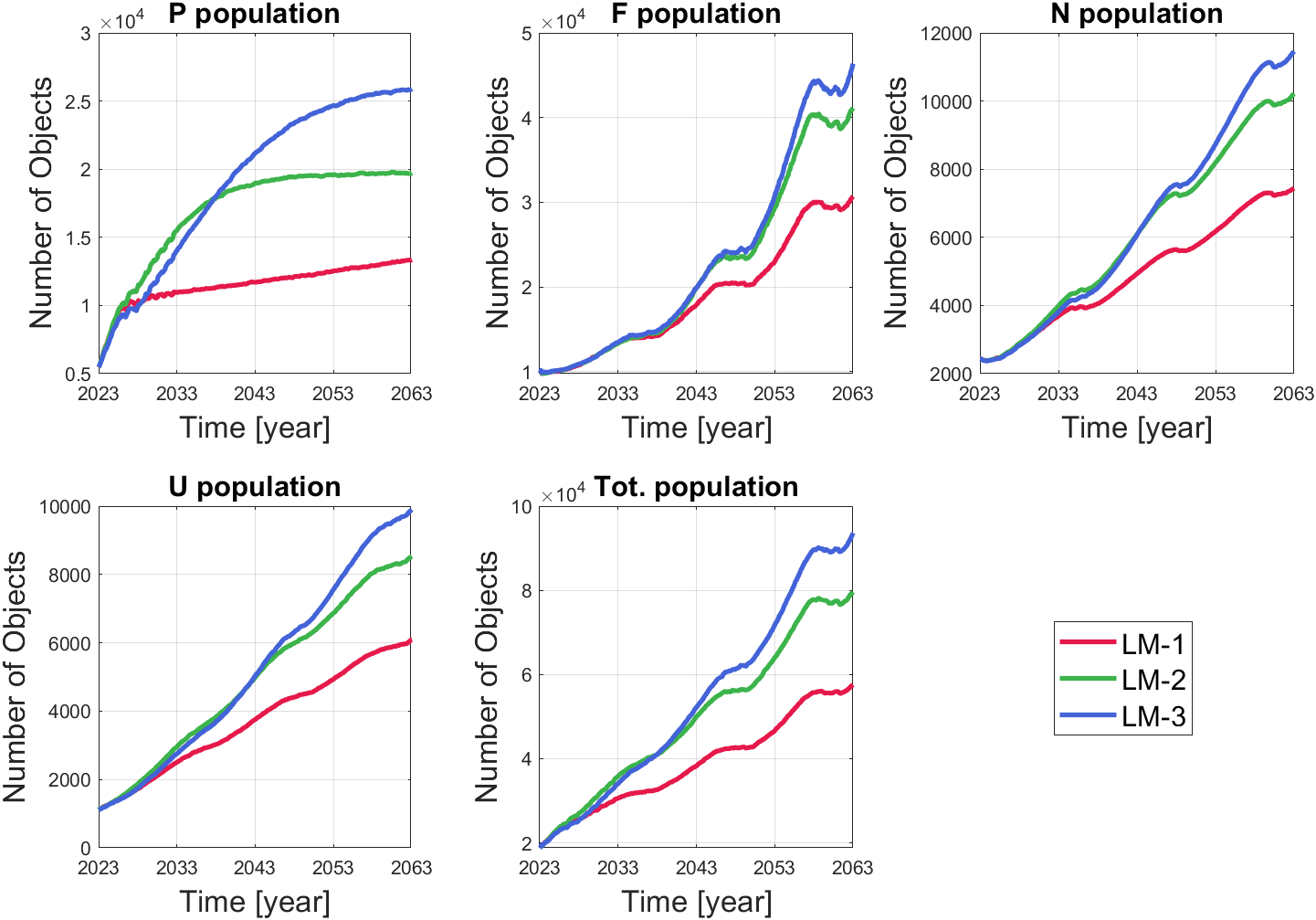}
    \caption{Evolution of the environment based on 3 different launch traffic forecasts.}
    \label{fig:comp_lm}
\end{figure}

\begin{figure}[H]
    \centering
    \includegraphics[width=0.5\textwidth]{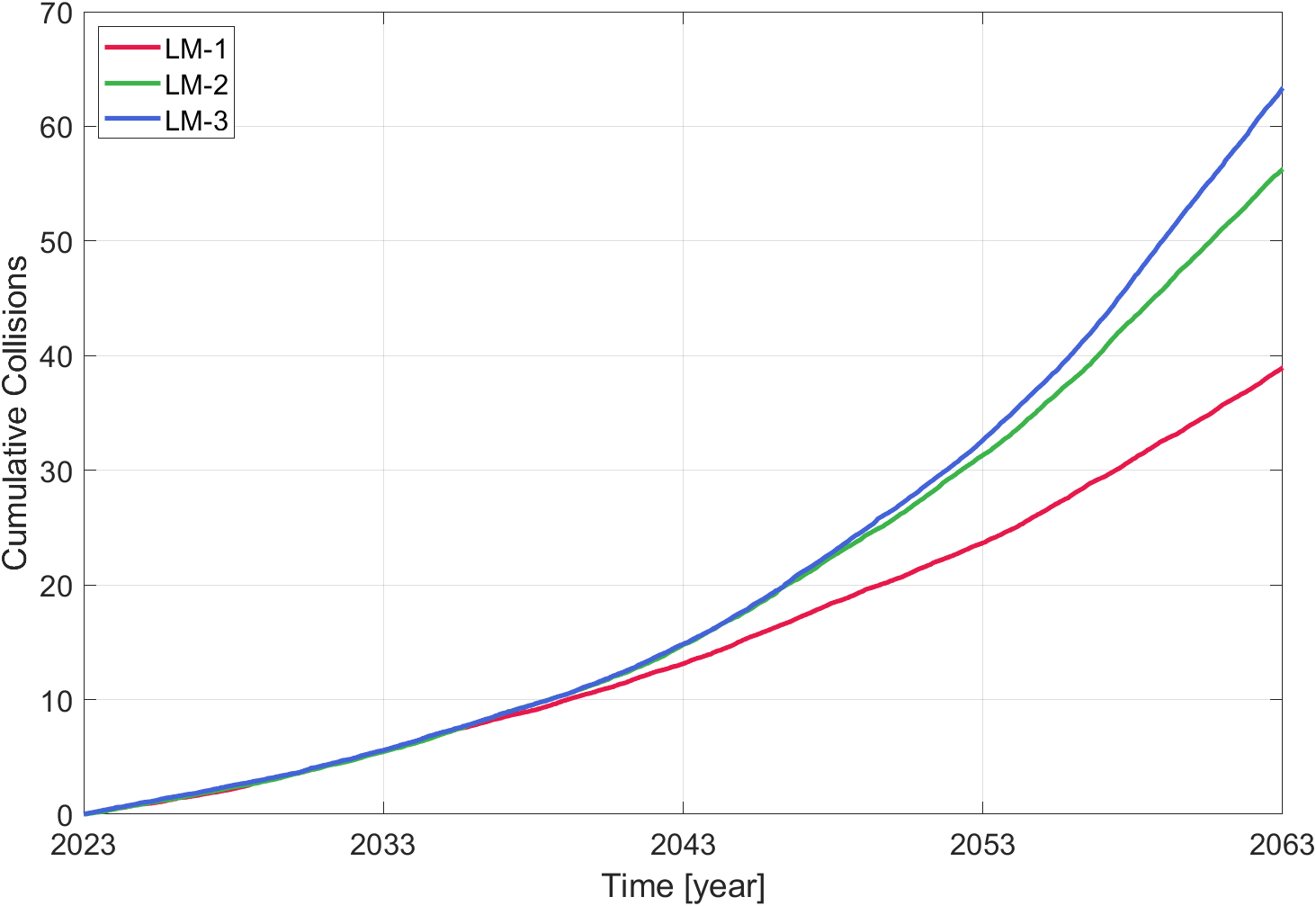}
    \caption{Cumulative number of catastrophic collisions based on 3 different launch traffic forecasts.}
    \label{fig:comp_lm_coll}
\end{figure}

 Figures \ref{fig:comp_lm} and \ref{fig:comp_lm_coll} compare the evolution of different species and the cumulative collisions over 40 years, for the three launch traffic forecasts (LM) introduced in Section \ref{sec:Launch Traffic Model}. Compared to the baseline case, the launch traffic introduces a significant increase in all species, with the increase in the number of collisions changing from linear to exponential. Since the launch model presented in this paper incorporate the addition of upper stages and non-manoeuvrable satellites, all three forecasts lead to an increase in the number of upper stages and non-manoeuvrable satellites, while the number of payloads eventually reaches either a plateau or continues with a linear increase. 

The time variation of the number of payloads can be further explained by computing the expected increment in the number of payloads computed with equations \eqref{eq: nodes dynamics}. Assuming that there is no decay of payloads due to atmospheric drag, the expected value of $\Delta x_{P_i}$ at time $t_k$ is
\begin{equation}
\mathbb{E}(\Delta x_{P_i}) = - \sum_j (1-s_{CAM})\left(\tau_{P_iU_j}\Delta t_k + \tau_{P_iN_j}\Delta t_k + \tau_{P_iF_j}\Delta t_k + (1-s_{CAM})\tau_{P_iP_j}\Delta t_k\right) - \sum_j \kappa \tau_{P_iF_j}\Delta t_k  - \Delta^{EOL}_{P_i} + \Lambda_{P_i},
\end{equation}
then the evolution equation of payloads can be simplified as
\begin{equation}
    \dot{x}_{P_i} = -\hat{\gamma}x_{P_i} + \lambda,
\end{equation}
where $\hat{\gamma}$ indicates the removal rate per payload due to PMD and collisions with other species. 
\begin{equation}\label{eq:payload_mean_eq}
    \hat{\gamma} = \frac{\sum_j (1-s_{CAM})\left(\tau_{P_iU_j} + \tau_{P_iN_j} + \tau_{P_iF_j} + (1-s_{CAM})\tau_{P_iP_j}\right) + \sum_j \kappa \tau_{P_iF_j} }{x_{P_i}} + \frac{\Delta^{EOL}_{P_i}}{\Delta t_k x_{P_i}},
\end{equation}
where the second term in $\hat{\gamma}$, representing the removal rate due to PMD, remains constant for a fixed operational lifetime. Meanwhile, $\lambda$ indicates the launch rate, given by $\lambda = \Lambda_{Pi} / \Delta t_k$. The analytical solution for the payload dynamics is therefore
\begin{equation}\label{eq:payoad_analytical}
    x_{P_i}(t) = \frac{\lambda}{\hat{\gamma}} - C e^{-\hat{\gamma} t}.
\end{equation}
This solution suggests that the payload population asymptotically approaches a steady-state value of $\lambda / \hat{\gamma}$, where $C$ represents the difference between the initial value and the steady-state value. In the baseline scenario, where the launch rate $\lambda$ is zero, the steady-state solution reduces to zero. Given that the successful rate of collision avoidance manoeuvrer is very high, with $s_{CAM} = 99.99\%$ as suggested in the previous section, the first term in the expression for $\hat{\gamma}$, which accounts for collisions, contributes less to the removal rate than the second term, which accounts for PMD. Consequently, it is reasonable to treat $\hat{\gamma}$ as a constant, and thus, the steady-state average number of payloads, $x_{P_i}$, becomes directly proportional to the launch rate $\lambda$. However, despite the high $s_{CAM}$ assumed in this paper, the collision rate among payloads and other objects is likely to increase as the number of objects grows over time. Consequently, $\hat{\gamma}$ will gradually increase, leading to a slightly lower steady-state value for $x_{P_i}$.

As shown in Figure \ref{fig:Model of number of objects launched to LEO}, both launch traffic forecast 2 (LM-2) and launch traffic forecast 3 (LM-3) show nearly constant launch rates $\lambda$ after around 2050, while the launch traffic forecast 1 (LM-1) shows a linear increase of launch rate $\lambda$ from 2023. Although payloads are assumed to perform CAM with a high success rate $s_{CAM} = 99.99\%$, launch-related objects, such as upper stages and mission-related objects, continue to contribute to the accumulation of non-manoeuvrable objects. Notably, the oscillations in the number of non-manoeuvrable satellites are slightly stronger than those in the number of upper stages, with the oscillation closely following those of the Jacchia-Bowman 2008 (JB2008) atmospheric density model over the years. This behaviour is attributed to the higher area-to-mass ratio of non-manoeuvrable satellites, which makes them more susceptible to atmospheric drag. Among the three launch traffic forecasts, LM-3 results in the largest growth across all object classes. LM-3 has the most significant impact on each species. Specifically, after 40 years, the number of payloads reaches approximately 25,500. Considering the introduction of mission-related objects and upper stages by the launch activities, the number of non-manoeuvrable satellites increases from around 2,450 to 11,350, and the number of upper stages increases from about 1,111 to 9,950. The number of fragments grows from 9,804 to about 46,140. The average cumulative number of catastrophic collisions reaches 63.2 in LM-3 over 40 years, compared to only about 40 over 100 years in a no-launch scenario.

\subsection{Impacts of Debris Mitigation Actions}

This section investigates the impact of debris mitigation actions, such as the Collision Avoidance Manoeuvres (CAM) and the Post-Mission Disposal (PMD), on the space environment. 

\begin{figure*}[htbp]
\centering
    \centering
    \includegraphics[width=0.7\textwidth]{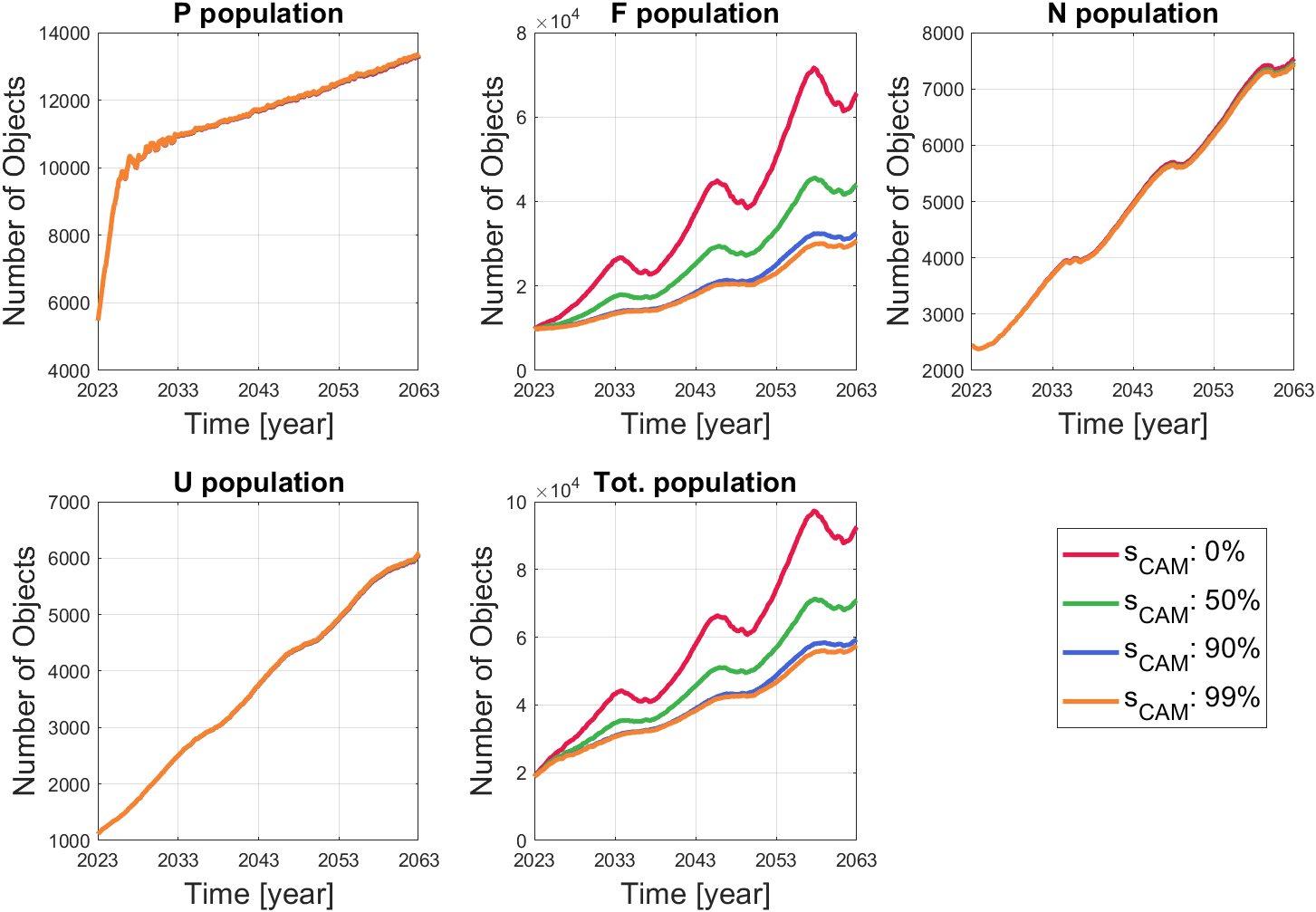}
    \caption{Evolution of the environment based on different $s_{CAM}$.}. 
    \label{fig:cam_ev}
\end{figure*}

Considering $s_{CAM}$ values of 0\%, 50\%, 90\% and 99.99\% with fixed $\gamma = 5\%$, Figure \ref{fig:cam_ev} compares the evolution of each species under launch traffic model \#1. Overall, simulation results suggest that a decrease in $s_{CAM}$ results in an increase in final population and cumulative catastrophic collisions. Compared to Figure \ref{fig:comp_lm} and Figure \ref{fig:comp_lm_coll}, decreasing $s_{CAM}$ from $99.99\%$ to $0\%$ results in an increase in final population of fragments from about 29,930 to around 64,320, and an increase in cumulative catastrophic collisions from 41 to 260. Specifically, as $s_{CAM}$ decreases, the probability of the collision involving payloads increases, leading to a reduction in the number of upper stages and an increase in fragments. Interestingly, the number of non-manoeuvrable satellites shows a slight increase. This is because the additional fragments generated by the collisions involving payloads further increase the probability of small collisions, that is, the increased likelihood that payloads will be converted into non-manoeuvrable satellites due to the increase in the number of small fragments between 1 cm and 10 cm. Moreover, the stronger oscillation behaviour observed in the case of lower $s_{CAM}$ is caused by the increased number of new fragments concentrated in lower-altitude regions (around 500 km),  where atmospheric drag significantly affects fragment dynamics.

\begin{figure*}[htbp]
\centering
    \centering
    \includegraphics[width=0.5\textwidth]{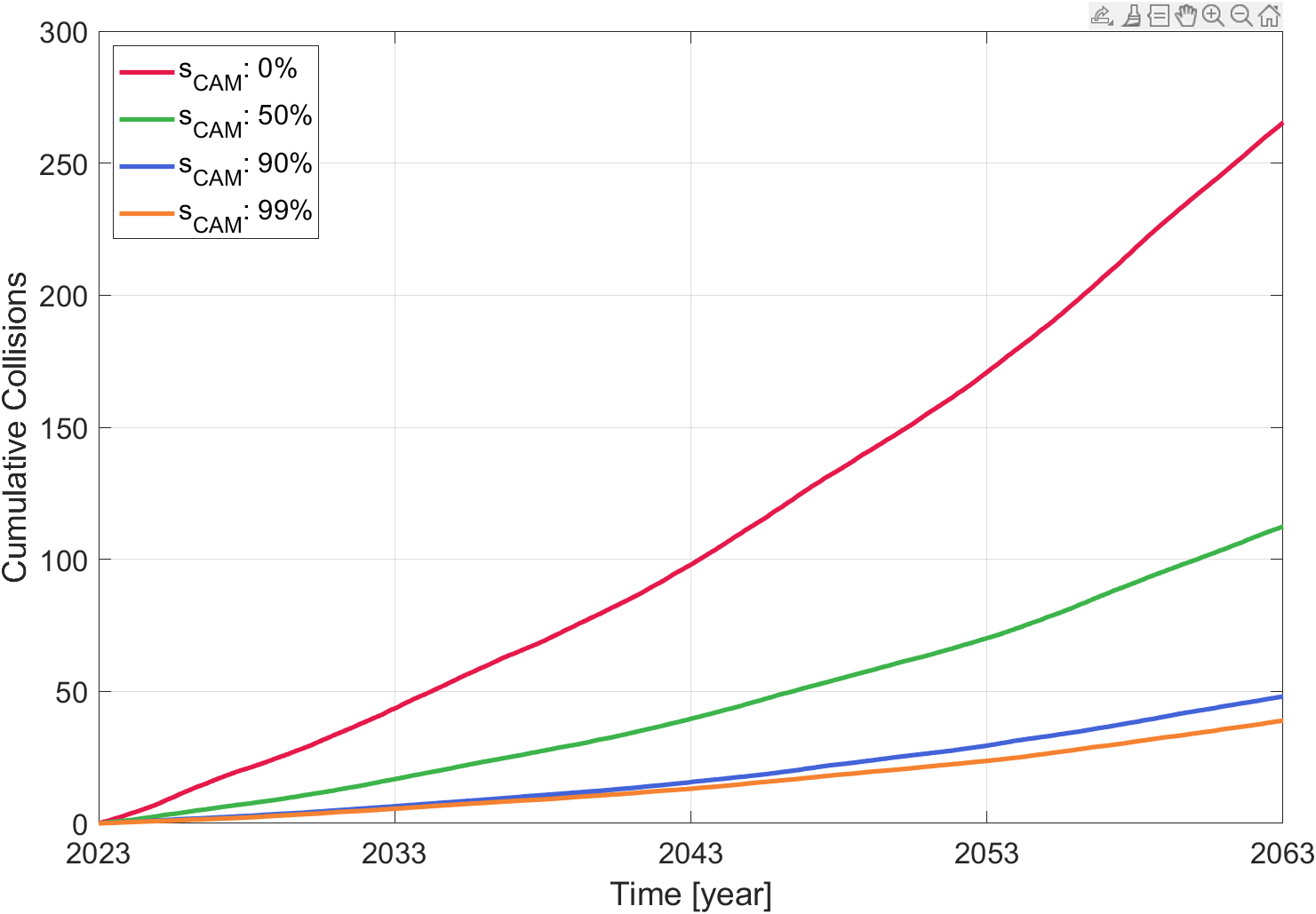}
    \caption{Cumulative number of catastrophic collisions based on different $s_{CAM}$}.
    \label{fig:cam_coll}
\end{figure*}

Considering $\gamma$ values of  0\%, 5\%, 20\%, and 40\% with fixed $s_{CAM} = 99.99\%$, Figure \ref{fig:cam_ev} compares the evolution of each species under launch traffic model \#1. Overall, simulation results suggest that an increase in $\gamma$ results in an increase in final population and cumulative catastrophic collisions. Compared to Figure \ref{fig:comp_lm} and Figure \ref{fig:comp_lm_coll}, increasing $\gamma$ from $0\%$ to $40\%$ results in an increase in the final population of fragments from about 28,130 to around 60,110, and an increase in the cumulative catastrophic collisions from 33 to 164. This trend is primarily driven by the fact that higher $\gamma$ values correspond to more payloads failing to perform PMD, thereby increasing the population of non-manoeuvrable satellites and, consequently, the number of fragments. Consistent with the aforementioned conclusions, the stronger oscillation behaviour observed in the fragment population at higher $\gamma$ is primarily due to the accumulation of new fragments at lower altitudes.

\begin{figure*}[htbp]
\centering
    \centering
    \includegraphics[width=0.7\textwidth]{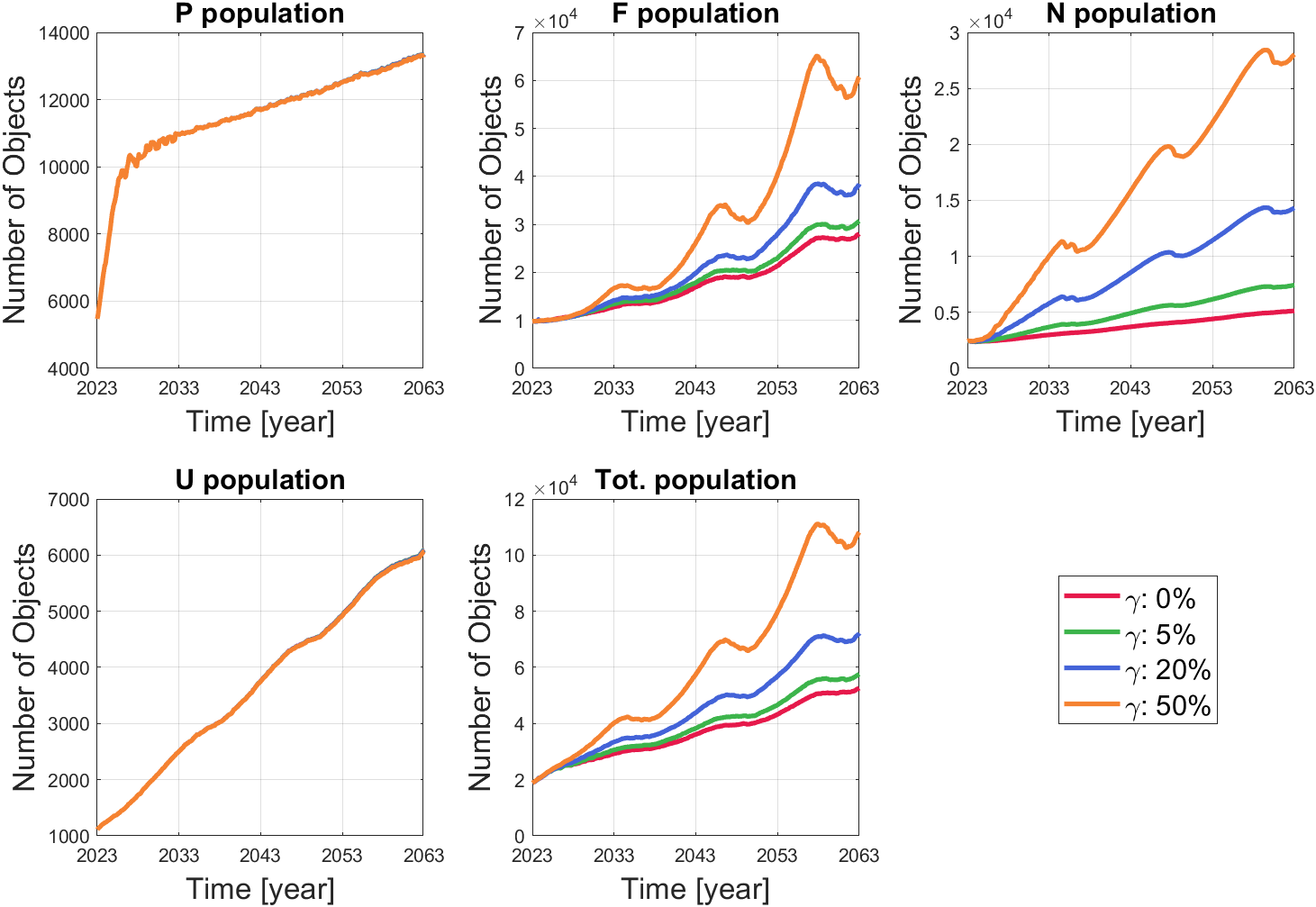}
    \caption{Evolution of the environment based on different $\gamma$.} 
    \label{fig:gamma_ev}
\end{figure*}

\begin{figure*}[htbp]
\centering
    \centering
    \includegraphics[width=0.5\textwidth]{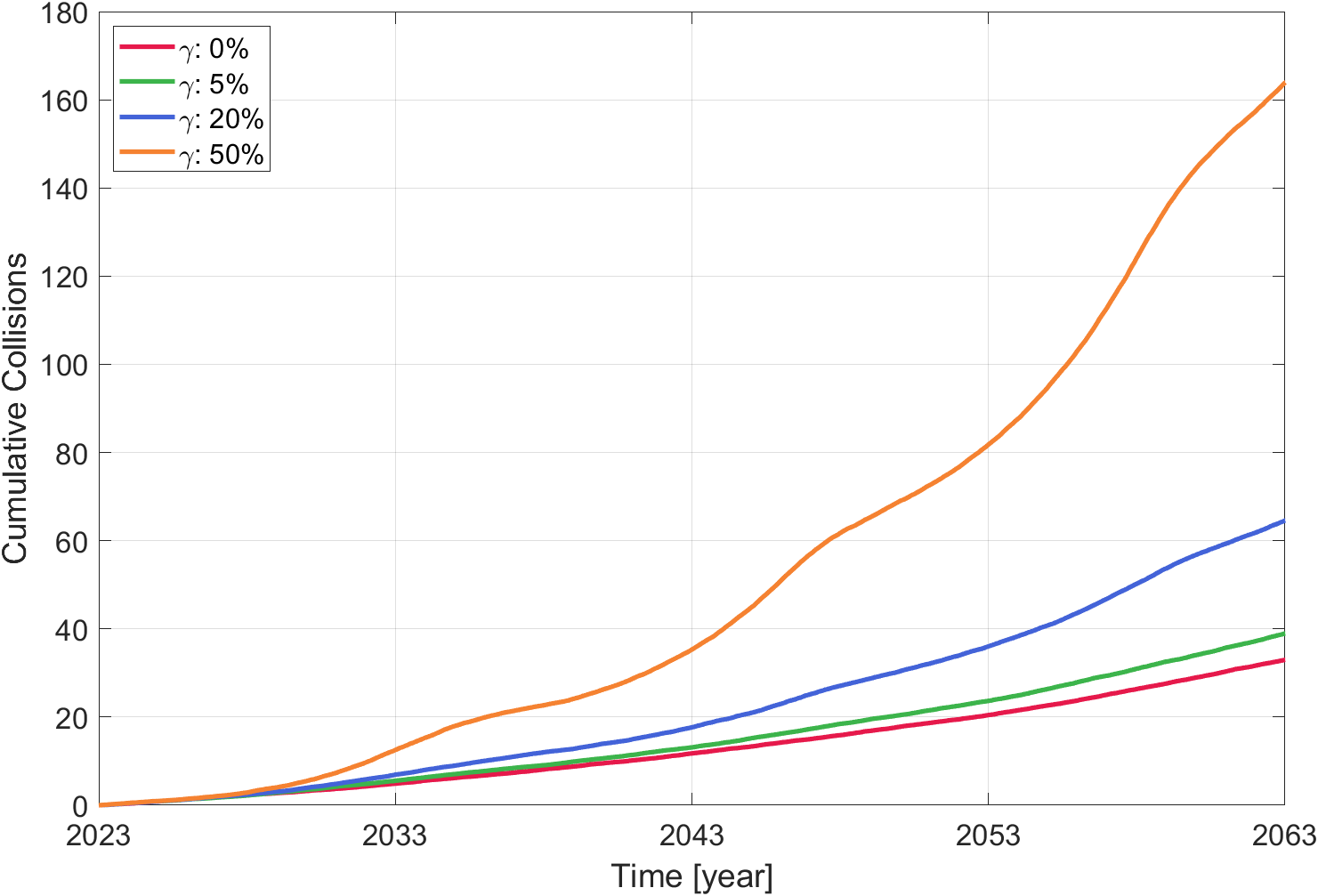}
    \caption{Cumulative number of catastrophic collisions based on different $\gamma$.} 
    \label{fig:gamma_coll}
\end{figure*}

In summary, the simulation results suggest that increasing $s_{CAM}$ or decreasing $\gamma$ can help reduce the final population of space objects. However, such mitigation strategies are insufficient to reverse the overall upward trend in fragment generation. In other words, while these measures can slow the degradation of the space environment, they are not capable of fully preventing it. This aligns with conclusions from previous studies \citep{liou2006risks}, which suggests that the mitigation measures will restrict the growth, but it is unlikely to stop the long-term debris population in LEO from increasing. Consequently, these results have led to the consideration of Active Debris Removal (ADR) as a necessary measure to remediate the space environment.

\section{Topological and Dynamical Properties of the Network}
The topological properties of a network are exclusively determined by the number of nodes and their connections, while the dynamical properties are defined by the dynamic evolution of links and nodes. Both properties can be studied by analysing an adjacency matrix derived from the evolutionary equations in Eq.\eqref{eq: nodes dynamics}, which encapsulates the connection probability among nodes. This section begins with the definition of the network links, followed by the centrality analysis of the network.

\subsection{Definition of the Network Links}\label{Definition of the links}
 
In our network model, a link between node $S_i$ and node $S_j$ is characterised by an event rate which quantifies the rate at which events affect the two nodes.  Different events happen with different rates according to Eq.\eqref{eq: nodes dynamics}, let $\chi^{\mathcal{C}}_{S_iS_j}$ represent the collision rate between node $S_i$ and node $S_j$, $\chi^{\mathcal{F}, F}_{S_iS_j}$ represent the rate of fragment flow event from node $S_i$ to node $S_j$, $\chi^{\mathcal{F}, D}_{S_iS_j}$ represent the decay rate from node $S_i$ to node $S_j$, $\chi^{\mathcal{F},SC}_{S_iS_j}$ represent the rate of transition from a payload node to a non-manoeuvrable satellite node due to small collisions, and $\chi^{\mathcal{F}, PMD}_{S_iS_j}$ represent the rate of failure PMD event leading the transition from a Payloads node a Non-manoeuvrable satellite node. Let $\delta_{\text{condition}}$ represent an indicator function which equals 1 if the condition is true, and 0 otherwise, then:
\begin{equation}
    \begin{cases}
        \chi^{\mathcal{C}}_{S_iS_j} = \tau^*_{S_iS_j} + \kappa \tau_{S_iS_j} \cdot \delta_{\left( S_i = P_i \land S_j = F_j \right) \lor \left( S_i = F_i \land S_j = P_j \right)}, \\
        \chi^{\mathcal{F},F}_{S_iS_j} = \sum_{S_l \in \mathcal{S}}\sum^n_{l=1} \tau^*_{S_iS_l} \xi^{S_j}_{S_iS_l} \cdot \delta_{\left( S_j = F_j\right)}, \\
        \chi^{\mathcal{F},SC}_{S_iS_j} = \sum^n_{l=1} \kappa\tau_{S_iS_l} \cdot \delta_{\left(i = j\right) \land \left(S_i = P_i\right) \land \left(S_j = N_j\right)},
        \\
        \chi^{\mathcal{F},D}_{S_iS_j} = \frac{\varepsilon_{S_iS_j}}{\Delta t_k},\\
        \chi_{S_i S_j}^{\mathcal{F}, PMD} = \Delta_{S_i}^\gamma \cdot \delta_{\left(i = j\right) \land \left(S_i = P_i\right) \land \left(S_j = N_j\right)},
    \end{cases}
\end{equation}
where $\mathcal{S}=\{P, U, N, F\}$ and $\tau^*_{{S_i}{S_j}}$ indicates the effective collision rate with the consideration of the manoeuvrability of payload nodes through $s_{CAM}$. For example, from the first equation in system Eq.\eqref{eq: nodes dynamics} one can see that the collision rate between a payload node and any other node must be scaled by a factor of $1 - s_{CAM}$, and by $(1 - s_{CAM})^2$ if both nodes involved are payload nodes. The quantity $\xi^{S_j}_{S_iS_l}$ indicates the fraction of fragments flowing from node $S_i$ to node $S_j$ after the collision between node $S_i$ and node $S_l$: 
\begin{equation}
   \xi^{S_j}_{S_iS_l} = \frac{\zeta^{S_j}_{S_iS_l}}{\zeta_{S_iS_l}},
\end{equation}
where $\zeta_{S_iS_l}$ indicates the total amount of fragments generated from the collision between node $S_i$ and node $S_l$. Given that $\xi^{S_j}_{S_iS_l}$ from a catastrophic collision can be estimated by a constant once the species and orbital region are fixed, $\xi^{S_j}_{S_iS_l}$ can be pre-calculated to save computational time. Assuming that there are two massive objects, with the mass of 20,000 kg for each, from node $S_i$ and node $S_j$ collide with each other within their respective orbital regions, the redistribution of new fragments is determined by NASA's breakup model and the network configuration. This redistribution process is repeated 10 times to obtain the average $\zeta^{S_j}_{S_iS_l}$ for each pair of nodes. These average values are stored in a pre-calculated tensor $\mathbf{T}$ with the corresponding element $T_{ilj}$. Quantity $\Delta^{\gamma}_{S_i}$ represents the amount of objects flowing from node $S_i$ to node $S_j$ due to the PMD failure and $\varepsilon_{S_iS_j}$ represents the amount of objects flowing from node $S_i$ to node $S_j$ due to the orbital decay, where $\varepsilon_{S_iS_j}$ is non-zero only if node $S_i$ and node $S_j$ are the same species and $j = i-1$. Note that  $\chi^{\mathcal{F}, F}_{S_iS_j}$ can be understood as the fraction of the collision events between $S_i$ and $S_j$ that produce fragments flowing to $S_j$. 
\\
We can then associate a probability $p$ that an event affects the nodes within a given time interval $\Delta t_k$.
Given a rate $\chi_{S_iS_j}$, the probability of at least one event affecting nodes $S_i$ and $S_j$ is derived by according to \cite{fleurence2007rates}:
\begin{equation}
    p_{S_iS_j} = 1-e^{-\chi_{{S_i}{S_j}}\Delta t_k}.
\end{equation}
Note that a node of type $F$ can be affected at the same time by the flow of fragments due to collisions or decay and by the collision with other nodes. Hence, in the following, the connection strength between a generic node $S_i$ and $F_j$, or $\tilde{p}_{S_iF_j}$, is defined as the maximum among the individual event probabilities
\begin{equation}
    \tilde{p}_{S_iF_j} = \max \{p^{\mathcal{C}}_{S_iF_j},p^{\mathcal{F}, F}_{S_iF_j}\},
\end{equation}
or
\begin{equation}
    \tilde{p}_{S_iF_j} = \max \{p^{\mathcal{F}, F}_{S_iF_j}, p^{\mathcal{F}, D}_{S_iF_j}\},
\end{equation}
where the former of the two is the situation in which node $S_i$ and $F_j$ belong to the same shell, while the latter is the situation in which node $S_i$ and $F_j$ belong to different shells.

\subsection{Centrality Analysis of the Network}\label{Sec:Centrality measurements}

In network analysis, in-degree and out-degree are key metrics used to describe the centrality of nodes within a directed graph. In-degree refers to the number of incoming connections that a node receives from other nodes, representing the extent to which a node is influenced by others in the network. A higher in-degree indicates that a node is central within the network, receiving more interactions. On the other hand, out-degree measures the number of outgoing connections that a node sends to other nodes. It reflects the level of influence of a node in terms of initiating interactions with others. Nodes with a high out-degree are typically more actively playing a critical role in diffusing the effects of an event across the network. The critical nodes can be identified according to their centrality. 

\begin{figure*}
\centering
    \centering
    \includegraphics[width=1\textwidth]{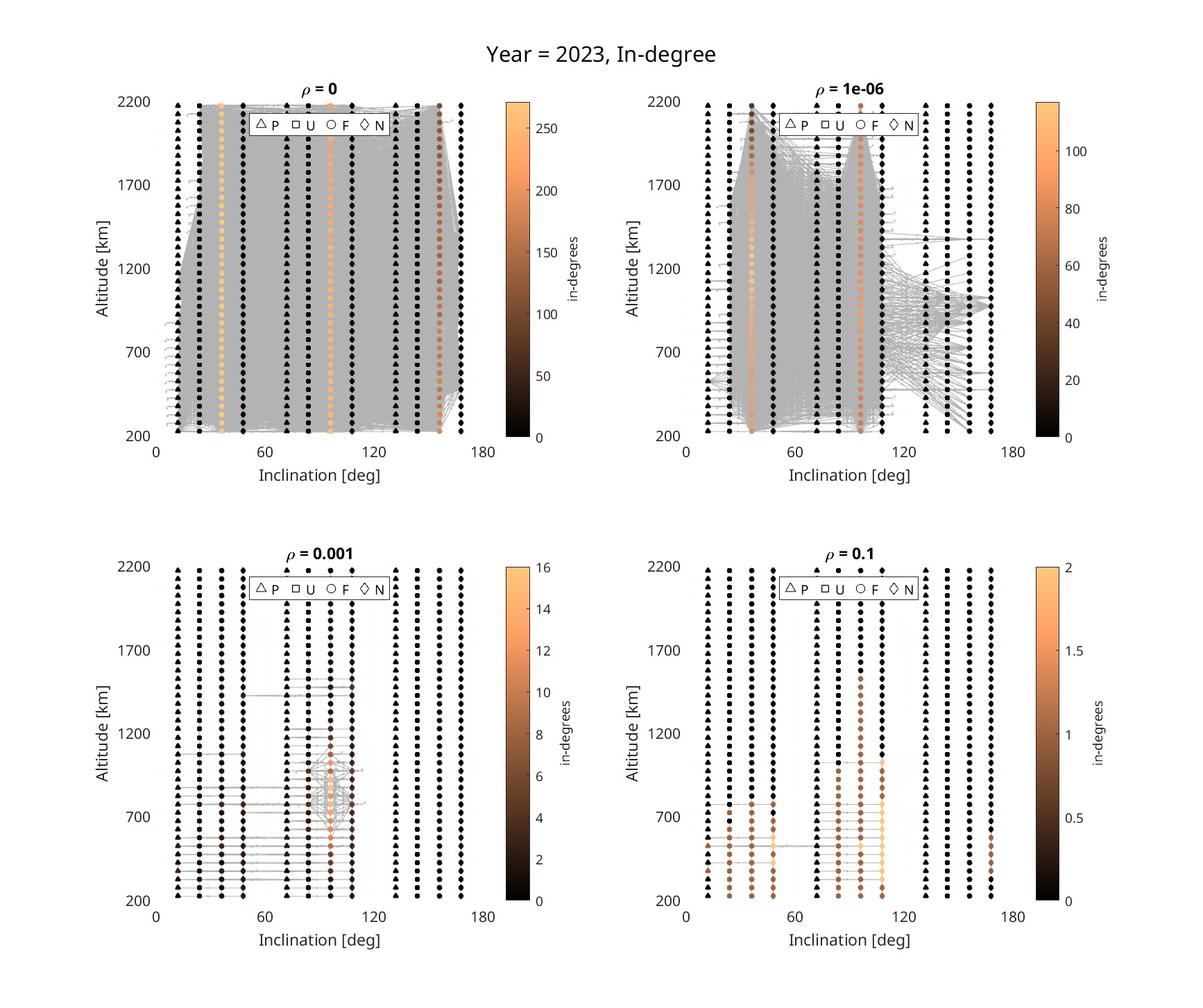}    \caption{\label{fig:network_nodes_lattice_indegree_threshold2023} Representation of the network of nodes with the lattice structure in year 2023 (In-degree). }
\end{figure*}

\begin{figure*}
\centering
    \centering
    \includegraphics[width=1\textwidth]{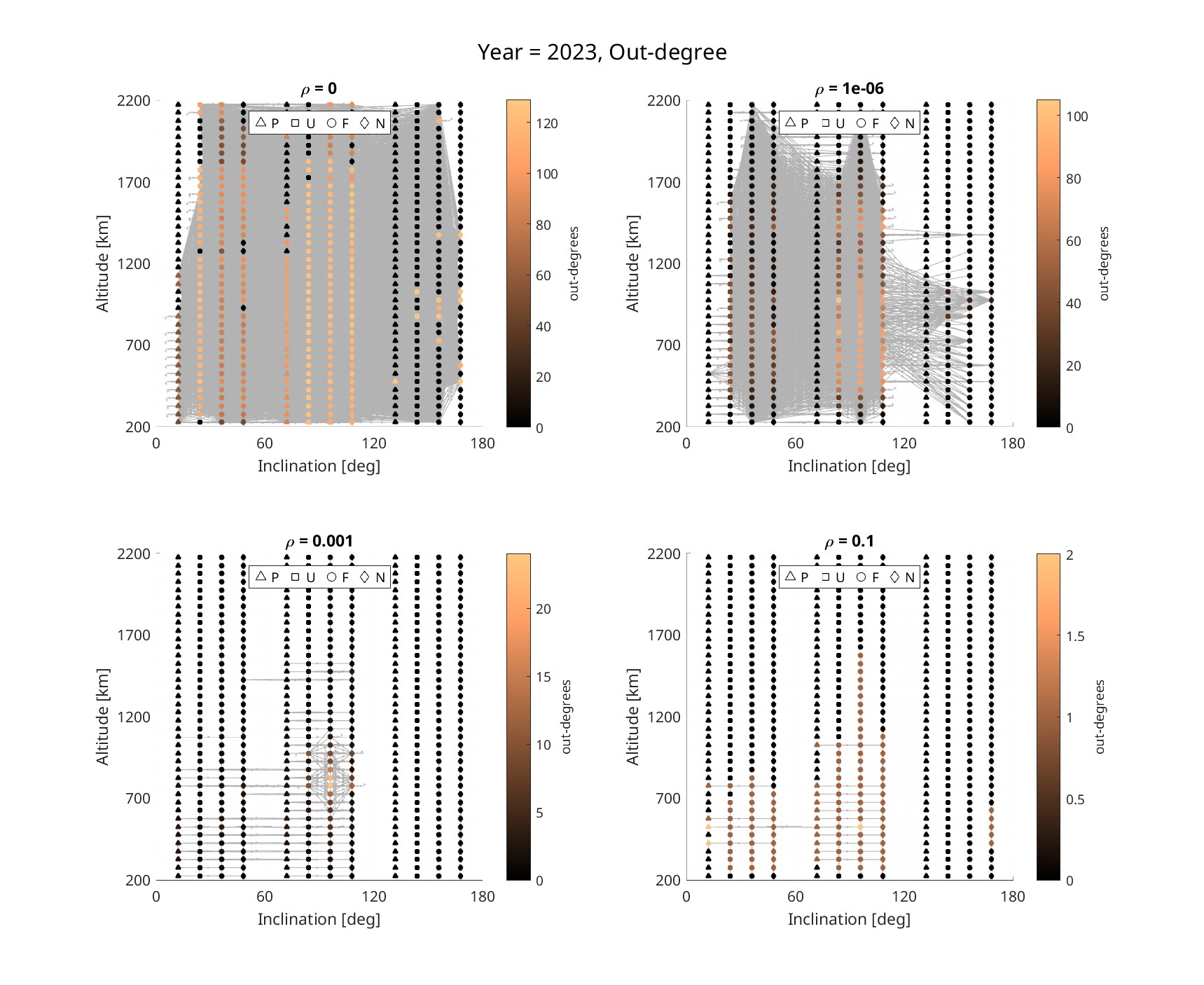}    \caption{\label{fig:network_nodes_lattice_outdegree_threshold2023} Representation of the network of nodes with the lattice structure in year 2023 (Out-degree).}
\end{figure*}

We call a subnetwork $G_{\rho}$ the set of nodes and links for which $\tilde{p}_{S_i S_j} \geq \rho$. Figure \ref{fig:network_nodes_lattice_indegree_threshold2023} shows the in-degree distribution of different subnetworks $G_{\rho}$ with $\rho$ = 0, 1e-6, 1e-3, 1e-1 at the snapshot of year 2023, where the network is represented in a lattice structure, while Figure \ref{fig:network_nodes_lattice_outdegree_threshold2023} shows the distribution of the out-degree. This visualisation provides an initial indication of the node centrality for different thresholds. As the threshold $\rho$ increases, the number of connections defining subnetwork $G_{\rho}$ decreases, resulting in changes in the criticality of the nodes. For instance, when $\rho$=0 — which includes links with non-zero $\tilde{p}$ — the $F$ nodes located within an altitude range of [1000, 1100] km and the inclination range of $[0, 60]^\circ$ show the highest in-degree, suggesting that these nodes act as strong sinks to attract fragments generated from collision events. On the other hand, the $U$ nodes in the altitude range of [1000, 1200] km and the inclination bin $[60, 120]^\circ$ show the highest out-degree, indicating that they are the most hazardous objects for spreading fragments to other nodes. However, the node with the highest in-degree/out-degree changes over $\rho$. For example, while a $U$ node has the highest out-degree in the subnetwork $G_{\rho=0}$, the associated link weights are relatively weak. As observed in the subnetwork $G_{\rho=1e-3}$, a $F$ node becomes the one with the highest out-degree. These observations highlight the need to consider both the direction and the weight of the links when conducting centrality analysis. Therefore, the weighted in-degree $\hat{d}_{in}$ and weighted out-degree $\hat{d}_{out}$ are used to evaluate the importance of each node, which considers both direction and weight of each link
\begin{equation}\label{eq:in_out_degree}
\begin{cases}
    \hat{d}_{in}(S_i) = \sum_{S_j \in \mathcal{S}}\sum^n_{j=1} \tilde{p}_{S_iS_j},\\
    \hat{d}_{out}(S_i) = \sum_{S_j \in \mathcal{S}}\sum^n_{j=1} \tilde{p}_{S_iS_j}.
\end{cases}
\end{equation}

Table \ref{Table: Top5 nodes of centrality measures_Baseline} presents a comprehensive summary of the top 5 nodes with the highest centrality values for the Baseline case of year 2023, 2063 and 2103. 
Specifically, in 2023, node 92, which represents a group of non-manoeuvrable satellites located within the altitude range of $[550, 600]$ km and the inclination bin $[60, 120]$ degrees —emerges as the most susceptible. This vulnerability is primarily attributed to payload failures in PMD and the decay of non-manoeuvrable satellites from higher orbital shells. A group of fragments located in the altitude range of $[500, 550]$ km and inclination bin $[60, 120]$ degrees appears to be the most hazardous, acting as a key source for propagating fragments throughout the network. Interestingly, the same top nodes maintain their ranks through 2063 and 2103 with only minor variation in centrality values. This consistency suggests the presence of strong self-connections, driven by self-collision processes in which the majority of resulting fragments are redistributed back into the originating node, thereby making these nodes increasingly critical over time.

\begin{table}[H]
\centering
\caption{\label{Table: Top5 nodes of centrality measures_Baseline} Top 5 nodes of centrality measures based on the snapshots of 2023, 2063 and 2103, Baseline case.}
\begin{tabular}{lcccccccccc}
\hline
\multirow{2}{*}{\textbf{Snapshot}} & \multirow{2}{*}{\textbf{Centrality}}                     & \multicolumn{3}{c}{\textbf{Top nodes}}                                               & \multirow{2}{*}{\begin{tabular}[c]{@{}c@{}}\textbf{Objects} \\ \textbf{number}\end{tabular}} & \multicolumn{4}{c}{\textbf{Location}}                                                                                                                                                & \multirow{2}{*}{\textbf{Class}} \\ \cline{3-5} \cline{7-10}
                          &                                                 & rank & node ID & \begin{tabular}[c]{@{}c@{}}Centrality\\ value\end{tabular} &                                                                            & \multicolumn{2}{c}{\begin{tabular}[c]{@{}c@{}}height bin\\ {[}km{]}\end{tabular}} & \multicolumn{2}{c}{\begin{tabular}[c]{@{}c@{}}inclination bin\\ {[}deg{]}\end{tabular}} &                        \\ \hline
\multirow{10}{*}{2023}    & \multirow{5}{*}{in-degree}                      & 1    & 92      & 1.98692695950322                                           & 124                                                                        & 550                                     & 600                                     & 60                                         & 120                                        & 'N'                    \\
                          &                                                 & 2    & 68      & 1.93500519104489                                           & 92                                                                         & 450                                     & 500                                     & 60                                         & 120                                        & 'N'                    \\
                          &                                                 & 3    & 80      & 1.80164839711700                                           & 73                                                                         & 500                                     & 550                                     & 60                                         & 120                                        & 'N'                    \\
                          &                                                 & 4    & 56      & 1.63303430119017                                           & 71                                                                         & 400                                     & 450                                     & 60                                         & 120                                        & 'N'                    \\
                          &                                                 & 5    & 79      & 1.53668845441548                                           & 186                                                                        & 500                                     & 550                                     & 60                                         & 120                                        & 'F'                    \\ \cline{2-11} 
                          & \multirow{5}{*}{out-degree}                     & 1    & 79      & 1.52485346012650                                           & 186                                                                        & 500                                     & 550                                     & 60                                         & 120                                        & 'F'                    \\
                          &                                                 & 2    & 139     & 1.18178525137863                                           & 1054                                                                       & 750                                     & 800                                     & 60                                         & 120                                        & 'F'                    \\
                          &                                                 & 3    & 151     & 1.14877931440353                                           & 1343                                                                       & 800                                     & 850                                     & 60                                         & 120                                        & 'F'                    \\
                          &                                                 & 4    & 73      & 1.10366201012284                                           & 2952                                                                       & 500                                     & 550                                     & 0                                          & 60                                         & 'P'                    \\
                          &                                                 & 5    & 91      & 1.09927242312778                                           & 305                                                                        & 550                                     & 600                                     & 60                                         & 120                                        & 'F'                    \\ \hline
\multirow{10}{*}{2063}    & \multirow{5}{*}{in-degree}                      & 1    & 151     & 1.09397358650202                                           & 1746                                                                       & 800                                     & 850                                     & 60                                         & 120                                        & 'F'                    \\
                          &                                                 & 2    & 139     & 1.09234148243687                                           & 1515                                                                       & 750                                     & 800                                     & 60                                         & 120                                        & 'F'                    \\
                          &                                                 & 3    & 187     & 1.06219846137200                                           & 840                                                                        & 950                                     & 1000                                    & 60                                         & 120                                        & 'F'                    \\
                          &                                                 & 4    & 163     & 1.06008867762470                                           & 1255                                                                       & 850                                     & 900                                     & 60                                         & 120                                        & 'F'                    \\
                          &                                                 & 5    & 127     & 1.05488332007186                                           & 1151                                                                       & 700                                     & 750                                     & 60                                         & 120                                        & 'F'                    \\ \cline{2-11} 
                          & \multicolumn{1}{l}{\multirow{5}{*}{out-degree}} & 1    & 139     & 1.12787968798332                                           & 1515                                                                       & 750                                     & 800                                     & 60                                         & 120                                        & 'F'                    \\
                          & \multicolumn{1}{l}{}                            & 2    & 151     & 1.12240859815599                                           & 1746                                                                       & 800                                     & 850                                     & 60                                         & 120                                        & 'F'                    \\
                          & \multicolumn{1}{l}{}                            & 3    & 187     & 1.06855298298059                                           & 840                                                                        & 950                                     & 1000                                    & 60                                         & 120                                        & 'F'                    \\
                          & \multicolumn{1}{l}{}                            & 4    & 127     & 1.05439108561404                                           & 1151                                                                       & 700                                     & 750                                     & 60                                         & 120                                        & 'F'                    \\
                          & \multicolumn{1}{l}{}                            & 5    & 163     & 1.04825800053616                                           & 1255                                                                       & 850                                     & 900                                     & 60                                         & 120                                        & 'F'                    \\ \hline
\multirow{10}{*}{2103}    & \multicolumn{1}{l}{\multirow{5}{*}{in-degree}}  & 1    & 151     & 1.08232922063680                                           & 1863                                                                       & 800                                     & 850                                     & 60                                         & 120                                        & 'F'                    \\
                          & \multicolumn{1}{l}{}                            & 2    & 139     & 1.07444292331339                                           & 1673                                                                       & 750                                     & 800                                     & 60                                         & 120                                        & 'F'                    \\
                          & \multicolumn{1}{l}{}                            & 3    & 187     & 1.07005626493219                                           & 1100                                                                       & 950                                     & 1000                                    & 60                                         & 120                                        & 'F'                    \\
                          & \multicolumn{1}{l}{}                            & 4    & 127     & 1.06066326660445                                           & 1306                                                                       & 700                                     & 750                                     & 60                                         & 120                                        & 'F'                    \\
                          & \multicolumn{1}{l}{}                            & 5    & 163     & 1.05581902179307                                           & 1455                                                                       & 850                                     & 900                                     & 60                                         & 120                                        & 'F'                    \\ \cline{2-11} 
                          & \multicolumn{1}{l}{\multirow{5}{*}{out-degree}} & 1    & 151     & 1.10387559309562                                           & 1863                                                                       & 800                                     & 850                                     & 60                                         & 120                                        & 'F'                    \\
                          & \multicolumn{1}{l}{}                            & 2    & 139     & 1.09485846928317                                           & 1673                                                                       & 750                                     & 800                                     & 60                                         & 120                                        & 'F'                    \\
                          & \multicolumn{1}{l}{}                            & 3    & 187     & 1.08359031135662                                           & 1100                                                                       & 950                                     & 1000                                    & 60                                         & 120                                        & 'F'                    \\
                          & \multicolumn{1}{l}{}                            & 4    & 127     & 1.06927172375704                                           & 1306                                                                       & 700                                     & 750                                     & 60                                         & 120                                        & 'F'                    \\
                          & \multicolumn{1}{l}{}                            & 5    & 163     & 1.04502793146937                                           & 1455                                                                       & 850                                     & 900                                     & 60                                         & 120                                        & 'F'                    \\ \hline
\end{tabular}
\end{table}

\newpage
\section{Space Carrying Capacity and Network Dynamics}
In this section, we study the relationship between network dynamics and the concept of space carrying capacity. The idea of looking for equilibrium points was already introduced in previous works, see \cite{talent1992} or \cite{CarryingCapacity2024}. Here we apply the same idea to the network dynamics.
\subsection{Carrying capacity of one-dimensional network dynamics}
In ecological networks, the carrying capacity of an environment in relation to the population $x$ living in that environment and the available resources is expressed by the Verhulst model \citep{verhulst1838}:
\begin{equation}\label{eq:verhulst}
\dot{x}=r x -\frac{r}{K}x^2,
\end{equation}
where $r$ is the birth rate minus the death rate of the population and $K$ is the carrying capacity of the environment.
Eq. (\ref{eq:verhulst}) has two equilibrium points: $x=0$ and $x=K$. Hence, the carrying capacity is an equilibrium of the population dynamics. In the case of the Verhulst model, it is easy to see that $x=K$ is stable, hence if the population is higher than $K$, it reduces to $K$ and if it is lower than $K$, it grows till reaching $K$.

In analogy with model  (\ref{eq:verhulst}) we can consider the case in which, asymptotically, in the case of no new launches, the space environment is occupied only by a single species, the fragments. From Eq.(\ref{eq: change in the nodes}) one can derive the discrete one dimensional average equation:
\begin{equation}\label{eq:1D_evo}
\Delta x_k=-\frac{\mathbb{E}(\Delta x_{decay})}{x_k}x_k+\frac{\mathbb{E}(\Delta x^{FF}_{collision})}{x_k^2}x^2_k,
\end{equation}
where $x_k$ is the average number of objects at step $k$, $\mathbb{E}(\Delta x_{decay})$ indicates the expected change in fragments due to the atmospheric drag
\begin{equation}
\label{eq: expected change in fragments due to the atmospheric drag}
    \mathbb{E}(\Delta x_{decay}) = \sum^n_{i=1} \left(\varepsilon^+_{F_i} - \varepsilon^-_{F_i}\right),
\end{equation}
and $\mathbb{E}(\Delta x^{FF}_{collision})$ indicates the expected change in fragments due to the self-collisions
\begin{equation}
\label{eq: expected change in fragments due to the self-collisions}
    \mathbb{E}(\Delta x^{FF}_{collision}) = \sum_i \left( -\sum_j \tau_{F_iF_j} \Delta t_k + \sum_j \sum_l \tau_{F_jF_l}\zeta^{F_i}_{F_jF_l}  \Delta t_k\right).
\end{equation}
To maintain the analogy with Eq. (\ref{eq:verhulst}) we will study the continuous version of Eq. (\ref{eq:1D_evo}). The mean evolution of the fragments can be written as: 
\begin{equation}\label{eq:mean_fragments}
\dot{x}=-a x +bx^2,
\end{equation}
where $a=\frac{\mathbb{E}(\Delta x_{decay})}{x\Delta t }$ is the rate of decay per unit time and object and $b=\frac{\mathbb{E}(\Delta x^{FF}_{collision})}{x^2 \Delta t}$ is the production of debris per unit time and object square. This form of the equation is equivalent to the result of \cite{talent1992}. If the cross-sectional area of the objects remains constant and the decay is not dependent on the solar cycle, then $a/b$ is a constant. 
In this case, Eq. (\ref{eq:mean_fragments}) has two equilibrium points, $x=0$ and $x=a/b$, where $x=0$ is stable and $x=a/b$ is unstable. More specifically, for $x>a/b$ the environment diverges to infinity and for $x<a/b$ it converges to 0.
In fact, Eq. (\ref{eq:mean_fragments}) is a Bernulli equation with solution:
\begin{equation}\label{eq:bernulli_solution}
x(t)=\frac{e^{-at}}{b\frac{e^{-at}}{a}+C_0},
\end{equation}
which has a vertical asymptote for $b\frac{e^{-at}}{a}+C_0=0$ if $x(0)>b/a$. We can, therefore, define $a/b$ as the carrying capacity since above this value the environment diverges to infinity. 

Given that different populations $x$ might result in different coefficients, due to the distribution of cross-section areas, and coefficients are not constant over time, here we take an average across multiple initial populations and an average over the propagation time.
Table \ref{table:2D_init_pop} presents the selected populations for the test cases, which were propagated using both NESSY and the analytical solution of Eq. (\ref{eq:mean_fragments}).

\begin{table}[h!]
\begin{center}
\caption{Initial populations test cases for carrying capacity analysis.}
\begin{tabular}{c c c} 
 \hline
 \textbf{Test case} & \textbf{$\#$ fragments} & \textbf{$\#$ payloads} \\ [0.5ex] 
 \hline
 Baseline & 9804 & 5471 \\ 
 
 Test 1 & 39020 & 13942 \\
 
 Test 2 & 68628 & 13942 \\
 
 Test 3 & 98040 & 19413 \\
 \hline

\end{tabular}
\label{table:2D_init_pop}
\end{center}
\end{table}

\begin{figure*}[htbp]
\centering
    \centering    \includegraphics[width=0.9\textwidth]{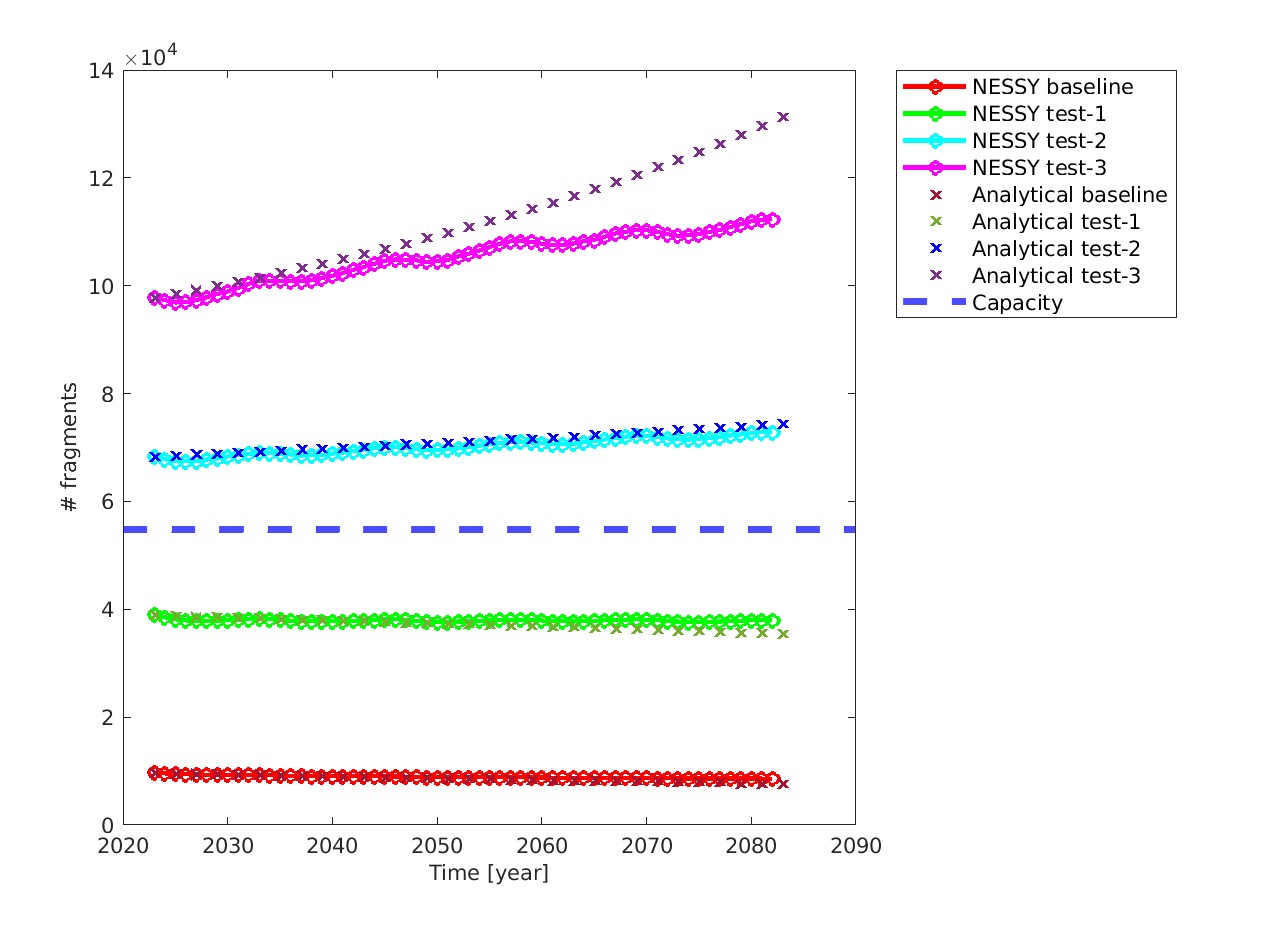}
    \caption{Carrying capacity analysis just with fragments population (test cases defined in Table \ref{table:2D_init_pop}).} \label{fig:1d_capacity}
\end{figure*}

Figure \ref{fig:1d_capacity} shows the result of the environment propagation with NESSY (dotted lines) for the different initial conditions, considering just one node for the network. Note that the capacity here is the mean value over several solar cycles. It is clear that when the initial population is higher than the capacity (horizontal dashed line), the population diverges, while if it is initially lower, it tends to zero, as expected. The propagation with NESSY was performed with a 1-year time step for 60 years and considering 30 independent runs to compute the mean trend. Note that, although all the simulations in this section were run with a longer time step, compared to previous sections, to speed up the analysis of multiple initial conditions, the results match rather well the analytical solution. The capacity is computed with respect to the baseline case, which corresponds to the actual population in 2023, averaging the coefficients across the random runs and time. The NESSY propagation is compared with the propagation performed with the analytical equations shown before (represented with crosses in the plot).

One can now ask which other equilibrium points the network dynamics has. In this paper, we limit our investigations to a simple mean-field case with two nodes representing the fragments and all new objects launched into space. More complex analyses over the whole network are deferred to future studies. Note that this case, although simple, is representative of a range of realistic situations, as it will be explained in the remainder of this section.

\subsection{Carrying capacity of two-dimensional network dynamics}

If one considers only fragments $x$ and payloads with new launches $y$, the evolutionary equations simplify to:

\begin{equation}\label{eq:2D_diff_case}
\begin{array}{l}
\dot{x}=b x^2-ax+c y^2+dxy,\\
\dot{y}=-ey^2-fxy+\lambda-\gamma y.
\end{array}
\end{equation}
Coefficient $a$ is the rate of fragments decay per unit time and object, $b$ is the rate of generation of fragments per unit time and fragment squared, $c$ is the rate of generation of fragments per unit time and payload squared, $d$ is the rate of fragment generation per unit time and object squared due to collisions between fragments and payloads, $e$ is the number of collisions between payloads per unit time and payload squared, $f$ is the number of collisions between payloads and fragments per unit time and object squared. The other two coefficients are $\lambda$, which is considered for the following simulations, 3000 payloads launched per year, and $\gamma$, which is instead computed as the rate of post-mission disposals per unit object. No small collision, $\kappa = 0$, is considered in this case. The equilibrium points are found by solving:
\begin{equation}\label{eq:2D_equi_case}
\begin{array}{l}
b x^2-a x+c y^2+d xy=0,\\
-ey^2-fxy+\lambda-\gamma y=0,
\end{array}
\end{equation}
where the value of the coefficients is:
\begin{equation}\label{eq:coeff_2d}
\begin{cases}
    a = \frac{\mathbb{E}(\Delta x_{decay})}{x \Delta t_k },\\
    
    b = \frac{\mathbb{E}(\Delta x^{FF}_{collision})}{x^2  \Delta t_k},\\
    
    c = \frac{\mathbb{E}(\Delta x^{PP}_{collision})}{y^{2}  \Delta t_k },\\
    
    d = \frac{\mathbb{E}(\Delta x^{PF}_{collision})}{xy \Delta t_k},\\
    
    e = \frac{\mathbb{E}(\Delta y^{PP}_{collision})}{y^{2}  \Delta t_k },\\
    
    f = \frac{\mathbb{E}(\Delta y^{PF}_{collision})}{x y   \Delta t_k},\\
    
    \gamma = \frac{\mathbb{E}(y_{PMD})}{y \Delta t_k}.
\end{cases}
\end{equation}
The expected values $\mathbb{E}(x_{decay})$ and $\mathbb{E}(\Delta x^{FF}_{collision})$ have the same meaning as in Eq.\eqref{eq: expected change in fragments due to the atmospheric drag} and Eq.\eqref{eq: expected change in fragments due to the self-collisions} under the assumption that there is no decay of payloads due to atmospheric drag. Quantity $\mathbb{E}(\Delta x^{PP}_{collision})$ indicates the expected change in fragments due to the self-collisions of payloads. 
\begin{equation}
    \mathbb{E}(\Delta x^{PP}_{collision}) = \sum_i \left( - \sum_{j} (1-s_{CAM})^2\tau_{P_iP_j}\Delta t_k + \sum_j \sum_l \tau_{P_jP_l}\zeta^{F_i}_{P_jP_l}  \Delta t_k \right),
\end{equation}
$\mathbb{E}(x^{PF}_{collision})$ indicates the expected change in fragments due to the collisions between payloads and fragments. 
\begin{equation}
    \mathbb{E}(\Delta x^{PF}_{collision}) = \sum_i \left( - \sum_{j} (1-s_{CAM})\tau_{P_iF_j}\Delta t_k + \sum_j \sum_l \tau_{P_jF_l}\zeta^{F_i}_{P_jF_l}  \Delta t_k \right),
\end{equation}
$\mathbb{E}(\Delta y^{PP}_{collision})$ indicates the expected change in payloads due to the self-collisions of payloads
\begin{equation}
    \mathbb{E}(\Delta y^{PP}_{collision}) = \sum_i \left( - \sum_{j} (1-s_{CAM})^2\tau_{P_iP_j}\Delta t_k\right ),
\end{equation}
$\mathbb{E}(\Delta y^{PF}_{collision})$ indicates the expected change in payloads due to the collisions between payloads and fragments
\begin{equation}
    \mathbb{E}(\Delta y^{PF}_{collision}) = \sum_i \left( - \sum_{j} (1-s_{CAM})\tau_{P_iF_j}\Delta t_k  \right).
\end{equation}
Quantity $\mathbb{E}(y_{PMD})$ indicates the expected change in payloads due to PMD. Similar to the one-dimensional case, we take an average across multiple initial populations and an average over the propagation time for each coefficient.

\begin{figure*}[htbp]
\centering
    \centering
    \includegraphics[width=1\textwidth]{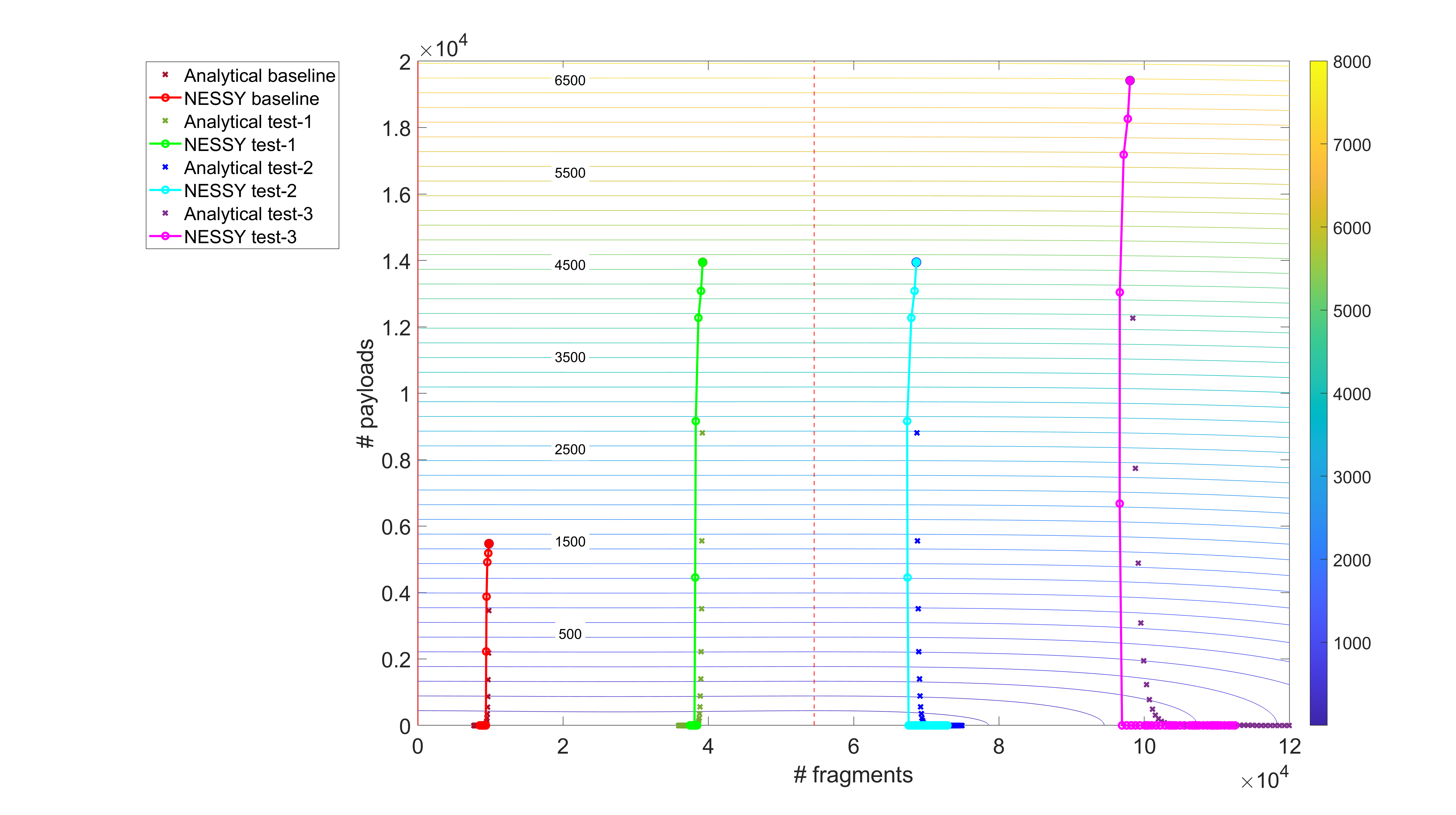}
    \caption{Carrying capacity analysis for the initial populations in Table \ref{table:2D_init_pop} without new launches. Vertical red lines identify the equilibrium points.} \label{fig:2d_frag}
\end{figure*}

\begin{figure*}[htbp]
\centering
    \centering
    \includegraphics[width=0.8\textwidth]{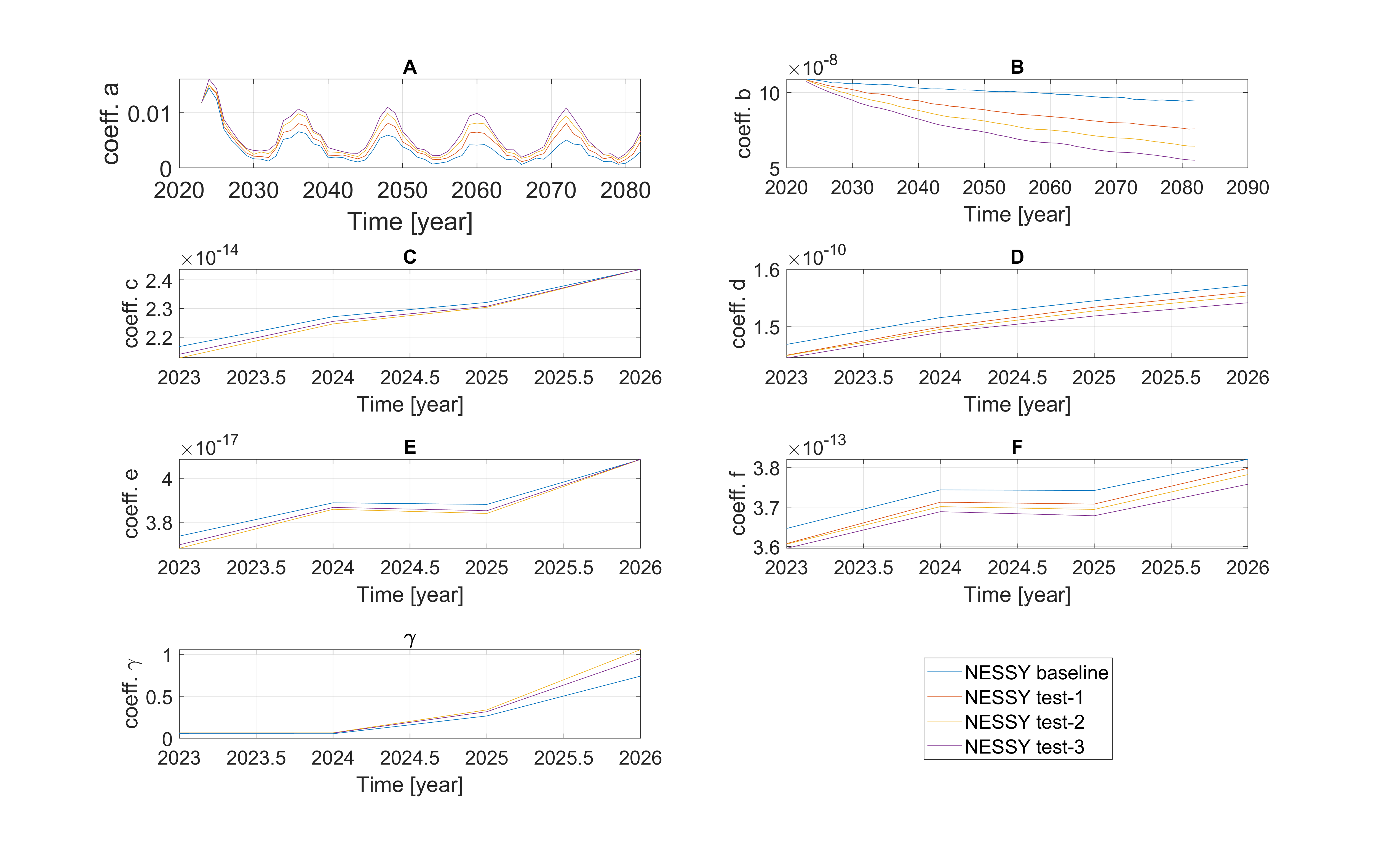}
    \caption{Coefficients evolution in time for the initial populations in Table \ref{table:2D_init_pop}, without new launches.} \label{fig:coeff_frag}
\end{figure*}

\begin{figure*}[htbp]
\centering
    \centering
    \includegraphics[width=1\textwidth]{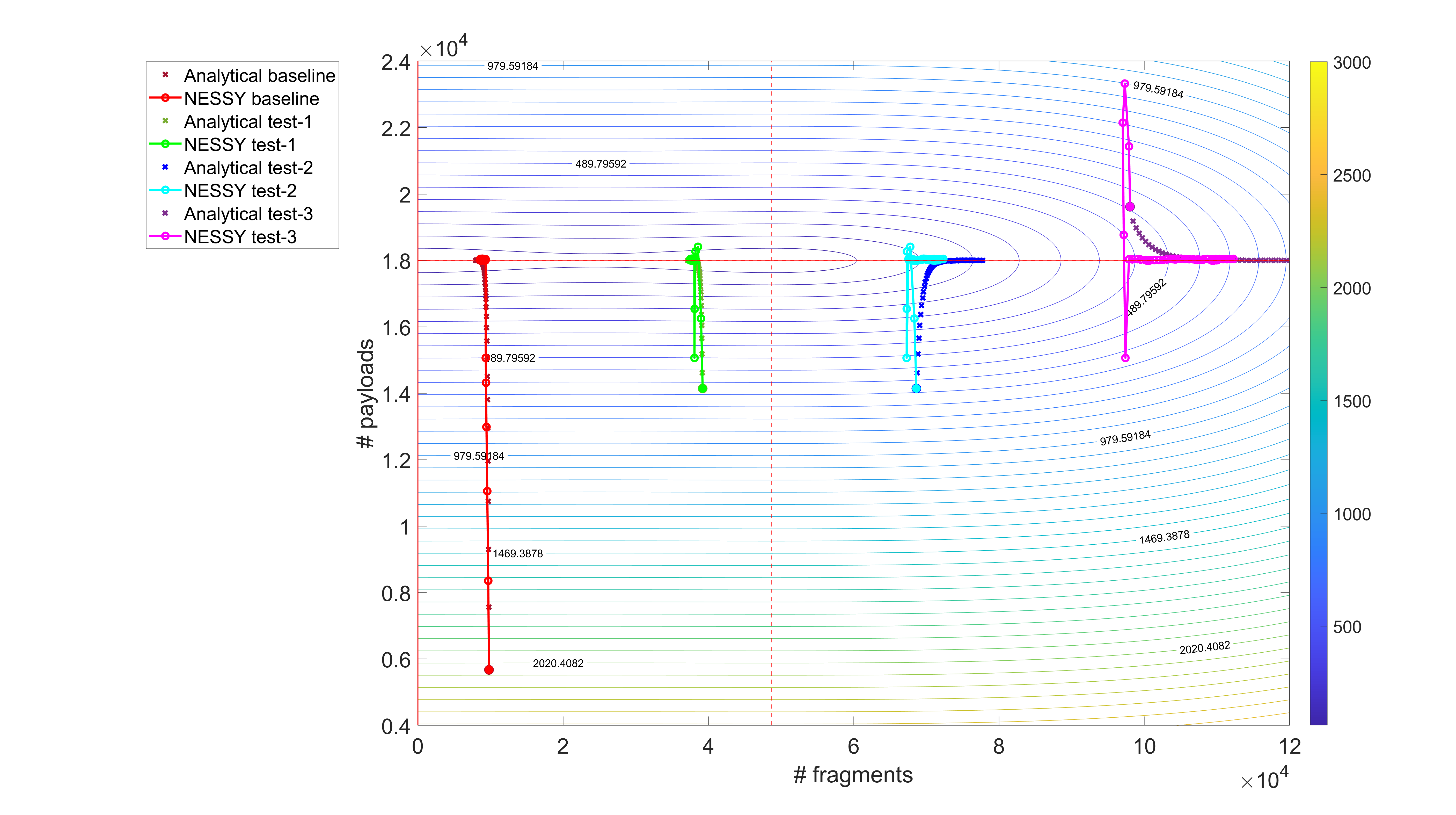}
    \caption{Carrying capacity analysis for the initial populations in Table \ref{table:2D_init_pop}, with new launches and $s_{CAM}$ = 99.99 \%. Vertical red lines identify the equilibrium points.} \label{fig:2d_01}
\end{figure*}

\begin{figure*}[htbp]
\centering
    \centering
    \includegraphics[width=0.8\textwidth]{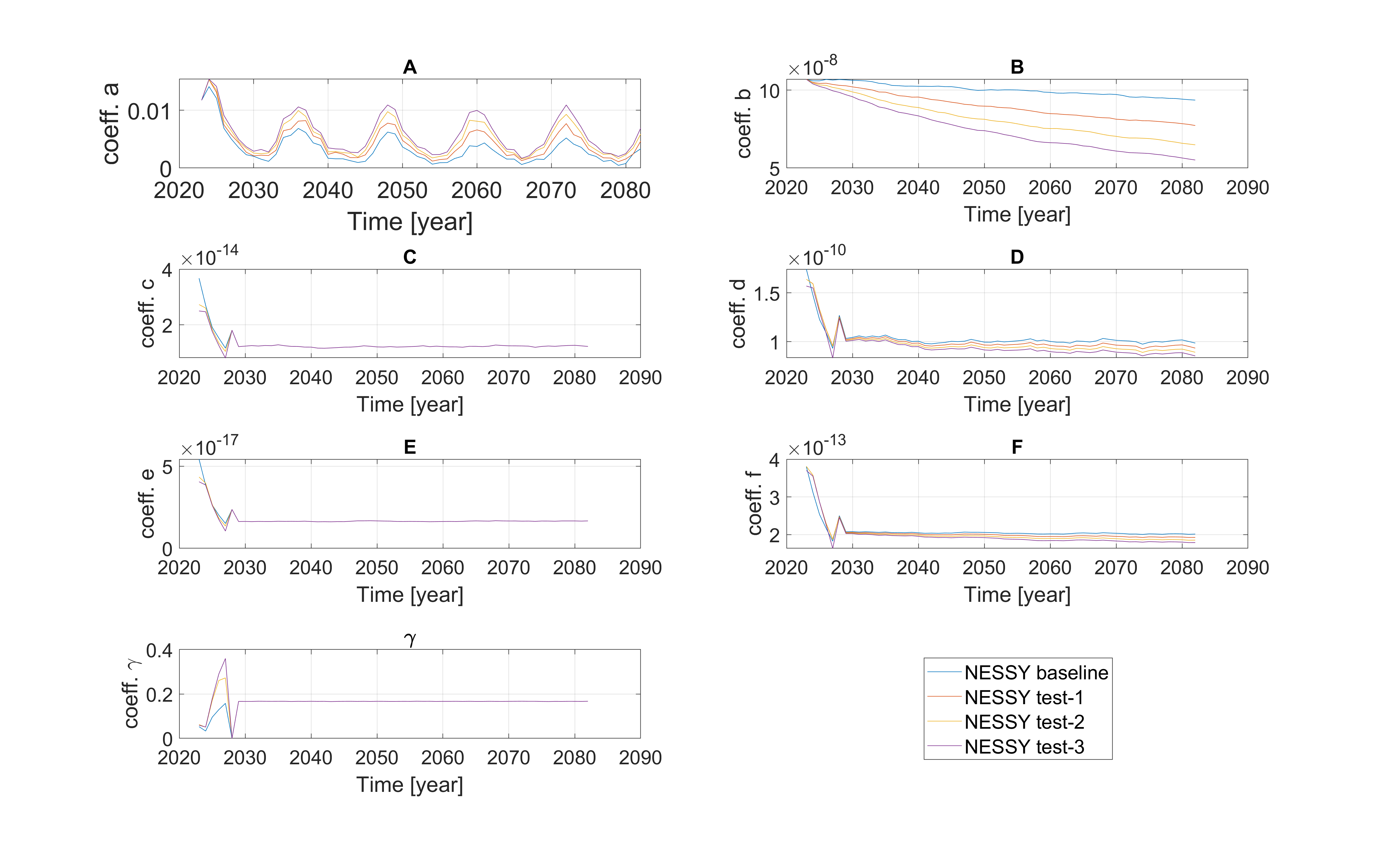}
    \caption{Coefficients evolution in time for the initial populations in Table \ref{table:2D_init_pop} with new launches, for $s_{CAM}$ = 99.99 \%.} \label{fig:coeff_01}
\end{figure*}

\begin{figure*}[htbp]
\centering
    \centering
    \includegraphics[width=1\textwidth]{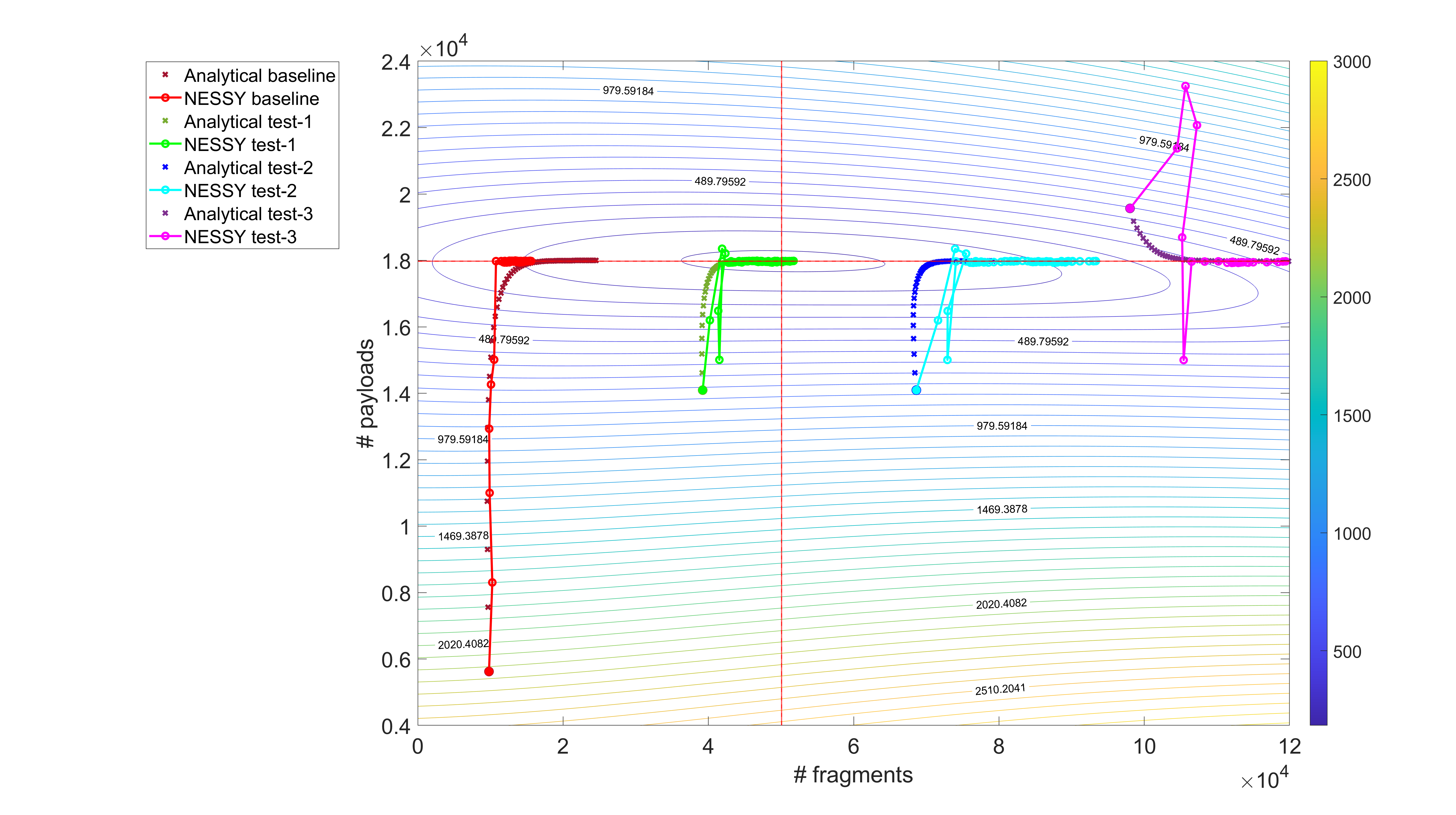}
    \caption{Carrying capacity analysis for the initial populations in Table \ref{table:2D_init_pop}, with new launches and $s_{CAM}$ = 0 \%.} \label{fig:2d_100}
\end{figure*}

\begin{figure*}[htbp]
\centering
    \centering
    \includegraphics[width=0.8\textwidth]{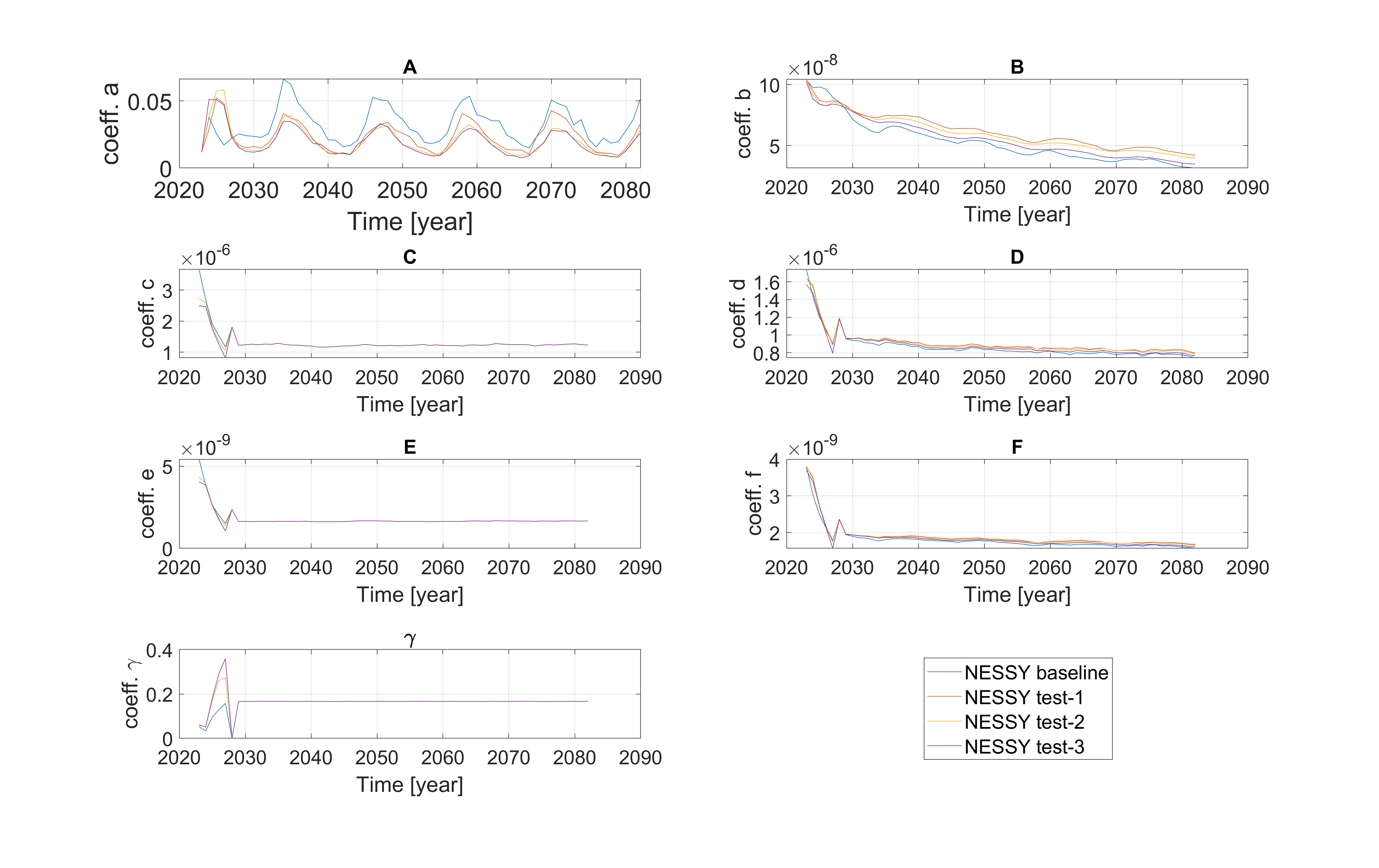}
    \caption{Coefficients evolution in time for the initial populations in Table \ref{table:2D_init_pop} with new launches, for $s_{CAM}$ = 0 \%.} \label{fig:coeff_100}
\end{figure*}

NESSY was run again from the initial populations in Table \ref{table:2D_init_pop}: one test does not consider new launches, while the others include
new launches. Table \ref{table:params_2D} reports the settings used to run NESSY on these cases.  
\begin{table}[h!]
\centering
\caption{NESSY  mission-related parameters.}
\begin{tabular}{c c c c} 
 \hline
 $s_{CAM}$ & $\gamma$ & Payload Mission lifetime & $\kappa$ \\ [0.5ex] 
 \hline
 99.99$\%$ & 0 $\%$ & 5 years & 0 \\ [1ex] 
 \hline
\end{tabular}
\label{table:params_2D}
\end{table}
Each initial population was propagated for 60 years with a one-year time step, and each run was repeated 30 times to compute the average variation of each species. The network is composed of only one node for each species. 
The results of the propagations are presented in Figure\ref{fig:2d_01} to Figure\ref{fig:coeff_100}. Figures \ref{fig:coeff_frag}, \ref{fig:coeff_01}, and \ref{fig:coeff_100} show the time evolution of the coefficients of Eqs. (\ref{eq:2D_diff_case}) derived from NESSY. Figure \ref{fig:2d_frag}, \ref{fig:2d_01}, and \ref{fig:2d_100} show the evolution of the two species computed with NESSY and compared to the propagation of Eq (\ref{eq:2D_diff_case}) (analytical tests and baseline). The level curves represent the value of $\sqrt{\dot{x}^{2} + \dot{y}^{2}}$. The equilibrium points were initially identified by finding the minima of $\sqrt{\dot{x}^{2} + \dot{y}^{2}}$ with the MATLAB \textit{fminunc} function.

For a zero launch rate, one equilibrium solution is $[0,0]$ while the other is $[a/b,0]$. From the Jacobian matrix, one can see that $[0,0]$ is linearly stable because the Jacobian has two real negative eigenvalues, while for $[a/b,0]$, it has one negative and one positive real eigenvalue.
Hence, there is a stable manifold that brings the space environment towards $[a/b,0]$ and an unstable manifold that diverges from $[a/b,0]$ and either connects to the stable manifold of $[0,0]$ or diverges to infinity.
Figure \ref{fig:2d_frag} shows the existence of these two equilibrium solutions. The propagation of system (\ref{eq:2D_diff_case}) (indicated with an "x" symbol in Figure \ref{fig:2d_frag}), with the average of the coefficients in Figure \ref{fig:coeff_frag}), shows the expected convergence to the stable point and the stable manifold along the y axis. The propagation with NESSY (continuous lines with circle markers) follows the same trend, confirming the theory. 

For a non-zero launch rate, we can distinguish two extreme cases: a) all new launches are manoeuverable and can avoid collisions with 99.99\% probability, or b) all new launches are non-manoeuvrable and can collide with each other and with fragments. 
In the former case, system (\ref{eq:2D_diff_case}) still has two equilibrium solutions: $[\epsilon,y_{eq}]$ and $[k,y_{eq}]$. Where again $[\epsilon,y_{eq}]$ is stable and  $[k,y_{eq}]$ is unstable with one stable manifold and one unstable manifold. The value $\epsilon$ is small for high collision avoidance rate and goes to zero when collision avoidance is 100\%.
Figure \ref{fig:2d_01} shows the case with manoeuvrable payloads. The corresponding time history of the coefficients can be found in Figure \ref{fig:coeff_01}. Coefficient $a$ follows the variation of atmospheric density as expected, while $b$ progressively reduces due to the decay of fragments and the reduction in cross-sectional area of the fragments. It is important to note that in this example, there are no collisions between lethal non-trackable fragments and other objects or fragments.

Continuous lines with circle markers in Figure \ref{fig:2d_01} 
 represent the propagation with NESSY while the propagation of the system (\ref{eq:2D_diff_case}) is represented by the "x" marker.
In this case, one can see two distinct equilibrium points for a non-zero value of $y$, one stable and one unstable. Note that when the simulations with NESSY start above $y_{eq}$, $y$ initially grows and then converges to the equilibrium value. This initial growth, observed only in the early years, is attributed to resident payloads. Most of these payloads are not removed by PMD because they remain within their operational lifetimes. 

Figure \ref{fig:2d_100} shows the case with non-manoeuvrable payloads. The time history of the coefficients is in Figure \ref{fig:coeff_100}. In this case, there is no equilibrium solution for the fragments, while the payloads converge to a stable equilibrium. 

Hence, we can say that the capacity of the environment to host new objects is given by the region of the $[x,y]$ space in between the two stable manifolds respectively of $[\epsilon,y_{eq}]$ and $[k,y_{eq}]$. When the two points collapse, there is zero capacity.
In general, if the removal rate, given by a combination of PMD and decay,  plus the collision avoidance rate are higher than the generation of new fragments, the system presents two distinct equilibria. 

The values of the coefficients and respective equilibria are reported in Table \ref{tab:equilibria_coeff}. Note that for the case with 0\% $s_{CAM}$ we report only the equilibrium value for $\dot{y}=0$ since there is no equilibrium $\dot{x}=0$. 

\begin{table}[]       
    \centering
    \caption{Equilibrium solutions and coefficients of system (\ref{eq:2D_diff_case})}    
    \begin{tabular}{ccc}
    \hline
         \textbf{Case} & \textbf{Coefficients} & \textbf{Equilibria}  \\
         \hline
         Baseline & 
         \begin{tabular}{l}
         a = 0.004728332083372\\
         b = 8.662467642990248e-08\\
         c = 1.175401267752297e-14\\
         d = 9.428437648428035e-10\\
         $\gamma$ = 0.368578793358788\\
         $\lambda$ = 0\\
         e = 2.003922397999517e-17\\
         f = 2.316928055993169e-13\\
         \end{tabular}
         & \begin{tabular}{l}
         x=0; y=0\\
x=5.4584e+04; y=0
\end{tabular}\\
         \hline
         99.99\% $s_{CAM}$ & 
         \begin{tabular}{l} 
         a = 0.004226706317436\\
         b = 8.676619156862889e-08\\
         c = 1.224456162356393e-15\\
         d = 9.606253682748494e-11\\
         $\gamma$ = 0.166677662732838\\
         $\lambda$ = 3000\\
         e = 1.647891773627737e-17\\
        f = 1.957039536003774e-13 
        \end{tabular} &  
        \begin{tabular}{l} 
        x=0.0;   y=1.8e4\\
        x=4.90e4;    y=1.8e4
        \end{tabular}\\
         \hline         
         0\% $s_{CAM}$     & 
         \begin{tabular}{l} 
         a= 0.022592560002365\\
         b = 6.214276689402071e-08\\
         c = 1.224460289412396e-07\\
         d = 8.664860267617623e-07\\
         $\gamma$ = 0.166621505691604\\
         $\lambda$ = 3000\\
         e = 1.647986448760241e-09\\
         f = 1.790973469946202e-09  
         \end{tabular} & 
         y=1.8e4  \\
         \hline         
    \end{tabular} 
\end{table}\label{tab:equilibria_coeff} 

\newpage
\section{Conclusions}
This paper proposed a novel approach to investigate the global behaviour of the space environment and study the dynamic evolution of the relationships among space objects. Different from previous works, we first developed a network model of the space environment and then derived the evolutionary equations from the structure of the network. 
The validation of the evolutionary equations showed that the proposed network model can provide results that are comparable to other state-of-the-art long-term evolutionary models.

Furthermore, it was shown how the network model offers the unique ability to identify critical regions of the space environment. For instance, the node with the highest impact on the rest of the environment can be identified by analysing the in-degree and out-degree of each node. The nodes responsible for spreading or accumulating fragments can be identified through the out-degree or in-degree centrality measures. From the perspective of space operators, instead of focusing on the entire space environment, limited resources such as monitoring and mitigation efforts can be strategically directed toward these identified priority space objects, optimising the allocation of resources based on their importance in comparison to others. 

It was also shown how to derive the equilibrium solutions of the network dynamics and how these equilibrium solutions can be stable or unstable depending on the removal rate or the manoeuvrability of the resident objects. It was found that the space carrying capacity can be directly quantified as the ratio between the decay rate and the collision rate if only fragments are present in the environment. In this case, an initial population above the capacity would diverge to infinity, while an initial population below the capacity would converge to zero.
If, on top of the fragments, new objects are launched into space, it was found that there are generally two equilibria, one stable and one unstable. One can then identify a region within which the environment would converge to a stable point. Outside this stability region, instead, the environment would diverge to infinity.
In this paper, the analysis was limited to only two species and did not include lethal non-trackable objects that would not lead to a catastrophic collision but would disable an active payload.

Future work will consider introducing other types of nodes, representing constellations, objects in highly elliptical orbits, and other mitigation and remediation actions. Furthermore, explosions and small non-trackable lethal fragments need to be included in the evolutionary model together with other orbit regimes.
New orbit regimes, like Geostationary Orbits or the cis-lunar environment, can be simply added as further sub-networks. This is the subject of current work. 

\section*{Acknowledgments}
This work was supported through the Horizon 2020 MSCA ETN Stardust-R (grant agreement number: 813644). ESA OSIP Multi-layer Network Model of the Space Environment, contract 4000129448/20/NL/MH, Project officer: Stijn Lemmens.
The authors would like to thank Callum Wilson for providing the launch traffic model used in this paper.

\bibliographystyle{jasr-model5-names}
\biboptions{authoryear}
\bibliography{nessy_main}

\begin{thebibliography}{47}
\expandafter\ifx\csname natexlab\endcsname\relax\def\natexlab#1{#1}\fi
\ifx\xfnm\relax \def\xfnm[#1]{\unskip,\space#1}\fi

\bibitem[{Acciarini \& Vasile(2020)}]{acciarini2020multi}
\bibinfo{author}{Acciarini, G.},  \& \bibinfo{author}{Vasile, M.}
  (\bibinfo{year}{2020}).
\newblock \bibinfo{title}{A multi-layer temporal network model of the space
  environment}.
\newblock In {\it \bibinfo{booktitle}{71st International Astronautical
  Congress}\/}.

\bibitem[{Acciarini \& Vasile(2021)}]{acciarini2021network}
\bibinfo{author}{Acciarini, G.},  \& \bibinfo{author}{Vasile, M.}
  (\bibinfo{year}{2021}).
\newblock \bibinfo{title}{A network-based evolutionary model of the space
  environment}.
\newblock In {\it \bibinfo{booktitle}{8th European Conference on Space
  Debris}\/} (pp. \bibinfo{pages}{1--9}).

\bibitem[{Adilov et~al.(2018)Adilov, Alexander \& Cunningham}]{2018An}
\bibinfo{author}{Adilov, N.}, \bibinfo{author}{Alexander, P.~J.},  \&
  \bibinfo{author}{Cunningham, B.~M.} (\bibinfo{year}{2018}).
\newblock \bibinfo{title}{An economic "kessler syndrome": A dynamic model of
  earth orbit debris}.
\newblock {\it \bibinfo{journal}{Economics Letters}\/},  {\it
  \bibinfo{volume}{166}\/}\bibinfo{issue}{(MAY)}, \bibinfo{pages}{79--82}.

\bibitem[{Boley \& Byers(2021)}]{2021Satellite}
\bibinfo{author}{Boley, A.},  \& \bibinfo{author}{Byers, M.}
  (\bibinfo{year}{2021}).
\newblock \bibinfo{title}{Satellite mega-constellations create risks in low
  earth orbit, the atmosphere and on earth}.
\newblock {\it \bibinfo{journal}{Nature Publishing Group}\/}, .

\bibitem[{Bowman et~al.(2008)Bowman, Tobiska, Marcos, Huang, Lin \&
  Burke}]{bowman2008new}
\bibinfo{author}{Bowman, B.}, \bibinfo{author}{Tobiska, W.~K.},
  \bibinfo{author}{Marcos, F.} et~al. (\bibinfo{year}{2008}).
\newblock \bibinfo{title}{A new empirical thermospheric density model jb2008
  using new solar and geomagnetic indices}.
\newblock In {\it \bibinfo{booktitle}{AIAA/AAS astrodynamics specialist
  conference and exhibit}\/} (p. \bibinfo{pages}{6438}).

\bibitem[{Bradley \& Wein(2009)}]{2009Space}
\bibinfo{author}{Bradley, A.~M.},  \& \bibinfo{author}{Wein, L.~M.}
  (\bibinfo{year}{2009}).
\newblock \bibinfo{title}{Space debris: Assessing risk and responsibility}.
\newblock {\it \bibinfo{journal}{Advances in Space Research}\/},  {\it
  \bibinfo{volume}{43}\/}\bibinfo{issue}{(9)}, \bibinfo{pages}{1372--1390}.

\bibitem[{D'Ambrosio et~al.(2022)D'Ambrosio, Lifson \& Linares}]{d2022capacity}
\bibinfo{author}{D'Ambrosio, A.}, \bibinfo{author}{Lifson, M.},  \&
  \bibinfo{author}{Linares, R.} (\bibinfo{year}{2022}).
\newblock \bibinfo{title}{The capacity of low earth orbit computed using
  source-sink modeling}.
\newblock {\it \bibinfo{journal}{arXiv preprint arXiv:2206.05345}\/}, .

\bibitem[{De~Marchi et~al.(2025)De~Marchi, Wang \& Vasile}]{PietroSPS}
\bibinfo{author}{De~Marchi, P.}, \bibinfo{author}{Wang, Y.},  \&
  \bibinfo{author}{Vasile, M.} (\bibinfo{year}{2025}).
\newblock \bibinfo{title}{Space environment impact of solar power satellites}.
\newblock In {\it \bibinfo{booktitle}{2025 AAS/AIAA Space Flight Mechanics
  Meeting}\/}.

\bibitem[{Dolado-Perez et~al.(2013)Dolado-Perez, Di~Constanzo \&
  Revelin}]{jc2013introducing}
\bibinfo{author}{Dolado-Perez, J.}, \bibinfo{author}{Di~Constanzo, R.},  \&
  \bibinfo{author}{Revelin, B.} (\bibinfo{year}{2013}).
\newblock \bibinfo{title}{Introducing medee--a new orbital debris evolutionary
  model}.
\newblock In {\it \bibinfo{booktitle}{Proceeding of the 6th European Conference
  on Space Debris}\/}.

\bibitem[{D’Ambrosio \& Linares(2024)}]{CarryingCapacity2024}
\bibinfo{author}{D’Ambrosio, A.},  \& \bibinfo{author}{Linares, R.}
  (\bibinfo{year}{2024}).
\newblock \bibinfo{title}{Carrying capacity of low earth orbit computed using
  source-sink models}.
\newblock {\it \bibinfo{journal}{Journal of Spacecraft and Rockets}\/},  {\it
  \bibinfo{volume}{0}\/}\bibinfo{issue}{(0)}, \bibinfo{pages}{1--17}.
  \URLprefix \url{https://doi.org/10.2514/1.A35729}.
  \DOIprefix\doi{10.2514/1.A35729}.
  \href{http://arxiv.org/abs/https://doi.org/10.2514/1.A35729}{\tt
  arXiv:https://doi.org/10.2514/1.A35729}.

\bibitem[{D’Ambrosio et~al.(2023)D’Ambrosio, Servadio, Siew \&
  Linares}]{dambrosio2023}
\bibinfo{author}{D’Ambrosio, A.}, \bibinfo{author}{Servadio, S.},
  \bibinfo{author}{Siew, P.~M.} et~al. (\bibinfo{year}{2023}).
\newblock \bibinfo{title}{Novel source–sink model for space environment
  evolution with orbit capacity assessment}.
\newblock {\it \bibinfo{journal}{Journal of SpaceCraft and Rockets}\/}, .

\bibitem[{Fleurence \& Hollenbeak(2007)}]{fleurence2007rates}
\bibinfo{author}{Fleurence, R.~L.},  \& \bibinfo{author}{Hollenbeak, C.~S.}
  (\bibinfo{year}{2007}).
\newblock \bibinfo{title}{Rates and probabilities in economic modelling:
  transformation, translation and appropriate application}.
\newblock {\it \bibinfo{journal}{Pharmacoeconomics}\/},  {\it
  \bibinfo{volume}{25}\/}, \bibinfo{pages}{3--6}.

\bibitem[{Henning et~al.(2019)Henning, Sorge, Peterson, Jenkin, McVey \&
  Mains}]{henning2019impacts}
\bibinfo{author}{Henning, G.}, \bibinfo{author}{Sorge, M.},
  \bibinfo{author}{Peterson, G.} et~al. (\bibinfo{year}{2019}).
\newblock \bibinfo{title}{Impacts of large constellations and mission disposal
  guidelines on the future space debris environment}.
\newblock {\it \bibinfo{journal}{IAC-19, A6, 2, 7, x50024}\/}, .

\bibitem[{Jang et~al.(2023)Jang, Gusmini, Siew, D’Ambrosio, Servadio, Machuca
  \& Linares}]{jang2023monte}
\bibinfo{author}{Jang, D.}, \bibinfo{author}{Gusmini, D.},
  \bibinfo{author}{Siew, P.~M.} et~al. (\bibinfo{year}{2023}).
\newblock \bibinfo{title}{Monte carlo methods to model the evolution of the low
  earth orbit population}.
\newblock In {\it \bibinfo{booktitle}{33rd AAS/AIAA Space Flight Mechanics
  Meeting, Austin, TX}\/}.
\newblock volume~\bibinfo{volume}{1}.

\bibitem[{Jang et~al.(2025)Jang, Gusmini, Siew, D’Ambrosio, Servadio, Machuca
  \& Linares}]{jang2025new}
\bibinfo{author}{Jang, D.}, \bibinfo{author}{Gusmini, D.},
  \bibinfo{author}{Siew, P.~M.} et~al. (\bibinfo{year}{2025}).
\newblock \bibinfo{title}{New monte carlo model for the space environment}.
\newblock {\it \bibinfo{journal}{Journal of Spacecraft and Rockets}\/},  (pp.
  \bibinfo{pages}{1--22}).

\bibitem[{Johnson et~al.(2001)Johnson, Krisko, Liou \&
  Anz-Meador}]{johnson2001nasa}
\bibinfo{author}{Johnson, N.~L.}, \bibinfo{author}{Krisko, P.~H.},
  \bibinfo{author}{Liou, J.-C.} et~al. (\bibinfo{year}{2001}).
\newblock \bibinfo{title}{Nasa's new breakup model of evolve 4.0}.
\newblock {\it \bibinfo{journal}{Advances in Space Research}\/},  {\it
  \bibinfo{volume}{28}\/}\bibinfo{issue}{(9)}, \bibinfo{pages}{1377--1384}.

\bibitem[{Kessler(1981)}]{KESSLER1981}
\bibinfo{author}{Kessler, D.~J.} (\bibinfo{year}{1981}).
\newblock \bibinfo{title}{Derivation of the collision probability between
  orbiting objects: the lifetimes of jupiter's outer moons}.
\newblock {\it \bibinfo{journal}{Icarus}\/},  {\it
  \bibinfo{volume}{48}\/}\bibinfo{issue}{(1)}, \bibinfo{pages}{39--48}.
  \URLprefix
  \url{https://www.sciencedirect.com/science/article/pii/0019103581901512}.
  \DOIprefix\doi{https://doi.org/10.1016/0019-1035(81)90151-2}.

\bibitem[{Kessler \& Cour-Palais(1978)}]{1978Collision}
\bibinfo{author}{Kessler, D.~J.},  \& \bibinfo{author}{Cour-Palais, B.~G.}
  (\bibinfo{year}{1978}).
\newblock \bibinfo{title}{Collision frequency of artificial satellites: The
  creation of a debris belt}.
\newblock {\it \bibinfo{journal}{Journal of Geophysical Research Space
  Physics}\/},  {\it \bibinfo{volume}{83}\/}\bibinfo{issue}{(A6)},
  \bibinfo{pages}{2637--2646}.

\bibitem[{Krisko(2011)}]{krisko2011proper}
\bibinfo{author}{Krisko, P.~H.} (\bibinfo{year}{2011}).
\newblock \bibinfo{title}{Proper implementation of the 1998 nasa breakup
  model}.
\newblock {\it \bibinfo{journal}{Orbital Debris Quarterly News}\/},  {\it
  \bibinfo{volume}{15}\/}\bibinfo{issue}{(4)}, \bibinfo{pages}{1--10}.

\bibitem[{Letizia et~al.(2017)Letizia, Colombo, Lewis \&
  Krag}]{letizia2017extending}
\bibinfo{author}{Letizia, F.}, \bibinfo{author}{Colombo, C.},
  \bibinfo{author}{Lewis, H.} et~al. (\bibinfo{year}{2017}).
\newblock \bibinfo{title}{Extending the ecob space debris index with
  fragmentation risk estimation}, .

\bibitem[{Letizia et~al.(2016)Letizia, Colombo, Lewis \&
  Krag}]{LETIZIA20161255}
\bibinfo{author}{Letizia, F.}, \bibinfo{author}{Colombo, C.},
  \bibinfo{author}{Lewis, H.~G.} et~al. (\bibinfo{year}{2016}).
\newblock \bibinfo{title}{Assessment of breakup severity on operational
  satellites}.
\newblock {\it \bibinfo{journal}{Advances in Space Research}\/},  {\it
  \bibinfo{volume}{58}\/}\bibinfo{issue}{(7)}, \bibinfo{pages}{1255--1274}.
  \URLprefix
  \url{https://www.sciencedirect.com/science/article/pii/S0273117716302393}.
  \DOIprefix\doi{https://doi.org/10.1016/j.asr.2016.05.036}.

\bibitem[{Lewis et~al.(2010)Lewis, Newland, Swinerd \& Saunders}]{lewis2010new}
\bibinfo{author}{Lewis, H.~G.}, \bibinfo{author}{Newland, R.~J.},
  \bibinfo{author}{Swinerd, G.~G.} et~al. (\bibinfo{year}{2010}).
\newblock \bibinfo{title}{A new analysis of debris mitigation and removal using
  networks}.
\newblock {\it \bibinfo{journal}{Acta Astronautica}\/},  {\it
  \bibinfo{volume}{66}\/}\bibinfo{issue}{(1-2)}, \bibinfo{pages}{257--268}.

\bibitem[{Liou et~al.(2004)Liou, Hall, Krisko \& Opiela}]{liou2004legend}
\bibinfo{author}{Liou, J.-C.}, \bibinfo{author}{Hall, D.},
  \bibinfo{author}{Krisko, P.} et~al. (\bibinfo{year}{2004}).
\newblock \bibinfo{title}{Legend--a three-dimensional leo-to-geo debris
  evolutionary model}.
\newblock {\it \bibinfo{journal}{Advances in Space Research}\/},  {\it
  \bibinfo{volume}{34}\/}\bibinfo{issue}{(5)}, \bibinfo{pages}{981--986}.

\bibitem[{Liou \& Johnson(2006)}]{liou2006risks}
\bibinfo{author}{Liou, J.-C.},  \& \bibinfo{author}{Johnson, N.~L.}
  (\bibinfo{year}{2006}).
\newblock \bibinfo{title}{Risks in space from orbiting debris}.

\bibitem[{Liou \& Johnson(2009)}]{2009A}
\bibinfo{author}{Liou, J.~C.},  \& \bibinfo{author}{Johnson, N.~L.}
  (\bibinfo{year}{2009}).
\newblock \bibinfo{title}{A sensitivity study of the effectiveness of active
  debris removal in leo}.
\newblock {\it \bibinfo{journal}{Acta Astronautica}\/},  {\it
  \bibinfo{volume}{64}\/}\bibinfo{issue}{(2-3)}, \bibinfo{pages}{236--243}.

\bibitem[{Liou et~al.(2003)Liou, Kessler, Matney \& Stansbery}]{liou2003new}
\bibinfo{author}{Liou, J.-C.}, \bibinfo{author}{Kessler, D.~J.},
  \bibinfo{author}{Matney, M.} et~al. (\bibinfo{year}{2003}).
\newblock \bibinfo{title}{A new approach to evaluate collision probabilities
  among asteroids, comets, and kuiper belt objects}.
\newblock In {\it \bibinfo{booktitle}{Lunar and Planetary Science
  Conference}\/} (p. \bibinfo{pages}{1828}).

\bibitem[{Lucken \& Giolito(2019)}]{lucken2019collision}
\bibinfo{author}{Lucken, R.},  \& \bibinfo{author}{Giolito, D.}
  (\bibinfo{year}{2019}).
\newblock \bibinfo{title}{Collision risk prediction for constellation design}.
\newblock {\it \bibinfo{journal}{Acta Astronautica}\/},  {\it
  \bibinfo{volume}{161}\/}, \bibinfo{pages}{492--501}.

\bibitem[{Martinusi et~al.(2015)Martinusi, Dell’Elce \&
  Kerschen}]{martinusi2015analytic}
\bibinfo{author}{Martinusi, V.}, \bibinfo{author}{Dell’Elce, L.},  \&
  \bibinfo{author}{Kerschen, G.} (\bibinfo{year}{2015}).
\newblock \bibinfo{title}{Analytic propagation of near-circular satellite
  orbits in the atmosphere of an oblate planet}.
\newblock {\it \bibinfo{journal}{Celestial Mechanics and Dynamical
  Astronomy}\/},  {\it \bibinfo{volume}{123}\/}, \bibinfo{pages}{85--103}.

\bibitem[{McKnight \& Lorenzen(1989)}]{mcknight1989collision}
\bibinfo{author}{McKnight, D.},  \& \bibinfo{author}{Lorenzen, G.}
  (\bibinfo{year}{1989}).
\newblock \bibinfo{title}{Collision matrix for low earth orbit satellites}.
\newblock {\it \bibinfo{journal}{Journal of Spacecraft and Rockets}\/},  {\it
  \bibinfo{volume}{26}\/}\bibinfo{issue}{(2)}, \bibinfo{pages}{90--94}.

\bibitem[{Muelhaupt et~al.(2019)Muelhaupt, Sorge, Morin \& Wilson}]{2019Space}
\bibinfo{author}{Muelhaupt, T.~J.}, \bibinfo{author}{Sorge, M.~E.},
  \bibinfo{author}{Morin, J.} et~al. (\bibinfo{year}{2019}).
\newblock \bibinfo{title}{Space traffic management in the new space era}.
\newblock {\it \bibinfo{journal}{Journal of Space Safety Engineering}\/},  {\it
  \bibinfo{volume}{6}\/}\bibinfo{issue}{(2)}, \bibinfo{pages}{80--87}.

\bibitem[{Murakami et~al.(2019)Murakami, Nag, Lifson \&
  Kopardekar}]{murakami2019space}
\bibinfo{author}{Murakami, D.~D.}, \bibinfo{author}{Nag, S.},
  \bibinfo{author}{Lifson, M.} et~al. (\bibinfo{year}{2019}).
\newblock \bibinfo{title}{Space traffic management with a nasa uas traffic
  management (utm) inspired architecture}.
\newblock In {\it \bibinfo{booktitle}{AIAA Scitech 2019 Forum}\/} (p.
  \bibinfo{pages}{2004}).

\bibitem[{Office(2024)}]{esa_space_debris_office_discosweb_2024}
\bibinfo{author}{Office, E. S.~D.} (\bibinfo{year}{2024}).
\newblock \bibinfo{title}{{DISCOSweb}}.
\newblock \URLprefix \url{https://discosweb.esoc.esa.int/}.

\bibitem[{Pardini \& Anselmo(2014)}]{pardini2014review}
\bibinfo{author}{Pardini, C.},  \& \bibinfo{author}{Anselmo, L.}
  (\bibinfo{year}{2014}).
\newblock \bibinfo{title}{Review of past on-orbit collisions among cataloged
  objects and examination of the catastrophic fragmentation concept}.
\newblock {\it \bibinfo{journal}{Acta Astronautica}\/},  {\it
  \bibinfo{volume}{100}\/}, \bibinfo{pages}{30--39}.

\bibitem[{Rodriguez-Fernandez et~al.(2024)Rodriguez-Fernandez, Sarangerel,
  Siew, Machuca, Jang \& Linares}]{rodriguez2024towards}
\bibinfo{author}{Rodriguez-Fernandez, V.}, \bibinfo{author}{Sarangerel, S.},
  \bibinfo{author}{Siew, P.~M.} et~al. (\bibinfo{year}{2024}).
\newblock \bibinfo{title}{Towards a machine learning-based approach to predict
  space object density distributions}.
\newblock In {\it \bibinfo{booktitle}{AIAA SCITECH 2024 Forum}\/} (p.
  \bibinfo{pages}{1673}).

\bibitem[{Romano et~al.(2024)Romano, Carletti \& Daquin}]{Romano_2024}
\bibinfo{author}{Romano, M.}, \bibinfo{author}{Carletti, T.},  \&
  \bibinfo{author}{Daquin, J.} (\bibinfo{year}{2024}).
\newblock \bibinfo{title}{The resident space objects network: A complex system
  approach for shaping space sustainability}.
\newblock {\it \bibinfo{journal}{The Journal of the Astronautical Sciences}\/},
   {\it \bibinfo{volume}{71}\/}\bibinfo{issue}{(4)}. \URLprefix
  \url{http://dx.doi.org/10.1007/s40295-024-00449-4}.
  \DOIprefix\doi{10.1007/s40295-024-00449-4}.

\bibitem[{Rossi et~al.()Rossi, Anselmo, Pardini, Valsecchi \&
  Jehn}]{anselmo_pardini_SDM3}
\bibinfo{author}{Rossi, A.}, \bibinfo{author}{Anselmo, L.},
  \bibinfo{author}{Pardini, C.} et~al. ().
\newblock \bibinfo{title}{Analysis of the space debris environment with sdm
  3.0}.
\newblock In {\it \bibinfo{booktitle}{55th International Astronautical Congress
  of the International Astronautical Federation, the International Academy of
  Astronautics, and the International Institute of Space Law}\/}.
\newblock \URLprefix
  \url{https://arc.aiaa.org/doi/abs/10.2514/6.IAC-04-IAA.5.12.1.06}.
  \DOIprefix\doi{10.2514/6.IAC-04-IAA.5.12.1.06}.
  \href{http://arxiv.org/abs/https://arc.aiaa.org/doi/pdf/10.2514/6.IAC-04-IAA.5.12.1.06}{\tt
  arXiv:https://arc.aiaa.org/doi/pdf/10.2514/6.IAC-04-IAA.5.12.1.06}.

\bibitem[{Rossi et~al.(2015)Rossi, Valsecchi \& Alessi}]{rossi2015criticality}
\bibinfo{author}{Rossi, A.}, \bibinfo{author}{Valsecchi, G.},  \&
  \bibinfo{author}{Alessi, E.} (\bibinfo{year}{2015}).
\newblock \bibinfo{title}{The criticality of spacecraft index}.
\newblock {\it \bibinfo{journal}{Advances in Space Research}\/},  {\it
  \bibinfo{volume}{56}\/}\bibinfo{issue}{(3)}, \bibinfo{pages}{449--460}.

\bibitem[{Stein(2022)}]{stein_l_leo_2022}
\bibinfo{author}{Stein, L.} (\bibinfo{year}{2022}).
\newblock \bibinfo{title}{{LEO} {SATCOM} report}, .
\newblock \URLprefix
  \url{https://www.primemoverslab.com/resources/ideas/leo-satcom.pdf}.

\bibitem[{Sturza \& Carretero(2021)}]{2021Design}
\bibinfo{author}{Sturza, M.},  \& \bibinfo{author}{Carretero, G.~S.}
  (\bibinfo{year}{2021}).
\newblock \bibinfo{title}{Design trades for environmentally friendly broadband
  leo satellite systems}, .

\bibitem[{Talent(1992)}]{talent1992}
\bibinfo{author}{Talent, D.~L.} (\bibinfo{year}{1992}).
\newblock \bibinfo{title}{Analytic model for orbital debris environmental
  management}.
\newblock {\it \bibinfo{journal}{Journal of Spacecraft and Rockets}\/},  {\it
  \bibinfo{volume}{29}\/}\bibinfo{issue}{(4)}, \bibinfo{pages}{508--513}.
  \URLprefix \url{https://doi.org/10.2514/3.25493}.
  \DOIprefix\doi{10.2514/3.25493}.
  \href{http://arxiv.org/abs/https://doi.org/10.2514/3.25493}{\tt
  arXiv:https://doi.org/10.2514/3.25493}.

\bibitem[{Verhulst(1838)}]{verhulst1838}
\bibinfo{author}{Verhulst, P.-F.} (\bibinfo{year}{1838}).
\newblock \bibinfo{title}{Notice sur la loi que la population suit dans son
  accroissement}.
\newblock {\it \bibinfo{journal}{Correspondance mathématique et physique}\/},
  {\it \bibinfo{volume}{10}\/}, \bibinfo{pages}{113--121}.

\bibitem[{Virgili(2016)}]{virgili2016delta}
\bibinfo{author}{Virgili, B.~B.} (\bibinfo{year}{2016}).
\newblock \bibinfo{title}{Delta debris environment long-term analysis}.
\newblock In {\it \bibinfo{booktitle}{Proceedings of the 6th International
  Conference on Astrodynamics Tools and Techniques (ICATT)}\/}.

\bibitem[{Walker et~al.(1999)Walker, Stokes, Wilkinson \&
  Swinerd}]{walker1999enhancement}
\bibinfo{author}{Walker, R.}, \bibinfo{author}{Stokes, P.},
  \bibinfo{author}{Wilkinson, J.} et~al. (\bibinfo{year}{1999}).
\newblock \bibinfo{title}{Enhancement and validation of the ides orbital debris
  environment model}.
\newblock {\it \bibinfo{journal}{Space Debris}\/},  {\it
  \bibinfo{volume}{1}\/}, \bibinfo{pages}{1--19}.

\bibitem[{Wang \& Liu(2019)}]{wang2019introduction}
\bibinfo{author}{Wang, X.-w.},  \& \bibinfo{author}{Liu, J.}
  (\bibinfo{year}{2019}).
\newblock \bibinfo{title}{An introduction to a new space debris evolution
  model: Solem}.
\newblock {\it \bibinfo{journal}{Advances in Astronomy}\/},  {\it
  \bibinfo{volume}{2019}\/}, \bibinfo{pages}{1--11}.

\bibitem[{Wang et~al.(2023)Wang, Wilson \& Vasile}]{wang2023multi}
\bibinfo{author}{Wang, Y.}, \bibinfo{author}{Wilson, C.},  \&
  \bibinfo{author}{Vasile, M.} (\bibinfo{year}{2023}).
\newblock \bibinfo{title}{Multi-layer temporal network model of the space
  environment}.
\newblock In {\it \bibinfo{booktitle}{2023 AAS/AIAA Astrodynamics Specialist
  Conference}\/}.

\bibitem[{Wilson et~al.(2024)Wilson, Vasile, Feng, McNally, Ant{\'o}n \&
  Letizia}]{wilson2024modelling}
\bibinfo{author}{Wilson, C.~J.}, \bibinfo{author}{Vasile, M.},
  \bibinfo{author}{Feng, J.} et~al. (\bibinfo{year}{2024}).
\newblock \bibinfo{title}{Modelling future launch traffic and its effect on the
  leo operational environment}.
\newblock In {\it \bibinfo{booktitle}{AIAA SCITECH 2024 Forum}\/} (p.
  \bibinfo{pages}{1814}).

\bibitem[{Zhang et~al.(2023)Zhang, Yuan, Yang \& Li}]{10121609}
\bibinfo{author}{Zhang, J.}, \bibinfo{author}{Yuan, Y.}, \bibinfo{author}{Yang,
  K.} et~al. (\bibinfo{year}{2023}).
\newblock \bibinfo{title}{Long-term evolution of the space environment
  considering constellation launches and debris disposal}.
\newblock {\it \bibinfo{journal}{IEEE Transactions on Aerospace and Electronic
  Systems}\/},  (pp. \bibinfo{pages}{1--17}).
  \DOIprefix\doi{10.1109/TAES.2023.3274097}.

\end{thebibliography}

\end{document}